%% file: main.tex
\documentclass[10pt]{article}
\usepackage{fullpage}
\usepackage{graphicx}
\usepackage{subcaption}
\usepackage{amsthm}
\usepackage{amsmath}
\usepackage{amsfonts}
\usepackage{amssymb}
\usepackage{tikz-cd}
\usepackage{tikz}
\usepackage{enumerate}
\theoremstyle{definition}
\newtheorem{thm}{Theorem}[section]
\newtheorem{lem}[thm]{Lemma}
\newtheorem{prop}[thm]{Proposition}
\newtheorem{cor}[thm]{Corollary}

\newtheorem{defn}[thm]{Definition}

\let\vec\mathbf

\newcommand{\y}{\mathbf{y}}
\DeclareMathOperator{\Gr}{Gr}

\usepackage[utf8]{inputenc}
\usepackage[english]{babel}
\usepackage{amssymb}
\usepackage{amsmath}
\usepackage{hyperref}
\usepackage{amsfonts}
\usepackage{mathrsfs}
\usepackage{dsfont}
\usepackage{commath}
\usepackage{graphicx}
\usepackage{text comp}

\theoremstyle{definition}
\newtheorem*{thm*}{Theorem}

\newtheorem*{lemma*}{Lemma}

\newtheorem*{corollary*}{Corollary}
\newtheorem*{prop*}{Proposition}

\newcommand{\R}{\mathbb{R}}

\newcommand{\1}{\mathbf{1}}

\usepackage[colorinlistoftodos]{todonotes}

\begin{document}

\title{Random \v{C}ech Complexes on Manifolds with Boundary}

\author{
  Henry-Louis de Kergorlay\\
  \texttt{s1562933@sms.ed.ac.uk} \\
  \and
  Ulrike Tillmann\\
  \texttt{tillmann@maths.ox.ac.uk}\\
  \and
  Oliver Vipond\\
  \texttt{vipond@maths.ox.ac.uk}\\
}







\maketitle

\begin{abstract}
Let $M$ be a compact, unit volume, Riemannian manifold with boundary.  In this paper we study the homology of a random \v Cech-complex generated by a homogeneous Poisson process in $M$. Our main results are two asymptotic threshold formulas, an upper threshold above which the \v{C}ech complex recovers the $k$-th homology of $M$ with high probability, and a lower threshold below which it almost certainly does not. These thresholds are close together in the sense that they have the same leading term. Here $k$ is positive and strictly less than the dimension $d$ of the manifold.

This extends work of Bobrowski and Weinberger in \cite{bobrowski_vanishing_2015} and Bobrowski and Oliveira \cite{Bobrowski2017} who establish similar formulas when $M$ is a torus and, more generally, is closed and has no boundary.  We note that the cases with and without boundary lead to different answers: The corresponding common leading terms for the upper and lower thresholds differ being $\log (n) $ when $M$ is closed and $(2-2/d)\log (n)$ when $M$ has boundary; here $n$ is the expected number of sample points. Our analysis identifies a special type of homological cycle, which we call a $\Theta$-like-cycle, which occur close to the boundary and establish that the first order term of the lower threshold is $(2-2/d)\log (n)$.

\end{abstract}


\section{Introduction}
\input{introduction}

\section*{Acknowledgements}
The authors would like to give recognition to The Alan Turing Institute through which the authors met, and this collaboration was initiated. They like to thank Omer Bobrowski for helpful comments on an earlier version of this paper.
HdK acknowledges support from the EPSRC studentship EP/L016508/1, and gratefully thanks The Alan Turing Institute for hosting him as an Enrichment student.
UT acknowledges support from The Alan Turing Institute through EPSRC grant EP/N510129/1.
OV gratefully acknowledges support from EPSRC studentship EP/N509711/1. UT and OV are members of the Oxford based Centre for Topological Data Analysis funded by  EPSRC grant EP/R018472/1.

\section{Background}
\input{background}

\section{Morse Theory}
\input{morse-theory}

\section{Asymptotic Coverage}
\input{coverage}

\section{Blaschke-Petkantschin Formulae}\label{Blaschke}
\input{blaschke-petkantschin}

\section{Upper Threshold}
\input{upper-threshold}

\section{Lower Threshold}
\input{lower-threshold}

\subsection{Lower Threshold Refined}
\input{lower-threshold-refined}

\section{Second Moment Calculations}
\input{second-moment-calculations}

\section{Conclusion}
\input{conclusion}

\clearpage
\section{List of Symbols}
\input{list-of-symbols}

\clearpage
\bibliography{main}
\bibliographystyle{alpha}

\end{document}

%% file: introduction.tex
Randomly generated simplicial complexes and their topology have recently attracted a lot of attention. 
The survey article by Bobrowski and Kahle \cite{bobrowski_topology_2014} collects together results 
and provides a wealth of open problems in this field. 
Here we focus on understanding the homology of random geometric complexes. While this topic was  first studied in \cite{Robins06} (see also \cite{Linial*2006}) our results build   directly on \cite{Bobrowski2017} and \cite{bobrowski_vanishing_2015}. 

Much of the current interest in the topology of random simplicial complexes is due to applications to topological data analysis, where random complexes can serve as null models when interpreting the topology of complexes on data sets.
In the context of manifold learning, a non-linear dimension reduction technique, one is interested in recovering the structure of low dimensional manifolds embedded in high dimensional space. Studying topological properties of the underlying manifold, such as homology, can inform the choices of hyperparameters in this dimension reduction technique \cite{Paul2017ASO, Tenenbaum2319}.
In the specific context of persistent homology, one of the main tools of topological data analysis, one adopts a multiscale approach to the study of the homology of randomly sampled data from a manifold. Understanding the conditions for which the topology of a complex built from a random sample coincides with the homology of the underlying manifold informs which bars in the multiscale barcode invariant relate to inherent features of the underlying manifold. For surveys of persistent homology see \cite{Zomorodian:2004:CPH:997817.997870,EdelsbrunnerHarer,Otter2017} and for work in stochastic persistent homology \cite{bobrowskiskraba2017}.

Various flavours of random simplicial complexes are present in the literature, which can include or exclude geometric considerations. Random simplicial complexes that extend the notion of a random graph in the sense of Erd\"os-R\'enyi to higher dimensional complexes are not constrained by an underlying geometry \cite{Kahle14}. In contrast, we shall work with the random simplicial complex realised by the \v{C}ech complex associated to a Poisson point process on a Riemannian manifold and ask the question when the topology of the simplicial complex approximates that of the manifold.

The question of recovering the topology of a space from a finite sample has been studied in \cite{niyogi_finding_2008,Chazal2009}  and \cite{Bobrowski2017,bobrowski_vanishing_2015} in differing contexts.  In \cite{niyogi_finding_2008} the authors consider submanifolds in Euclidean space and use the metric of the ambient Euclidean space when building the \v{C}ech complex. They provide explicit conditions to recover the homology of the manifold with high confidence. Naturally these explicit conditions are dependent on the curvature and nearness to self-intersection of the embedded manifold. In contrast we are in the context of \cite{Bobrowski2017,bobrowski_vanishing_2015}, we base the construction of the \v{C}ech complex on an intrinsic metric of the manifold independent of any embedding. Our work studies asymptotic properties and the phase transition for which one can recover the homology with high probability when increasing the size of point sample and decreasing the radius over which the associated complex is constructed. Specifically, we will give an answer to the problem posed in the survey article \cite{bobrowski_topology_2014} about extending the homological connectivity theory established for closed Riemannian manifolds in \cite{Bobrowski2017}, to Riemannian manifolds with boundary.

A principal advantage of studying asymptotics is that results rely on fewer assumptions on the underlying manifold from which the point process is sampled. Our main result like that of \cite{Bobrowski2017} thus only has dependence on the dimension of the underlying manifold and the homological dimension one wishes to recover. Our argument follows a similar framework to the argument presented in \cite{Bobrowski2017} and \cite{bobrowski_vanishing_2015}. We will need to develop completely new arguments to take into account the effect of the boundary. 

Our main result is stated with respect to the term $\Lambda := n \omega_d r^d $, the expected number of points of a uniform Poisson process of intensity $n$ lying in a $d$-dimensional radius $r$ ball. Depending on the asymptotic behaviour of the $\Lambda$, the associated \v{C}ech complex at scale $r$ built on the point process exhibits different behaviours. 
There are three distinct regimes of behaviour for $\Lambda$ as $n\to \infty, r\to 0$. In the subcritical regime ($\Lambda \to 0$) the connectivity of the \v{C}ech complex is very sparse and mostly disconnected, with the number of connected components growing at the same rate as the number of points. In the critical regime ($\Lambda \to \lambda\in (0,\infty)$) the \v{C}ech complex is sufficiently connected to exhibit non-trivial homology. However the number of connected components still grows linearly with the number of points. In the supercritical regime ($\Lambda \to \infty$) for sufficiently large $\Lambda$ the \v{C}ech complex is connected, and for even larger $\Lambda$ the point cloud covers the underlying manifold with high probability. 

Analysis in the supercritical regime yields a sequence of increasing thresholds, (\textit{homological connectivity thresholds}), such that if $\Lambda$ is greater than the $k^\text{th}$ threshold the \v{C}ech complex recovers the $k^\text{th}$ homology of the underlying closed manifold with high probability. The intermediate homological connectivity thresholds interpolate between the thresholds for more commonly studied properties, from the $0^\text{th}$ homology which detects connectivity up to the $d^\text{th}$ homology which detects coverage. We produce homological connectivity thresholds in the supercritical regime for which the \v{C}ech complex recovers the homology of a smooth compact manifold with non-trivial boundary.

The non-trivial boundary has a significant impact on the homological connectivity thresholds. As far as we know, our study of manifolds with boundary has uncovered a previously unobserved phenomenon occuring close to the boundary. Our analysis shows that close to the boundary a large number of spurious $k$-cycles appear which are not homological cycles inherent to the $k^\text{th}$ homology of the underlying manifold. This phenomenon determines that the homological connectivity thresholds for manifolds with boundary are larger than those for a closed manifold. We attain the following result:

\begin{thm*}(Homological Connectivity of Riemannian Manifold \textit{with} Boundary) 

Let $M$ be a unit volume compact Riemannian manifold with smooth non-empty boundary. Let $d\geq 2$ be the dimension of $M$, $\Lambda = \omega_d n r^d$ and $\mathcal{P}_n$ a Poisson process of intensity $n$ on $M$. Suppose $w(n)$ is any function with $w(n) \to \infty$ as $n \to \infty$. Then for $1 \leq k \leq d-1$ 

$$ \lim_{n\to \infty} \mathbb{P}(H_k(\mathcal{C}(n,r)) \cong H_k(M)) = \begin{cases}
1 & \Lambda = (2-\frac{2}{d})\log n + 2k \log \log n + w (n),\\
0 & \Lambda = (2-\frac{2}{d})\log n + 2(k - 2 -(k+1- \frac{1}{d})) \log \log n - w (n), 
\end{cases}$$
\end{thm*}

Note that when simplified the coefficient of the second order term in the lower threshold, $2(\frac{1}{d}-3)$, is independent of the homological dimension $k$. We state the coefficient of the second order term of the lower threshold of our theorem in unsimplified form to make an easy comparison to the lower thresholds established by Bobrowski and Weinberger in \cite{bobrowski_vanishing_2015} and by Bobrowski  and Oliveira in \cite{Bobrowski2017}, who studied the case when $M$ is a $d$-dimensional torus and
when $M$ is an arbitrary compact closed Riemannian manifold respectively. With the same setup for our Theorem on a closed Riemannian manifold the corresponding thresholds for $1\leq k \leq d-1$ are computed to be:

$$ \lim_{n\to \infty} \mathbb{P}(H_k(\mathcal{C}(n,r)) \cong H_k(M)) = \begin{cases}
1 & \Lambda = \log n + k \log \log n + w (n),\\
0 & \Lambda = \log n + (k-2) \log \log n - w (n),
\end{cases}$$

We see that the leading term for both the upper and lower bounds are nearly twice as large as in the case of manifolds with boundaries. At the end of the paper we provide an intuitive explanation as to why the presence of a boundary results in differences in the homological connectivity thresholds.

\subsection{Outline}

Our argument follows the same structure as the arguments presented in \cite{Bobrowski2017, bobrowski_vanishing_2015}. There are several key ideas in this framework. The first essential idea is to bound the number of homological cycles of a complex by counting the critical points of an associated Morse function. This simplifies the task of counting global phenomena of homological cycles to the purely local considerations which determine critical points of a Morse function.

In order to compute bounds for the number of critical points we require a change of variables integral formula, the Blaschke-Petkantschin formula. This change of variable formula facilitates computing the expected number of critical points induced by the distance function of a Poisson point process.

Within our argument we adapt results which apply to closed manifolds to manifolds with boundary. We use the double manifold as a canonical closed manifold in which our manifold with boundary is embedded. This trick allows us to translate results for closed manifolds to manifolds with boundary.

In Section \ref{Background} we collate results from the theory of point processes, define our asymptotic notation and also collect Riemannian volume estimates which we use in order to produce bounds in the later sections.

Section \ref{Morse Theory} provides a brief introduction to classical Morse theory. Since the distance function induced by a point process is not necessarily smooth we also provide the necessary results from \cite{Gershkovich1998} which describe Morse theory for a wider class of functions to which the distance function belongs, the so-called min-type functions.

In Section \ref{Asymptotic Coverage} we derive a coverage result for Riemannian manifolds with boundary. This coverage result is required to establish the upper homological connectivity thresholds.

We introduce in Section \ref{Blaschke-Petkantschin Formulae} the change of variables formulae we require in later sections to compute bounds on the number of critical points. These formulae are also used in Section \ref{Second Moment Calculations} in order to bound the variance of the number of critical points.

Section \ref{Upper Threshold} is dedicated to computing an upper bound for the expected number of critical points using the tools and results provided in previous sections. This upper bound on the number of critical points is used to produce an upper threshold for the homological connectivity in terms of $\Lambda$.

In Section \ref{Lower Threshold} we identify a special class of critical point which induce erroneous homological cycles, (homological cycles which are not inherent to the underlying manifold). We produce a lower bound for the expected number of this type of critical points. This lower bound is used to produce a lower threshold for the homological connectivity in terms of $\Lambda$.

Section \ref{Second Moment Calculations} bounds the variance of the number of occurences of the special type of homological cycles established in Section \ref{Lower Threshold}. We show that the variance is small and so with high probability the number of erroneous cycles behaves like the expected number of erroneous cycles.

Finally in Section \ref{Conclusion} we calculate homological connectivity thresholds given the computed bounds on the number of critical points. We compare and contrast the homological connectivity thresholds for Riemannian manifolds with boundary we have calculated to those in \cite{Bobrowski2017} for closed manifolds. We indicate how the geometric difference in the two situations inform the difference in the thresholds. 

%% file: background.tex
\label{Background}
In this section we shall provide a summary of the background theory on which we build our proofs. Results in this section are well documented so will mostly be stated without proof although references are provided for those who seek further details. 

\subsection{Homology and Complexes}

Homology is a measure of complexity of a topological space. It is an algebraic invariant that has proved a powerful tool in the study of geometry.
A thorough introduction can be found in \cite{Hatcher} far more comprehensive than the brief overview of basic homology theory we present for the unfamiliar reader.

At it's most basic level, the homology of a space can be thought of as a sequence of abelian groups, where the $k^{\text{th}}$ group summarises topological information about the $k$-dimensional subspaces. 
This algebraic summary is an incomplete invariant for homotopic spaces, that is to say,  homotopic spaces have isomorphic homology groups, although it is possible for non-homotopic spaces to have isomorphic homology groups.

For the purposes of this paper we shall consider homology over coefficients in a field, in which case the algebraic summary is a sequence of vector spaces. Equally in this setting one can consider our algebraic summary as a sequence of integers corresponding to the dimensions of these vector spaces, known as the Betti numbers of the space.

In order to make a topological space amenable to computation, we introduce a purely combinatorial object known as a simplicial complex. A simplicial complex prescribes the construction of a space out of simplices in which one glues together vertices, edges, triangles, tetrahedra and their higher dimensional analogues. The following definition describes a simplicial complex constructed from a point process on a metric space.

\begin{defn}(\v Cech Complex)
Let $\mathcal{P}$ be a collection of points in a metric space $(M,\rho)$. Define a one parameter family $\mathcal {C} (\mathcal P ,r)$ of simplicial complexes on vertex set $\mathcal{P}$ associated to this collection of points as follows: For $r \in [0 ,\infty)$, 
$$ \sigma = [p_{0},...,p_{k}] \in \mathcal{C}(\mathcal{P},r) \iff  \cap_{j=0}^k B_r(p_{j}) \neq \emptyset.$$
\end{defn}	

Here $B_r( p)$ denotes the open ball in $M$  of radius $r$ and  with centre $p$.

\begin{lem}(Nerve Lemma)
Let $\mathcal{U} = \{B_r(p) : p\in \mathcal{P}\}$ be an open cover of the metric space $(M,\rho)$. Suppose that the finite intersections of sets in $\mathcal{U}$ are empty or contractible, then the \v Cech complex $\mathcal{C}(\mathcal{P},r)$ is homotopy equivalent to $M$.
\end{lem}

The Nerve Lemma gives us a guarantee that for a suitably dense point process on a metric space and a well chosen radius $r$ we will be able to recover the homology of the underlying space.

\subsection{Poisson Point Processes}

We shall introduce in this section the notion of a general Poisson point process. Let us remark here that although later discussions will use uniformly distributed processes (with respect to the volume measure), this is merely a point of convenience. Since we are considering compact manifolds, more general distributions will only effect our results up to some constant factor.

\begin{defn}(General Poisson Point Process)\cite{Baddeley2007}
Let $(M,\mathcal{F},\mu)$ be a measure space with $M$ a compact metric space, 
 $\mathcal {F} $  a collection of measurable sets, and measure $\mu$ which is finite on compact sets and with no atoms. The Poisson process on $M$ of intensity measure $\mu$ is a point process on $M$ such that:
\begin{enumerate}
\item For every compact set $K\subset M$ the number of points $N(K)$ lying in $K$ follows a Poisson distribution with mean $\mu(K)$;
\item If $K_i \subset M$ are disjoint and compact then $N(K_i)$ are independent.
\end{enumerate} 
\end{defn} 


A salient feature of Poisson point processes is the following independence result, Theorem \ref{thm:PalmTheory1}. We shall use this result in our calculations to compute bounds on the number of critical points of the distance function from the point process which correspond to simplices of the associated \v{C}ech complex.

\begin{thm}(Palm Theory)\cite{Penrose}
Let $(X,\rho)$ be a metric space, $f: X \to \mathbb{R}$ a probability density and $\mathcal{P}_n$ a Poisson process on $X$ with intensity $\lambda_n = nf$. If $h(\mathcal{Y},\mathcal{X})$ is a measurable function for all finite subsets $\mathcal{Y} \subset \mathcal{X} \subset X^{d+1}$ with $|\mathcal{Y}| = k+1$ then:
$$ \mathbb{E}\left[\sum_{|\mathcal{Y}| = k+1} h(\mathcal{Y},\mathcal{P}_n)\right] = \frac{n^{k+1}}{(k+1)!}\mathbb{E}[h(\mathcal{Y}' ,\mathcal{Y}' \cup\mathcal{P}_n)]$$
where $\mathcal{Y}'$ is a set of $k+1$ i.i.d points in $X$ with density $f$.
\label{thm:PalmTheory1}
\end{thm}

\begin{cor}
Let $(X,\rho)$ be a metric space, $f: X \to \mathbb{R}$ a probability density and $\mathcal{P}_n$ a Poisson process on $X$ with intensity $\lambda_n = nf$. If $h(\mathcal{Y},\mathcal{X})$ is a measurable function for all finite subsets $\mathcal{Y} \subset \mathcal{X} \subset X^{d+1}$ with $|\mathcal{Y}| = k+1$ then:
\label{cor: PalmTheory2}

$$ \mathbb{E}\left[\sum_{\substack{|\mathcal{Y}_1| = |\mathcal{Y}_2| = k+1, \\ |\mathcal{Y}_1 \cap \mathcal{Y}_2| = j}} h(\mathcal{Y}_1,\mathcal{P}_n)h(\mathcal{Y}_2,\mathcal{P}_n)\right] = \frac{n^{2k-j}}{j!(k-j)!}\mathbb{E}[h(\mathcal{Y}_1' ,\mathcal{Y}' \cup\mathcal{P}_n)h(\mathcal{Y}_2' ,\mathcal{Y}' \cup\mathcal{P}_n)]$$
where $\mathcal{Y}' = \mathcal{Y}'_1 \cup \mathcal{Y}'_2$ is a set of $2k-j$ i.i.d points in $X$ with density $f$, with $|\mathcal{Y}'_1 \cap \mathcal{Y}'_2| = j$ and $\mathcal{Y}'$ independent of $\mathcal{P}_n$.

\end{cor}

\subsection{Asymptotic Notation}

We shall use the following set of notation to denote different asymptotic behaviours of functions $f,g$:

\begin{enumerate}
    \item $f(n) = O(g(n))$ if there is a constant $C$ and $n_0 \in \mathbb{N}$ such that $|f(n)| \leq C |g(n)|$ for all $n> n_0$
    \item $f(n) = o(g(n))$ if $\lim_{n \to \infty} \frac{|f(n)|}{|g(n)|} = 0 $
    \item $f(n) = \Omega(g(n))$ if if there is a constant $C$ and $n_0 \in \mathbb{N}$ such that $|f(n)| \geq C |g(n)|$ for all $n> n_0$
\end{enumerate}

Moreover we shall use the $f(n) = O(g_1(n),g_2(n)) := O(g_1(n) + g_2(n))$, to emphasise when our bounding functions are coming from distinct calculations. 

\subsection{Riemannian Volumes}
\label{subsec:RiemannianVolumes}
A major portion of our later proofs require us to bound various Riemannian volumes by their Euclidean counterparts in order to control the asymptotic behaviour of the Betti numbers. In this section we shall provide these approximations of Riemannian volumes. Let us denote a smooth Riemannian manifold as the pair $(M,g)$ and consider the case when $M$ is compact, of dimension $d$, and the metric $g$ is smooth. An introduction to Riemannian Geometry can be found for example in \cite{Lee2012}.

A smooth metric is a smoothly varying inner product on the tangent space $g : T_pM\times T_pM \to \mathbb{R}$ and therefore endows the tangent space at each point with a norm. We define the length of a path using $g$ by integrating over the norm of its velocity. We shall use $\rho(p_1,p_2)$ to denote the shortest path length between two points $p_1,p_2$ on our manifold. 

Let us denote the open ball of radius $r$ about a point $p$ on our manifold as $B_r(p)$, and the sphere of radius $r$ about $p$ as $S_r(p)$. If we are considering a collection of points $\mathcal{P}$ then we denote the union of open balls of radius $r$ centred at each point as $B_r(\mathcal{P}) := \bigcup_{p\in \mathcal{P}} B_r(p)$. 

The exponential map $\text{exp}_p : T_pM \to M$ is defined by $\text{exp}_p(\vec{v}) = \gamma(1)$ where $\gamma$ is the unique geodesic in $M$ with $\dot{\gamma}(0)=\vec{v}$. Since $\text{exp}_p$ is a local diffeomorphism, an orthonormal basis of $T_pM$ induces local coordinates about $p$, which we shall denote as $(x^1,x^2,...,x^d)$ and are called geodesic normal coordinates. Using the Taylor expansion we can write the metric $g$ in terms of these coordinates where $R_{iklj}$ is known as the Riemann curvature tensor:

$$ g = \sum_{i,j} g_{ij} dx^i \otimes dx^j \text{ with } g_{ij} = \delta_{ij} + \frac{1}{3}\sum_{k,l}R_{iklj}x^kx^l + O(|x|^3).$$

The Euclidean metric on $M= \mathbb R^d$ is the simple case where $g_{ij} = \delta_{ij}$. Given a point $p$ of a Riemannian manifold $(M,g)$ with local neighbourhood $U$ and geodesic normal coordinates $(x^1,...,x^d)$ in the neighbourhood $U$, for sufficiently small radius $r$ we can consider the intrinsic Euclidean ball $B_r^E(p)$ as the radius $r$ ball with respect to the metric $g_E$ on $U$ where $g_E = \sum_{i,j} \delta_{ij}dx^i \otimes dx^j$. Let us denote the $d-1$ dimensional unit round sphere as $\mathbb{S}^{d-1}$.
%
%
The canonical measure induced by the Riemannian density on the manifold $M$ can be expressed in terms of the Euclidean measure $|\text{dvol}_{g_E}|$ associated to the Euclidean metric:
$$ |\text{dvol}_g| = \sqrt{|\det (g_{ij})|} |\text{dvol}_{g_E}|.$$
The Ricci curvature tensor at a point $p$ is given by $Ric_{ij} = - \sum_k R_{ikkj}$, and we can calculate that:
$$ \sqrt{|\det (g_{ij})|} = 1 - \frac{Ric_{ij}}{3}x^ix^j +O(|x|^3).  $$

Our first volume approximations also apply to balls and spheres within a manifold with boundary which are wholly contained within the manifold, that is the centres lie far enough from the boundary.

\begin{lem}\cite{Bobrowski2017}
Let $(M,g)$ be a closed compact Riemannian manifold of dimension $d$, and let $\omega_d$ denote the volume of the $d$-dimensional unit Euclidean ball. Let $\text{Vol}$ denote the Riemannian volume  on $M$. Then:
$$ \text{Vol}(B_r(p)) = \omega_d \, r^d\left( 1 - \frac{s(p)}{6(d+2)}r^2 + O(r^3)\right),$$
$$ \text{Vol}(S_r(p)) = d\, \omega_d \, r^{d-1}\left( 1 - \frac{s(p)}{6d}r^2 + O(r^3)\right),$$
where $s(p) = \sum_i Ric_{ii}$ is the scalar curvature at $p$, and $ \text{Vol}(S_r(p))$ denotes the volume with respect to the induced metric on $S_r(p)$.
\end{lem}


%
%
\begin{lem}\cite{Bobrowski2017}
Let $(M,g)$ be a closed compact Riemannian manifold of dimension $d$. Let $|Ric_p| = \sup_{v\in T_pM\setminus 0 } \frac{|Ric(\vec{v},\vec{v})|}{|\vec{v}|^2}$ denote the norm of the Ricci tensor at $p$. For any $\nu>0$ there is a continuous choice of $r_\nu>0$ such that for any smaller radius $r\leq r_\nu$ and any $p\in M$ the following bounds hold on $B_r(p)$:
$$ r^{d-1}\left(1 - \frac{|Ric_p| + \nu}{3}r^2\right)|\text{dvol}_{\mathbb{S}^{d-1}}| \leq |\text{dvol}_{S_r(p)}| \leq  r^{d-1}\left(1 + \frac{|Ric_p| + \nu}{3}r^2\right)|\text{dvol}_{\mathbb{S}^{d-1}}| $$

\end{lem}

Using the polar decomposition of the volume of a ball we attain the following corollary:

\begin{cor}\cite{Bobrowski2017}
Let $s(p)$ denote the scalar curvature and define $s_{\text{min}}(\nu) = \inf_{p\in M} \frac{s(p)}{6(d+2)} - \nu$, $s_{\text{max}}(\nu) = \sup_{p\in M} \frac{s(p)}{6(d+2)} + \nu$. Then for all $\nu >0$ there is a continuous choice of $r_\nu >0$ such that for all $r\leq r_\nu$
$$ \omega_dr^d(1-s_{\text{max}}(\nu)r^2) \leq \text{Vol}(B_r(p)) \leq \omega_dr^d(1-s_{\text{min}}(\nu)r^2)$$
\label{cor: Ball Volume}
\end{cor}

\begin{lem}\cite{Bobrowski2017}
Let $(M,g)$ be a compact Riemannian manifold of dimension $d$. For all $\nu>0$ there is a continuous choice of $r_\nu >0$ such that for all $r\leq r_\nu$, $p\in M$ we have the following bounds hold on $B_r(p)$

$$ (1-\nu r^2)|\text{dvol}_{g_E}| \leq |\text{dvol}_{g}| \leq (1+\nu r^2)|\text{dvol}_{g_E}|$$
\end{lem}

\begin{lem}\cite{Bobrowski2017}
Let $(M,g)$ be a compact Riemannian manifold of dimension $d$. Then there is some $\nu>0$ and $r_\nu >0$ such that for all $r<r_\nu$ and any two point with $\text{dist}(p_1,p_2) <2r$ we have:
$$ \left(B^E_{(1-\nu r)r}(p_1) \cup B^E_{(1-\nu r)r}(p_2) \right) \subset (B_r(p_1) \cup B_r(p_2)) \subset \left(B^E_{(1+\nu r)r}(p_1) \cup B^E_{(1+\nu r)r}(p_2) \right)$$
\label{lem:IntersectingBalls}
\end{lem}

The proofs of these approximations can be found in the Appendix of \cite{Bobrowski2017}.

\vskip .2in

Let us now consider how we must adapt these approximations for balls intersecting the boundary. We may consider our manifold to be embedded in $\mathbb{R}^D$, for some $D$. Indeed, Nash's Imbedding Theorem \cite{nash1956} guarantees existence of an isometric embedding. 

\begin{defn}(Reach of a Manifold)
Let the medial axis $\text{Med}(M)$ of a manifold $M$ embedded in $\mathbb{R}^D$, be the set of points in $\mathbb{R}^D$ which do not have a unique nearest element in $M$. The reach $\tau_M$ of a manifold is defined to be $\tau_M = \inf_{p\in \text{Med}(M)} \text{dist}(M,p)$.
\end{defn}

\begin{figure}[ht]
\centering
\includegraphics[scale=0.4]{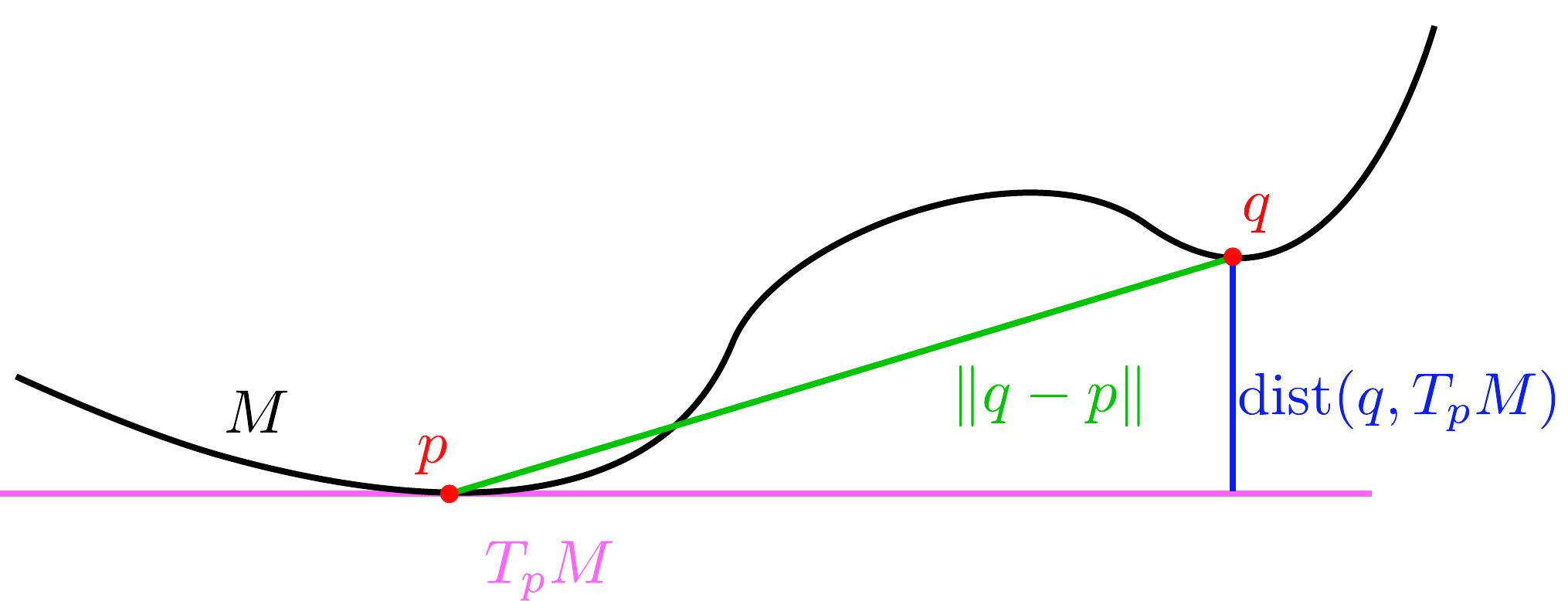}
\caption{We sketch the distances involved in Theorem \ref{thm:Federer} for which we can calculate the reach. }
\label{fig:Reach-Sketch}
\end{figure}

\begin{thm}(Theorem 4.18 \cite{federer_curvature_1959}) 
Let $M$ be a submanifold of $\mathbb{R}^D$ and $\tau_M$ the reach of $M$. Then the reach is realised as:
$$ \tau_M = \inf_{q\neq p \in M} \frac{\|q-p\|^2}{2\text{dist}(q,T_pM)}$$

See Figure \ref{fig:Reach-Sketch}.
\label{thm:Federer}
\end{thm}

The manifolds we consider are smooth and compact and so in particular the reach of our manifolds is non-zero. We shall frequently consider the double manifold $DM = M \cup_{\partial M} M'$ associated to a manifold with boundary, since this gives us a canonical compact closed manifold in which $M$ is embedded; $M'$ denotes here an identical copy of $M$. The following proposition provides an estimate for the volume of a ball centred near the boundary.

\begin{prop}
Let $B_r^{DM}(c)$ be a ball of radius $r$ centred at $c\in M$ in the double manifold $DM$, and let $\delta$ denote the distance from the centre $c$ to the boundary. Let $B_r^M(c) = B_r^{DM}(c) \cap M$ denote the portion of this ball contained in $M$. Then the volume of $B_r^M(c)$ can be expressed as:
 $$\textrm{Vol}(B_r^{M}(c)) = (\frac{1}{2} + \alpha)\textrm{Vol}(B_r^{DM}(c))$$
where $\alpha = O(\frac{\delta}{r}) \pm O(r)$.
 \label{prop:Half-BallVolumes}
\end{prop}
\begin{proof}
Let us assume that the Riemannian metric is Euclidean, since we have seen above that this will only change the volume of our ball by a factor of order $(1-O(r^2))$. Let $p$ denote the closest point on the boundary to $c$. We shall first bound the volume of the ball $B_r^{DM}(c)$ capped by the tangent plane $T_pM$,  so that we are dealing with a Euclidean ball from which we have removed a hemispherical cap. The volume of a $d$-dimensional ball of radius $r$ capped at height $h$, $B_{r}^h$, can be expressed as a fraction of the volume of the full ball $\omega_d r^d$:

\[
 \text{Vol}(B_{r}^h) = \omega_dr^d\frac{1}{2}\left(1+\frac{G_d(q)}{G_d(1)}\right)=:\omega_d r^d(\frac{1}{2} + \alpha),
\]
where $G_d(u) = \int_0^u(1-t^2)^{(d-1)/2}dt$, $q = \frac{\delta}{r}$, and $\delta = r-h$. For sufficiently small $q$ we have the trivial bounds $\frac{1}{2}q \leq G_d(q)\leq q$, and so we see that:
\[
  \alpha = O(\frac{\delta}{r}).
\]
Let us now estimate the error volume bounded between the boundary and $T_pM$. Consider the local parameterization of $\partial M$ projected on $T_pM$ using the local diffeomorphism $\textrm{exp}_p$. Denote this local parameterization by $e:B_{r_2}^{(d-1)}(0)\to \R$. The error volume induced by the non-flatness of the boundary is then:
\[
 \int_{B_{r_2}^{(d-1)}(0)}e(x)dx.
\]
Using Theorem \ref{thm:Federer} we observe that $e(x) = O(r^2)$ and so the error term is of order $O(r^{d+1})$.
Thus we attain the desired result that:
$$\textrm{Vol}(B_r^{M}(c)) = (\frac{1}{2} + \alpha)\textrm{Vol}(B_r^{DM}(c))$$
Where $\alpha = O(\frac{\delta}{r}) \pm O(r)$.

\end{proof}

We can combine this result with Lemma \ref{lem:IntersectingBalls} to similarly  estimate the volume of intersecting balls lying close to the boundary.

If $M$ is a smooth manifold with non-empty boundary then we call a neighbourhood of the boundary $\partial M$ a \textit{collar neighbourhood} if it is the image of a smooth embedding $\partial M \times [0,1) \to M$ and the embedding restricts to the identification $\partial M \times \{0\} \to \partial M$. A standard result (sometimes called the Collar Neighbourhood Theorem \cite{Lee2012}) guarantees any smooth manifold with non-empty boundary has a collar neighbourhood. For small $r$ the $r$-neighbourhood of the boundary denoted $\partial M_r$, is a collar neighbourhood. It is straight forward to show that the volume of this collar neighbourhood $\text{Vol}(\partial M_r) \sim r \text{Vol}(\partial M)$.

%% file: morse-theory.tex
\label{Morse Theory}

\subsection{Morse Theory on Manifolds with Boundary}

We start with a terse summary of Morse theory for manifolds with boundary. 
The upshot of this section is that 
given a suitable smooth function on the manifold  we define the Morse complex in the same manner as in   \cite{Milnor1965}. This complex will compute relative homology $H_*(M, \partial M)$ if the function attains a minimum on the boundary, and will compute absolute homology $H_*(M)$ if the function attains a maximum on the boundary. 

\begin{defn}(Non-Degenerate Critical Point)\cite{Milnor1965}
Let $M$ be a smooth $d$-manifold with boundary and $f: M \to \mathbb{R}$ a smooth function. A point $p\in M$ is critical if $\nabla f(p)=0$. A critical point is non-degenerate if the Hessian $H(p)$ is non-singular.

\end{defn}

\begin{defn}(Smooth Manifold Triad)\cite{Milnor1965}
We say $(M;V_0,V_1)$ is a smooth manifold triad if $M$ is a smooth manifold with boundary, and the boundary $\partial M$ is the disjoint union of the open and closed submanifolds $V_0,V_1$.
\end{defn}

\begin{defn}(Morse Function)\cite{Milnor1965}
A smooth function $f:(M;V_0,V_1) \to [a,b]$ is a Morse function on a smooth manifold triad if $f^{-1}(a)=V_0,f^{-1}(b)=V_1$ and all the critical points of $f$ lie in the interior of $M$, ($\text{int}(M)=M \setminus \partial M$).
\end{defn}

In this setting we still have the classical Morse Lemma \cite{Smiley1964a} for closed manifolds that asserts the existence of a coordinate system about each critical point for which the Morse function has a  diagonal quadratic form. The index of a critical point is again given by the dimension of the negative eigenspace of the Hessian matrix $H(p)$.

If we have a cobordism $c$ represented by the smooth manifold triad $(M;V_0,V_1)$ and a Morse function $f: (M;V_0,V_1)\to \mathbb{R}$ then we can factor the cobordism as $c = c_0c_1...c_d$ such that each $c_\lambda$ admits a Morse function with critical points all of index $\lambda$, and such that the critical points of $c_\lambda$ are in one to one correspondence with the index $\lambda$ critical points of $f$.

Let the composite of the first $\lambda$ cobordisms be represented by the manifold $M_\lambda$, and for the edge case set $M_{-1} = V_0$. Define $C_k = H_k(M_k,M_{k-1})$ and let $\partial : C_k \to C_{k-1}$ be the boundary homomorphism in the LES of the triple $(M_{k-2},M_{k-1},M_k)$.

\begin{thm}(Theorem 7.4)\cite{Milnor1965}
With $C_\ast$ and $\partial$ defined as above, $(C_k,\partial)$ is a chain complex and moreover $H_k(C_\ast) \cong H_k(M,V_0)$.
\label{thm: Morse Theory}
\end{thm}

In particular we will be interested in the case where we decompose the boundary trivially and so recover either $ H_k(M)$ or  $H_k(M,\partial M)$. The use of Morse theory in our arguments will not rely on any knowledge of the boundary maps, and will solely be used to bound Betti numbers.


\subsection{Morse Theory for Min Type Functions}

Whilst regular Morse theory is concerned with smooth functions, the distance function associated to a point cloud 
on a Riemannian manifold is generally not 
smooth. However, the square of such a distance function is a min-type function for which a version of Morse theory has been developed \cite{Gershkovich1998}. 

\begin{defn}(Min-type function)\cite{Gershkovich1998}
Let $f: \mathbb{R}^d \to \mathbb{R}$ be a germ of a continuous function at $p\in \mathbb{R}^d$. Then $f$ is a germ of a min-type function at $p$ if there exist germs of smooth functions $\alpha_i$ at $p$ such that locally around $p$ we have that $f= \min_{i=1}^m \alpha_i$.
A function on a $d$-manifold $f: M \to \mathbb{R}$ is min-type if for all $p\in M$ the germ at $p$ is a germ of a min-type function.
\end{defn}

A couple of technical Lemmas allow us to show that each germ of a min-type function has an essentially unique minimal representation. We use this canonical minimal representation to build the min-type version of Morse Theory.  

\begin{defn}(Non-degenerate regular point)\cite{Gershkovich1998}
Let $f: \mathbb{R}^d \to \mathbb{R}$ be a germ of a min-type function. The point $p\in \mathbb{R}^d$ is a non-degenerate regular (NDR) point if $f$ has a minimal representation at $p$, $f= \min_{i=1}^m \alpha_i$, such that:

\begin{enumerate}
\item $\{ \nabla(\alpha_i - \alpha_{i+1})(p) \}_{i=1}^{m-1}$ is linearly independent
\item $0 \notin Grad(f) := \text{Conv}\{ \nabla \alpha_i(p) \}_{i=1}^m$ where Conv denotes the convex hull
\item $f\vert_{\mathcal{G}_f}$ is a germ of a Morse function on the boundary set $\mathcal{G}_f := \{ x : \alpha_1(x) = ... = \alpha_m(x)\}$
\item Any $m-1$ gradients among $\{ \nabla \alpha_i(p) \}_{i=1}^m$ are linearly independent
\end{enumerate}

\end{defn}

Condition $1$ ensures that $\mathcal{G}_f$ is a smooth $(d-m+1)$-submanifold. If the gradients $\{ \nabla \alpha_i(p) \}_{i=1}^m$ are linearly dependent then the convex hull $\text{Conv}\{ \nabla \alpha_i(p) \}_{i=1}^m = Grad(f)$ can be thought of as an $(m-1)$-simplex in the tangent space.

\begin{defn}(Non-degenerate critical point)\cite{Gershkovich1998}
Let $f: \mathbb{R}^d \to \mathbb{R}$ be a germ of a min-type function. The point $p\in \mathbb{R}^d$ is a non-degenerate critical (NDC) point if $f$ has a minimal representation at $p$, $f= \min_{i=1}^m \alpha_i$, such that:

\begin{enumerate}
\item $\{ \nabla(\alpha_i - \alpha_{i+1})(p) \}_{i=1}^{m-1}$ is linearly independent
\item $0 \in Grad(f) = \text{Conv}\{ \nabla \alpha_i(p) \}_{i=1}^m $
\item $f\vert_{\mathcal{G}_f}$ is a germ of a Morse function on the boundary set $\mathcal{G}_f = \{ x : \alpha_1(x) = ... = \alpha_m(x)\}$
\item Any $m-1$ gradients among $\{ \nabla \alpha_i(p) \}_{i=1}^m$ are linearly independent
\end{enumerate}

\end{defn}

\begin{defn}(Morse min-type function)
We say that a min-type function $f: M \to \mathbb{R}$ is a Morse min-type function if every point is either an NDR or an NDC point.
\end{defn}

If $p$ is an NDC point of $f:M\to \mathbb{R}$ then $p$ is an NDC point of the smooth Morse function $f\vert_{\mathcal{G}_f}$ in the usual Morse theoretic sense. Consequently we define the index of an NDC point as follows:

\begin{defn}(Index)\cite{Gershkovich1998}
Let $f: \mathbb{R}^n \to \mathbb{R}$ be a Morse min-type function with NDC point $p$ and associated minimal representation $f= \min_{i=1}^m \alpha_i$. The index of $p$ is defined to be $\text{Ind}_p(f) =  (m-1) + \text{Ind}_p(f\vert_{\mathcal{G}_f})$
\end{defn}

One of the main results from \cite{Gershkovich1998} relates Morse min-type functions to smooth Morse functions. In essence the following Theorem says that a Morse min-type function can be $\varepsilon$-approximated by a classical smooth Morse function with the same number of critical points of the same index, and these points are arbitrarily close to the original critical points.

\begin{thm}(Morse Min-Type Approximation Theorem)\cite{Gershkovich1998}
Let $f$ be a Morse min-type function on a compact closed Riemannian manifold $(M,g)$ with critical points $y_1,...,y_m$ and corresponding indices $q_1,...,q_m$. For any $\varepsilon>0$ there is a smooth Morse function $f_\varepsilon$ satisfying the following properties:
\begin{enumerate}
\item $f_\varepsilon$ is an $\varepsilon$-approximation of $f$ in the $C^0$- metric

\item $f_\varepsilon$ is an $\varepsilon$-approximation of $f$ in the $C^2$- metric for any neighbourhood where $f$ is smooth

\item $f_\varepsilon$ has the same number of critical points $y_1^\varepsilon,...,y_m^\varepsilon$ with $\text{Ind}_{y_i}(f) = \text{Ind}_{y_i^\varepsilon}(f_\varepsilon)$

\item $\rho(y_i,y_i^\varepsilon) \leq \varepsilon$
\end{enumerate}
\label{Approx}
\end{thm}


\subsection{Morse Theory for the Distance Function on a Compact Manifold with Boundary}

In this section we seek to recover a Morse function on the triad $(M,\emptyset,\partial M)$ induced by the distance function of a point cloud on $M$.  For $x,y\in M$ and $ \mathcal P$ a finite subset of $M$ define
distance functions
$$
\rho^2 _x(y) := \rho^2 (x,y) \quad \text { and } \quad  \rho^2_{\mathcal P} (y):= \min_{x\in \mathcal P} \rho^2_x(y).
$$


\begin{lem}\cite{Bobrowski2017}
Let $(M,g)$ be a compact Riemannian manifold with boundary, then there exists a positive real number  $r_{\text{mt}}>0$ such that for every finite sample of points $\mathcal{P}\subset M \setminus \partial M$  the distance function $\rho_\mathcal{P}^2$ is a Morse min-type function on the $r_{\text{mt}}$ neighbourhood of these points $B_{r_{\text{mt}}}(\mathcal{P})$.
\end{lem}

\begin{proof}
Given any point $x\in M$ the distance function $\rho_x^2$ is smooth, Morse,  and strictly convex on some neighbourhood $B_{r_x}(x)$. Since our metric $g$ is smooth, $r_x$ may be chosen continuously. By compactness there is some positive $r_{\text{mt}} \leq r_x$ for all $x\in M$.

\label{Min-Type}
\end{proof}

\begin{prop}
Let $(M,g)$ be a compact Riemannian manifold with boundary and let $\mathcal P$ be a finite subset of $M$ such that $M \subset B_r(\mathcal{P})$ for some $r < \frac{r_{\text{mt}}}{2}$.
 Then  the critical points of the Morse min-type distance function $\rho_{\mathcal{P}}^2$ are in one to one index preserving correspondence with the critical points of a Morse function on the smooth manifold triad $(M;\emptyset,\partial M)$.
\label{prop: Morse Min Type Dist}
\end{prop}

\begin{proof}
Consider the double $DM = M\cup_{\partial M} M'$  and let $\mathcal P' \subset M'\setminus \partial M'_{\frac{r_{mt}} {2}} $ be a finite set of points in $M'$  such that
$$
B_{r_{mt} } (\mathcal P \cup \mathcal P' ) = DM.
$$
Then the distance function $\rho^2_{\mathcal P \cup \mathcal P'}$ is of Morse min-type on all of $DM$ and hence, by the Approximation Theorem, has an $\epsilon$-approximation 
by a smooth Morse function $f$ on $DM$.
As $\mathcal P'$ is bounded away from the common boundary $\partial M= \partial M'$ by $\frac{r_\text{mt}}{2} $ and $M \subset B_r(\mathcal{P})$ with $r < \frac{r_{\text{mt}}}{2}$, the distance function $\rho^2 _{\mathcal P \cup \mathcal P'}$, and thus also $f$, increases on $\partial M$ in the direction of the normal pointing into $M'$. Hence $f|_M$ can be extended to a smooth Morse function on $M_\delta := M \cup (\partial M \times [0, \delta])$ which attains its maximum on
the boundary $\partial M \times \{\delta \}$. Identifying $M_\delta$ with $M$ gives the required result.
\end{proof}

\begin{cor}
Let $\mathcal{P}_n$ be a Poisson process on $M$ a compact Riemannian manifold with boundary. Then with high probability the distance function $\rho^2_{\mathcal P_n}$ induces a Morse function on the smooth manifold triad $(M;\emptyset,\partial M)$.
\end{cor}
\begin{proof}
Observe that for a Poisson process $\mathcal{P}_n$ on $M$ and fixed $r < \frac{r_\text{mt}}{2}$, it follows that $M \subset B_r(\mathcal{P}_n)$ with high probability. Hence using Proposition \ref{prop: Morse Min Type Dist} the point process $\mathcal{P}_n$ induces a Morse min-type distance function $\rho^2 _{\mathcal{P}_n}$ whose critical points recover the homology of $M$ with high probability.
\end{proof}

Let us consider the conditions for which a point of $M$ is a $k$-critical point or non-degenerate regular point of the distance function induced by a finite set  $\mathcal{P} \subset M$.

For $y\in \mathcal{P}$ the local behaviour of $\rho^2_\mathcal{P}$ close to  $y$ coincides with $\rho^2_y(x) =\rho^2(x,y)$. Thus the minimal representation of $\rho^2_\mathcal{P}$ at $y$ is $\rho^2(-, y)$, the submanifold $\mathcal{G}_f$ is just $M$, and  $\nabla \rho^2_y(y) = 0$. So $y$ is an NDC point of index $0$.

For any $p \in M \setminus \mathcal{P}$ if there is a point  $y\in\mathcal{P}$ strictly closer to $p$ than any other point in $\mathcal P$, then locally about $p$ we have that $\rho^2_\mathcal{P}(x) = \rho^2_y(x)$ and  $\nabla \rho^2_y(p) \neq 0$. So $p$ is an NDR point.

Let $p \in M$ be such a point with  $ \min_{y\in \mathcal{P}}\rho(p,y)$ achieved precisely at all $p \in  \mathcal Y = \{y_0,...,y_k\}\subset \mathcal{P}$. Then locally about $p$ we have the minimal representation $\rho^2_\mathcal{P}(x) = \min_{i=0}^k\rho^2_{y_i}(x)$. 
In order for $p$ to be critical we require linear independence of the set $\{\nabla(\rho^2_{y_i} -\rho^2_{y_{i+1}})(p)\}_{i=0}^{k-1}$ which corresponds to saying that the set $\{y_0,...,y_k\}$ is generic. Further we require that at $p$ we have $0\in Grad(f)$, which corresponds to $p$ lying in the convex hull of the set $\{y_0,...,y_k\}$. If these conditions are met, $p$ is critical and the index of such a critical point will be $k$, as the point $p$ attains a minimum of the distance function restricted to the submanifold $\mathcal{G}_f$ at $p$.

In terms of the \v{C}ech complex construction built on $\mathcal{P}$, an index $k$ NDC critical point with critical value $r$ occurs at the point $x\in M$ if $x$ is the point of intersection of closed radius $r$ balls about $k+1$ points of $\mathcal{P}$, and lies in their convex hull. We can thus identify index $k$ critical points of the distance function with the introduction of $k$ simplices to the \v{C}ech complex at their respective critical values. 

Let $\mathcal{Y} = \{y_0,...,y_k\} \subset M$  and define the $\mathcal{Y}$-equidistant sets:

$$ E(\mathcal{Y}) := \{ x \in M \, | \, \rho_{y_0}(x) = ... = \rho_{y_k}(x)  \} \ \ ,\ \ E_r(\mathcal{Y}) := E(\mathcal{Y}) \cap B_r(\mathcal{Y}).$$
The following result is a generalisation to Riemannian manifolds of the fact that $k+1$ generic points in Euclidean space lie on a sphere of dimension $k-1$. In particular we can associate a centre and radius to a collection of points which are sufficiently close together.

\begin{lem}\cite{Bobrowski2017}
There is a positive $r_{\text{max}}< r_{\text{mt}}$ such that if $\mathcal{Y}$ is generic with $E_{r_{\text{max}}}(\mathcal{Y}) \neq \emptyset$, then the set $\mathcal{Y}$ has a unique point $c(\mathcal{Y})\in M$ such that for all $p \in \mathcal Y$   
$$
\rho_p(c(\mathcal{Y}))  = \inf_{x\in E(\mathcal{Y})} \rho_{\mathcal{Y}}(x).
$$

\label{lem: Unique Centres}
\end{lem}
In this case, we will refer to the point  $c(\mathcal{Y})$ as the {\it centre} of $\mathcal Y$ and
to $\rho(\mathcal Y):= \rho_{\mathcal Y} (c(\mathcal Y))$ as its {\it radius}. The set $\mathcal{Y} = \{y_0,...,y_k\}$ corresponds to an index $k$ critical point of the distance function $\rho^2 _{\mathcal P}$ with critical value $\rho^2(\mathcal{Y})$.

%% file: coverage.tex
\label{Asymptotic Coverage}
In forming the upper threshold
we need to understand the asymptotic coverage of a Riemannian manifold with boundary. We would like to understand the conditions under which $B_r(\mathcal{P})$ covers $M$ w.h.p. 

The paper Random Coverings \cite{Flatto} gives a comprehensive treatment of the asymptotic behaviour of the number of randomly chosen radius $r$ balls required to cover a compact \textit{closed} Riemannian manifold $m$ times. The result is translated into a sharp coverage threshold in \cite{bobrowski_vanishing_2015}, but again this result applies only to \textit{closed} Riemannian manifolds.
We use the coverage result of \cite{bobrowski_vanishing_2015} to derive conditions (not necessarily sharp) for asymptotic coverage of a compact Riemannian manifold with boundary.
A sharp coverage result for manifolds with boundary found in \cite{Wei18} has come to our attention. Nevertheless, we include our coverage result which is sufficiently strong for our purposes and arises from a natural geometric argument.

\begin{thm}(Sharp Coverage Threshold)\cite{bobrowski_vanishing_2015}

Let $M$ be a compact, closed, unit volume Riemannian manifold. Let $\mathcal{P}_n$ be a uniform Poisson process on $M$ of intensity $n$. Let $B_{r}(\mathcal{P}_n)$ denote the $r$ neighbourhood of the points $\mathcal{P}_n$ and $w(n)\to \infty$, then we yield:
$$ \lim_{n\to \infty} \mathbb{P}(M \subset B_{r}(\mathcal{P}_n)) = \begin{cases}
1 & \Lambda = \log n + (d-1) \log \log n + w (n),\\
0 & \Lambda = \log n + (d-1) \log \log n - w (n).
\end{cases}$$
\label{SharpCoverage}
\end{thm}





Our coverage result below is not proposed as a sharp threshold but merely a threshold that suits our needs in subsequent proofs.

\begin{cor}(Coverage of Compact Manifold with Boundary)

Let $M$ be a compact Riemannian manifold with boundary and $C >1$ a constant. Let $\mathcal{Q}_n$ be a uniform Poisson process on $M$ with intensity $n$. If $\Lambda = C\log n $, then the $2r$ neighbourhood of the Poisson process covers $M$ w.h.p:

$$ \lim_{n\to \infty} \mathbb{P}(M \subset B_{2r}(\mathcal{Q}_n)) = 1 . $$
\label{Coverage}
\end{cor}

\begin{proof}
Our proof will be an application of the sharp threshold developed in \cite{bobrowski_vanishing_2015}. Let $C = 1 + \delta$ then let us form the rescaled double manifold $DM^{\delta} = M \cup_{\partial M} M^{\delta}$ where the metric on $DM$ is smoothly rescaled so that $\text{Vol}(M^{\delta})=\delta$ and the metric on $M$ is unchanged.
Now consider a Poisson process $\mathcal{P}_{Cn}$ of intensity $Cn$ on $DM^{\delta}$. This restricts to a Poisson process $\mathcal{Q}_{n}=\mathcal{P}_{Cn} \cap M $ of intensity $n$ on $M$. Let $\Lambda^{\delta} = (1+\delta)^{-1}n\omega_dr^d = (1+\delta)^{-1}\Lambda$. Then applying Theorem \ref{SharpCoverage} we see that if $\Lambda^{\delta} = \log n + (d-1) \log \log n + \omega (n)$
$$ \lim_{n\to \infty} \mathbb{P}(DM^{\delta} \subset B_{r}(\mathcal{P}_{Cn})) = 1  $$
Let us define the following notation $M_r = M \setminus (\partial M \times [0,r) )\subset DM^{\delta}$.
For any point $x \in M_r$ the $r$-neighbourhood of $x\in DM^{\delta}$ coincides with the $r$-neighbourhood of $x\in M$. For points in the collar $\partial M \times [0,r)$ we cannot make the same statement since the attachment of $M^{\delta}$ may have introduced $r$-geodesics between points in the collar that were previously separated by a larger distance.
The theorem follows immediately from the observation that: $$DM^{\delta} \subset B_{r}(\mathcal{P}_{Cn}) \implies M_r \subset B_{r}(\mathcal{P}_{Cn}\cap M)\subset DM^{\delta} \implies M_r \subset B_{r}(\mathcal{Q}_{n})\subset M \implies M \subset B_{2r}(\mathcal{Q}_{n})\subset M.$$

\end{proof}

%% file: blaschke-petkantschin.tex
\label{Blaschke-Petkantschin Formulae}
A key component of our later arguments and the arguments found in \cite{Bobrowski2017} is an integral formula which facilitates calculating bounds on the expected number of critical points of a distance function associated to a Poisson process on a Riemannian manifold. In this section we explain this change of variables formula and make appropriate adaptations to the case for Riemannian manifolds with boundary.
 
\subsection{The Blaschke-Petkantschin formula in the Euclidean case}
We first recall a derivation of the classical Blaschke-Petkantschin formula in the Euclidean case. Our derivation and Proposition \ref{prop:Blaschke-Pentkanschin Formula Euclidean Case}, roughly follow Sections $2$ and $3$ of \cite{Miles1971}.\\ 

Let $E_d$ be a $d$-dimensional Euclidean space, and let $(e_i)_{i=1}^d$ be an orthonormal moving frame in $E_d$, where for an infinitesimal rotation of the frame,
\[
 e_i\cdot de_i = 0\text{, }\forall i\in [d];
\]
and set
\[
 \omega_{ij}:= e_i\cdot d e_j=  -\omega_{ji}\text{, }\forall i,j\in [d].
\]
Given $r$ points $\{x_i:i\in [r]\}\subset E_d$, Miles derives heuristically the associated volume form
\[
 \bigwedge_{j=1}^r d V(x_j) = \bigwedge_{i=1}^d\bigwedge_{j=1}^r e_i.d x_j.
\]
Furthermore, given the Grassmannian manifold $\text{Gr}(r,d)$ with invariant measure $d\mu_{r,d}(V)$, Miles also derives
\[
 d\mu_{r,d} = \bigwedge_{i=1}^r\bigwedge_{j=r+1}^d \omega_{ij}.
\]
Using the above, the Blaschke-Petkantschin formula expresses the Euclidean volume form $dV(x_i^d)$ on $\{x_i:i\in [r]\}$ in terms of the volume element associated to the $r$-plane containing $\{x_i:i\in [r]\}$, denoted by $dV(x_i^r)$.

 \begin{prop}[Blaschke-Petkantschin Formula Euclidean Case (\cite{Miles1971})]
 Let $\{x_i \ | \ i\in [r]\}$ be a linearly independent set of vectors spanning $V = \text{Span}(\{e_i\ | \ i\in [r]\})\in \text{Gr}(r,d)$. For each $j\in [r]$, let $(\lambda_{jk})_{k\in [r]}\in E^r$ be such that
 $$  x_j=\sum_{k=1}^r \lambda_{jk} e_k, \text{ and let } \Upsilon := \abs{ det(\lambda_{jk})}>0.$$
Then 
$$ \bigwedge_{i=1}^r dV(x_i^d) = \Upsilon^{d-r}d\mu_{r,d}\bigwedge_{i=1}^r dV(x_i^r).$$
\label{prop:Blaschke-Pentkanschin Formula Euclidean Case}
\end{prop}
\begin{proof}
Given $j\in [r]$ and $i\in \{r+1,\dots,n\}$, we have
 \begin{align*}
  dx_j &= \sum_{k=1}^r(d\lambda_{jk}e_k + \lambda_{jk}d e_k)\\
  e_i.dx_j &= \sum_{k=1}^r\lambda_{jk}\omega_{ik}\\
  \bigwedge_{j=1}^re_i.dx_j &= \abs{det(\lambda_{jk})}\bigwedge_{k=1}^r\omega_{ik}\\
  \bigwedge_{i=r+1}^d\bigwedge_{j=1}^re_i.dx_j &= \Upsilon^{d-r}\bigwedge_{i=r+1}^d\bigwedge_{k=1}^r\omega_{ik}\\
  \bigwedge_{i=r+1}^d\bigwedge_{j=1}^re_i.dx_j &= \Upsilon^{d-r}d\mu_{r,d},
 \end{align*}
and multiplying on both sides above by $\bigwedge_{i=1}^r\bigwedge_{j=1}^r e_i.dx_j = \bigwedge_{j=1}^rdV(x_j^r)$, we obtain the desired result.\\
\end{proof}

\subsection{Blaschke-Petkantschin formula for Riemannian Manifolds}

Following \cite{Bobrowski2017}, we obtain a Riemannian generalization of the Blaschke-Petkantschin formula. 
The formula is valid for functions with support close to the diagonal of $M^{k+1}$. It enables us to reparametrise a $(k+1)$-tuple of points near the diagonal of $M^{k+1}$ into local coordinates about their centre.
The change of variables has the following form:
$$ M^{k+1} \longleftrightarrow M \times \mathbb{R} \times Gr(k,d) \times (S^{(k-1)})^{k+1}$$
$$ \vec{y} \longleftrightarrow (c(\vec{y}),u,V,\vec{w})$$
A $(k+1)$-tuple of points in $\vec{y}\in M^{k+1}$, is reparametrised by the centre of this tuple $c(\vec{y})$, the distance of the points from their centre $u$, the $k$-plane in which the pre-image of the points lie in the tangent space at the centre $V$, and the $k+1$ points of the $(k-1)$-sphere upon which they lie $\vec{w}$. 

Suppose that $M\subset \R^d$ is a closed Riemannian manifold and let $\vec{y}=(y_i)_{i=1}^{k+1}\in M^{k+1}$, with centre $c=c(\vec{y})$, with radius $\rho(\vec{y})\leq r$, and local normal coordinates $(x^1,\dots,x^d)$. 
For sufficiently small $r$, we can write for all $i\in [k+1]$
\[
 y_i=\exp_c(v_i),
\]
with
\[
 v_i = \sum_{j=1}^d x^j(y_i)\left(\frac{\partial}{\partial x^i}\right)_c.
\]
Let $\1_r(\vec{y}) = \1\{E_{r_\text{max}}(y_0,...,y_k) \neq \emptyset \text{ and } \rho(\vec{y}) \leq r \}$. Note that this indicator function has support near to the diagonal of $M^{k+1}$ and each tuple in the support of this function has a unique centre. 
It is shown in \cite{Bobrowski2017}, that $\{v_i:i\in [k+1]\}$ have linear dependency and span a $k$-dimensional subspace $V\subset T_{c(\vec{y})} M$ when $c(\vec{y})$ is a critical point. We yield the following change of variable formula.

\begin{lem}[\cite{Bobrowski2017}]
 Let $M$ be a compact closed Riemannian manifold with $M'\subset M$ a submanifold with or without boundary. Let $r_{\max}$ be as in Lemma \ref{lem: Unique Centres}, and $r<r_{\max}$. Then there exists an invariant measure $d\mu_{k,d}(V)$ on $Gr(k,T_cM) = Gr(k,d)$, such that for every $f\in C^{\infty}(M^{k+1};\R)$
 \begin{align*}
  \int_{M^{k+1}} f(\vec{y})\1_r(\vec{y})\1\{c(\vec{y})\in M'\}\abs{d vol_g(\vec{y})} &= \\
  &\int_{M'}\abs{d vol_g(c)}\int_0^r du u^{dk-1}\int_{Gr(k,T_cM)}d\mu_{k,d}(V)\\
  &\times\left(\int_{\mathcal S_1^{k+1}}\Upsilon_1^{d-k}(w)f(\exp_c(uw))\prod_{i=1}^{k+1}\sqrt{\det(g_{\exp_c(uw_i)})}\abs{d vol_{\mathcal S_1(V)}(w_i)}\right).
 \end{align*}
\label{lem: BP Change of Variables}
\end{lem}
\begin{proof}

First note that if $M'$ has positive codimension then both expressions are zero, so assume $M'$ has zero codimension. If $\vec{y}\in M^{k+1}$ is such that $E_{r_{\max}}(\vec{y})\neq \emptyset$, then the induced centre $c(\vec{y})$ is uniquely defined (Lemma \ref{lem: Unique Centres}); hence:
\[
 \{\vec{y}\in M^{k+1}\mid c(\vec{y})\in M'\text{ and }E_{r_{\max}}(\vec{y})\neq \emptyset\} = \bigcup_{c\in M'}\mathcal Y(c),
\]
where $\mathcal Y(c):=\{\vec{y}\in M^{k+1}\mid E_{r_{\max}}(\vec{y})\neq \emptyset \text{ and } c(\vec{y}) = c\}$, and this union is disjoint (by uniqueness of the centre).\\
Thus
\[
\int_{M^{k+1}} f(\vec{y})\1_r(\vec{y})\1\{c(\vec{y})\in M'\}\abs{d vol_g(\vec{y})} = \int_{M'}\abs{d vol_g(c)}\int_{\mathcal Y(c)}f(\vec{y})\1\{\rho(\mathcal Y)\leq r\}\abs{d vol_g(\vec{y})}.
\]
Now fix $c\in M'$ with local normal coordinates $(x^1,\dots,x^d)$; for $\vec{y}\in \mathcal Y(c)$ with $\rho(\vec{y})\leq r<r_{\max}$ (this last condition ensures that $y_i$ can be written as $\exp_c(v_i)$, $v_i\in T_cM$ and $\abs{v_i} = u\leq r$, for all $i\in [k+1]$), we find:\\
\begin{align*}
 \abs{d vol_g(\vec{y})} &= \abs{\wedge_{i=1}^{k+1} d vol_g(y_i)}\\
 &= \abs{\wedge_{i=1}^{k+1}\sqrt{\abs{det(g_{y_i})}}dx^1(y_i)\wedge\dots\wedge d x^d(y_i)}\\
 &= \prod_{i=1}^{k+1}\sqrt{\abs{det(g_{y_i})}}\abs{\wedge_{i=1}^{k+1}dx^1(y_i)\wedge\dots\wedge dx^d(y_i)}\\
 &= \prod_{i=1}^{k+1}\sqrt{\abs{det(g_{y_i})}}\abs{\wedge_{i=1}^{k+1}d vol_{g_{E_d}}(v_i)}\\
 &= \prod_{i=1}^{k+1}\sqrt{\abs{det(g_{y_i})}}\abs{d vol_{g_{E_d}}(\vec{v})},\\
\end{align*}
Note that we have the polar decomposition:
\[
 \abs{dvol_{g_{E_d}}(\vec{v})} =  du \abs{d vol_{\mathcal S_u(E_d)(\vec{v})}},
\]
and so by the Blaschke-Petkantschin formula, since $\{v_i:i\in[k+1]\}$ lies in a $k$-dimensional subspace $V\subset T_c M$, we have
\[
 \abs{d vol_{\mathcal S_u(E_d)}(\vec{v})} = \Upsilon_u(\vec{v})^{d-k} d\mu_{k,d}(V)\abs{d vol_{\mathcal S_u(V)}(\vec{v})},
\]
hence, we deduce that
\[
 \abs{d vol_g(\vec{y})} = \prod_{i=1}^{k+1}\sqrt{\abs{det(g_{y_i})}} du \Upsilon_u(\vec{v})^{d-k} d\mu_{k,d}(V)\abs{d vol_{\mathcal S_u(V)}(\vec{v})}.
\]
This shows that for $c\in M'$
\begin{align*}
 \int_{\mathcal Y(c)}f(\vec{y})\1\{\rho(\mathcal Y)\leq r\}\abs{d vol_g(\vec{y})} &=\\
 &\int_0^r du u^{dk-1}\int_{Gr(k,T_cM)}d\mu_{k,d}(V)\\
  &\times\left(\int_{\mathcal S_1^{k+1}}\Upsilon_1^{d-k}(\vec{w})f(\exp_c(u\vec{w}))\prod_{i=1}^{k+1}\sqrt{\det(g_{\exp_c(uw_i)})}\abs{d vol_{\mathcal S_1(V)}(w_i)}\right),
\end{align*}
and thus the result follows.

\end{proof}

\subsection{The Blaschke-Petkantschin formula for Compact Riemannian Manifold with Non-Empty Boundary}

Using the change of variable formula established in Lemma \ref{lem: BP Change of Variables} for compact closed Riemannian manifolds we attain a formula for Riemannian manifolds with non-empty boundary. 
We note that our formula is altered near the boundary since we must restrict our integral in the tangent space to those points whose image under the exponential map remain in the manifold.
\begin{lem}
Suppose that $M$ is a compact Riemannian manifold with non-empty boundary, let $DM$ be its double manifold and let $M'\subset M$ be a submanifold. Then we have:
\begin{align*}
\int_{M^{k+1}}f(\vec{y})\1_r(\vec{y})&\1\{c(\vec{y})\in M'\}|dvol_g(\vec{y})|\\
&=\int_{c\in M'}|dvol_g(c)|\int_0^rdu\int_{Gr(k,T_cDM)}d\mu_{k,d}(V)\\
&\times\int_{\mathcal (S_u(V)\cap \exp_c^{-1}(M))^{k+1}}\Upsilon_u^{d-k}(v)f(\exp_c(v))\prod_{i=1}^{k+1}\sqrt{|det(g_{\exp_c(v_i)})|}|d vol_{\mathcal S_u^{k+1}(V)}(v)|.
\end{align*}
\end{lem}
\begin{proof}
We have
\begin{align*}
\int_{M^{k+1}}f(\vec{y})\1_r(\vec{y})&\1\{c(\vec{y})\in M'\}|dvol_g(\vec{y})|\\
    &=\int_{DM^{k+1}}f(\vec{y})\1\{y\in M^{k+1}\}\1_r(\vec{y})\1\{c(\vec{y})\in M'\}|dvol_g(\vec{y})|;\\
\end{align*}
the double manifold $DM$ is closed, hence applying the change of variables formula in Lemma \ref{lem: BP Change of Variables} to the function 
$$
\vec{y}\longmapsto f(\vec{y})\1\{y\in M^{k+1}\},
$$
we find
\begin{align*}
&\int_{DM^{k+1}}f(\vec{y})\1\{y\in M^{k+1}\}\1_r(\vec{y})\1\{c(\vec{y})\in M'\}|dvol_g(\vec{y})|=\\
&\int_{c\in M'}|dvol_g(c)|\int_0^rdu\int_{Gr(k,T_cDM)}d\mu_{k,d}(V)\\
&\times\int_{\mathcal S_u(V)^{k+1}}\Upsilon_u^{d-k}(v)f(\exp_c(v))\1\{\exp_c(v)\in M^{k+1}\}\prod_{i=1}^{k+1}\sqrt{|det(g_{\exp_c(v_i)})|}|d vol_{\mathcal S_u^{k+1}(V)}(v)|,
\end{align*}
which gives the lemma.

\end{proof}

\subsection{The Multivariable Blaschke-Petkantschin Formula}

We require another change of variable formula in Section \ref{section: Second Moments} in order to bound the variance of the number of critical points induced by a point process.
We show how to bound a change of variable formula when integrating over two variables in $M^{k+1}$ where $M$ has non-empty boundary. This formula is already used (without proof) in \cite{Bobrowski2017} in the case where $M$ is closed.

\begin{lem}
Let $M$ be a compact Riemannian manifold with non-empty boundary and let $\y_1,\y_2\in M^{k+1}$. Denote the respective centres by $c_1,c_2$.
Define
\[
 \Omega:=\left\lbrace(\y_1,\y_2)\in \left(M^{k+1}\right)^2\mid a\le \rho_M(c_1,c_2)\le b\right\rbrace.
\]
Then there is some constant $C_M$ dependent solely on $M$ for which the following bound holds:
 \begin{align*}
  &\int_{\Omega}f_1(\bold y_1)f_2(\bold y_2)\1_r(\bold y_1,\bold y_2)\abs{dvol_g(\bold y_1,\bold y_2)}\leq\\
  &C_M \int_M\abs{dvol_g(c_1)}\int_a^{b}ds\int_{\mathcal S_1(T_{c_1}M)}s^{d-1}\abs{ dvol_{\mathcal S_1(T_{c_1M})}(w)}\\
  &\times\prod_{i=1}^2\int_{0}^r du_iu_i^{kd-1}\int_{Gr(k,d)}d\mu_{k,d}(V)\int_{\mathcal S_1(V)^{k+1}}\abs{dvol_{\mathcal S_1(V)^{k+1}}(\bold w_i)}f_i(\exp_{c_i}(u_i\bold w_i)). 
 \end{align*}
\label{lem: Multivariable BP Formula}
\end{lem}
\begin{proof}

We use the Blaschke-Petkantschin formula for integrals over one variable in $M^{k+1}$ to attain the result. Given $\y_1\in M^{k+1}$ and $c_1$ the induced center, let 
 \[
  \Omega(\y_1):=\left\lbrace \y_2\in M^{k+1}\mid a\le \rho_M(c_1,c_2)\le b \right\rbrace.
 \]

Denote the above integral on the LHS by $I$, we have
\begin{align*}
 I & = \int_{M^{k+1}}f_1(\y_1)\1_r(\y_1)\left(\int_{\Omega(\y_1)}f_2(\y_2)\1_r(\y_2)\abs{dvol_g(\y_2)}\right)\abs{dvol_g(\y_1)}.
\end{align*}
We first compute the inner integral, for fixed $c_1\in M$.\\
Note that
$$
\int_{\Omega(\y_1)}f_2(\y_2)\1_r(\y_2)\abs{dvol_g(\y_2)} = \int_{M}f_2(\y_2)\1_r(\y_2)\1\{c_2\in A_a^b(c_1)\} \abs{dvol_g(\y_2)},
$$
where $A_a^b(c_1):= B_b(c_1)\setminus B_a(c_1)^o$.\\
By the Blaschke Petkantschin formula for manifolds with non-empty boundary, we then find
\begin{align*}
 &\int_{\Omega(\y_1)}f_2(\y_2)\1_r(\y_2)\abs{dvol_g(\y_2)}\\
 &=\int_{A_a^b(c_1)}\abs{d vol_g(c_2)}\int_0^r du_2\int_{Gr(k,T_c DM)}d\mu_{k,d}(V)\\
 &\times\left(\int_{(\mathcal S_u(V)\cap \exp_{c_2}^{-1}(M))^{k+1}}\abs{dvol_{\mathcal S_u(V)^{k+1}}(\bold v_2)}f_2(\exp_{c_2}(\bold v_2))\Upsilon_u^{d-k}(\bold v_2)\prod_{j=1}^{k+1}\sqrt{\det(g_{\exp_{c_2}(v_j)})}\right)\\
 &\leq C \int_{A_a^b(c_1)}\abs{d vol_g(c_2)}\int_0^r du_2 u_2^{dk-1}\int_{Gr(k,d)}d\mu_{k,d}(V)\\
 &\times\left(\int_{\mathcal S_1(V)^{k+1}}\abs{dvol_{\mathcal S_1(V)^{k+1}}(\bold w_2)}f_2(\exp_{c_2}(u_2\bold w_2))\Upsilon_1^{d-k}(\bold w_2)\prod_{j=1}^{k+1}\sqrt{\det(g_{\exp_{c_2}(u_2w_j)})}\right).
\end{align*}
Using the compactness of $M$, the above is
$$
    \leq C \int_{A_a^b(c_1)}\abs{d vol_g(c_2)}\int_0^r du_2 u_2^{dk-1}\int_{Gr(k,d)}d\mu_{k,d}(V)\int_{\mathcal S_1(V)^{k+1}}\abs{dvol_{\mathcal S_1(V)^{k+1}}(\bold w_2)}f_2(\exp_{c_2}(u_2\bold w_2)).
$$
Furthermore, using the Riemannian approximation results and polar decomposition, we have
\begin{align*}
 \int_{A_a^b(c_1)}\abs{d vol_g(c_2)} &\leq C \int_{a}^{b}ds\int_{\mathcal S_s(T_{c_1}DM)}\abs{d vol_{\mathcal S_s(T_{c_1}DM)}(w)}\\
 &=\int_{a}^{b}ds\int_{\mathcal S_1(T_{c_1}DM)}s^{d-1}\abs{d vol_{\mathcal S_1(T_{c_1}DM)}(w)}.
\end{align*}
The outer integral in the expression of $I$ is estimated again with the Blaschke-Petkantschin formula for manifolds with non-empty boundary. Doing so and combining the result with the above expression for the inner integral of $I$ yields the claimed formula.

\end{proof}

%% file: upper-threshold.tex
\label{Upper Threshold}

In this section we produce an upper bound on the expected number of critical points induced by a Poisson process on our manifold. Similar to the argument presented in \cite{Bobrowski2017} we utilise an auxiliary radius $r_0$ and count critical points with critical value in the range $[r,r_0)$. The auxiliary radius is chosen to be sufficiently large so that the $r_0$ neighbourhood of the Poisson process covers the manifold w.h.p., but is simultaneously sufficiently small so that the number of critical points with critical values in the range $[r,r_0)$ is asymptotically zero.

Let $\Lambda  = n\omega_d r^d $ denote the expected number of points of a Poisson process with intensity $n$ inside a $d$-dimensional ball of radius $r$. Let $\beta_k(r)$ denote the $k^\text{th}$ Betti number of the \v Cech complex associated to the submanifold $B_r(\mathcal{P}_n)$. We shall count critical points of the associated Morse min-type function treating points close to the boundary and far from the boundary separately.

\begin{prop}(Betti Number Upper Bound)
Let $M$ be a $d$-dimensional manifold with boundary. Let $n \to \infty$ and $r,r_0\to 0$ such that $\Lambda \to \infty$, $\Lambda r_0 \to 0$, $\Lambda_{r_0}r_0^2 \to 0$ and $r_0\geq r(\frac{\omega_d}{\kappa}(1+ |\log r|))^{1/d}$ (where $\kappa$ is a constant associated to the manifold apparent in the proof). For all $1\leq k \leq d-1$:
$$ \mathbb{E}[\beta_k(r)] \leq \beta_k(M) + O( n \Lambda^k e^{-\Lambda}, n^{1-\frac{1}{d}}\Lambda^{k}e^{-\frac{1}{2}\Lambda})$$
\label{BettiUpper}
\end{prop}

To prove the above let us follow \cite{Bobrowski2017} in attempting to bound the number of $k$-critical points of the distance function associated to a Poisson process on $M$ with critical value in the range $(r, r_0]$. Let us denote this set of critical points by $C^{\rho_M}_k(r,r_0)$.

\begin{lem}
Let $n \to \infty$ and $r,r_0\to 0$ such that $r= o(r_0)$, $\Lambda \to \infty$, $\Lambda_{r_0}r \to 0$ and $\Lambda_{r_0}r_0^2 \to 0 $ where $\Lambda _{r_0} = \omega_d n r_0^d$. Then for all $1\leq k \leq d$:
$$  \mathbb{E}[C^{\rho_M}_k(r,r_0)] = O(n^{1-\frac{1}{d}}\Lambda^{k-1}e^{-\frac{1}{2}\Lambda}, n\Lambda^{k-1}e^{-\Lambda})$$
\label{Critical point upper bound}
\end{lem}

In what follows let $\mathcal{P}_n$ denote a Poisson process of intensity $n$ on $DM$ with density function $f = \mathbf{1}\{x\in M\}$. Let us use the shorthand notation $\partial M_{r_0}$ for an $r_0$ neighbourhood of the boundary, and $M_{r_0} = M \setminus \partial M_{r_0}$.

Let us define the following indicator functions which we use to count critical points:
\begin{enumerate}
\item $h(\mathcal{Y}) = \mathbf{1}\{0 \in Grad(\mathcal{Y})\}$
\item $h_{r,r_0}(\mathcal{Y}) = h(\mathcal{Y})\mathbf{1}\{r \leq \rho(\mathcal{Y}) < r_0\}$

\item $g_{r,r_0}(\mathcal{Y},\mathcal{P}_n) = h_{r,r_0}(\mathcal{Y})\mathbf{1}\{B(\mathcal{Y}) \cap \mathcal{P}_n = \emptyset \}$; where $B(\mathcal Y)$ is the open ball $B_{\rho(\mathcal Y)} (c(\mathcal Y))$
\item $g_{r,r_0}^{M_{r_0}}(\mathcal{Y},\mathcal{P}_n) = g_{r,r_0}(\mathcal{Y},\mathcal{P}_n)\mathbf{1}\{c(\mathcal{Y})\in M_{r_0}\}$
\item $g_{r,r_0}^{\partial M_{r_0}}(\mathcal{Y},\mathcal{P}_n) = g_{r,r_0}(\mathcal{Y},\mathcal{P}_n)\mathbf{1}\{c(\mathcal{Y})\in \partial M_{r_0}\}$
\end{enumerate}

Then we can express the number of critical points as the following sum over
generic $\mathcal{Y}\subset \mathcal{P}_n$:

$$ |C^{\rho_M}_k(r,r_0)| = \sum_{|\mathcal{Y}| = k+1} |g_{r,r_0}^{M_{r_0}}(\mathcal{Y},\mathcal{P}_n)| + |g_{r,r_0}^{\partial M_{r_0}}(\mathcal{Y},\mathcal{P}_n) |$$

We can count the critical points in $M_{r_0}$ and $\partial M_{r_0}$ separately. Let us denote the critical points in $\partial M_{r_0}$ with critical value in the interval $[r,r_0)$ as $C_{\partial M_{r_0}}$, and the critical points in $M_{r_0}$ with critical value in the interval $[r,r_0)$ as $C_{ M_{r_0}}$. The counting of the critical points $C_{M_{r_0}}$ is not affected by the presence of a non-trivial boundary and so we can use the analysis of $\cite{Bobrowski2017}$ to bound the term $\mathbb{E}[|C_{M_{r_0}}|] $  by an order of $O(n\Lambda^{k-1}e^{-\Lambda})$. 

We require new analysis to calculate an upper bound for the expected number of critical points $\mathbb{E}[|C_{\partial M_{r_0}}|] $. Applying Palm Theory (Theorem \ref{thm:PalmTheory1}) to the uniform Poisson Process on $M$ with intensity $n$ we yield:

$$\mathbb{E}[|C_{\partial M_{r_0}}|]  =  \frac{n^{k+1}}{(k+1)!}\mathbb{E}[g_{r,r_0}^{\partial M_{r_0}}(\mathcal{Y}' ,\mathcal{Y}' \cup\mathcal{P}_n)] $$

Let us now condition on a given sample $\mathcal{Y}'$ of $k+1$ points of $M$ with a view to then integrating over the whole double manifold. Recall $B(\mathcal{Y}')$ denotes the ball centred at $c(\mathcal{Y}')$ of radius $\rho(\mathcal{Y}')$ in $DM$, and that the number of points in a given subset of a uniform Poisson process has Poisson distribution with parameter proportional to the volume of the given subset.

$$ \mathbb{E}[\mathbf{1}\{B(\mathcal{Y}') \cap  \mathcal{P}_n = \emptyset\}| \mathcal{Y}'] = \mathbb{P}(\mathcal{P}_n(B(\mathcal{Y}'))=0 | \mathcal{Y}') = e^{-n\text{Vol}(B(\mathcal{Y}')\cap M)}$$

$$ \mathbb{E}[|C_{\partial M_{r_0}}|] = \frac{n^{k+1}}{(k+1)!}
 \int_{{DM}^{k+1}} \mathbf{1}\{c(\mathbf{y})\in {\partial M_{r_0}}\}h_{r,r_0}(\mathbf{y})e^{-n\text{Vol}(B(\mathbf{y})\cap M)} |\text{dvol}_g(\mathbf{y})|$$

For sufficiently small $r_0$ we have the following lower bound for the volume of a ball with centre at distance $\delta$ from the boundary via Proposition \ref{prop:Half-BallVolumes}:

$$\text{Vol}(B(\mathbf{y})\cap M) \geq (\frac{1}{2}-O(r))\text{Vol}(B(\mathbf{y})) + \frac{1}{2}\omega_d \delta^d $$.

\begin{figure}[ht]
\centering
\includegraphics[scale=0.4]{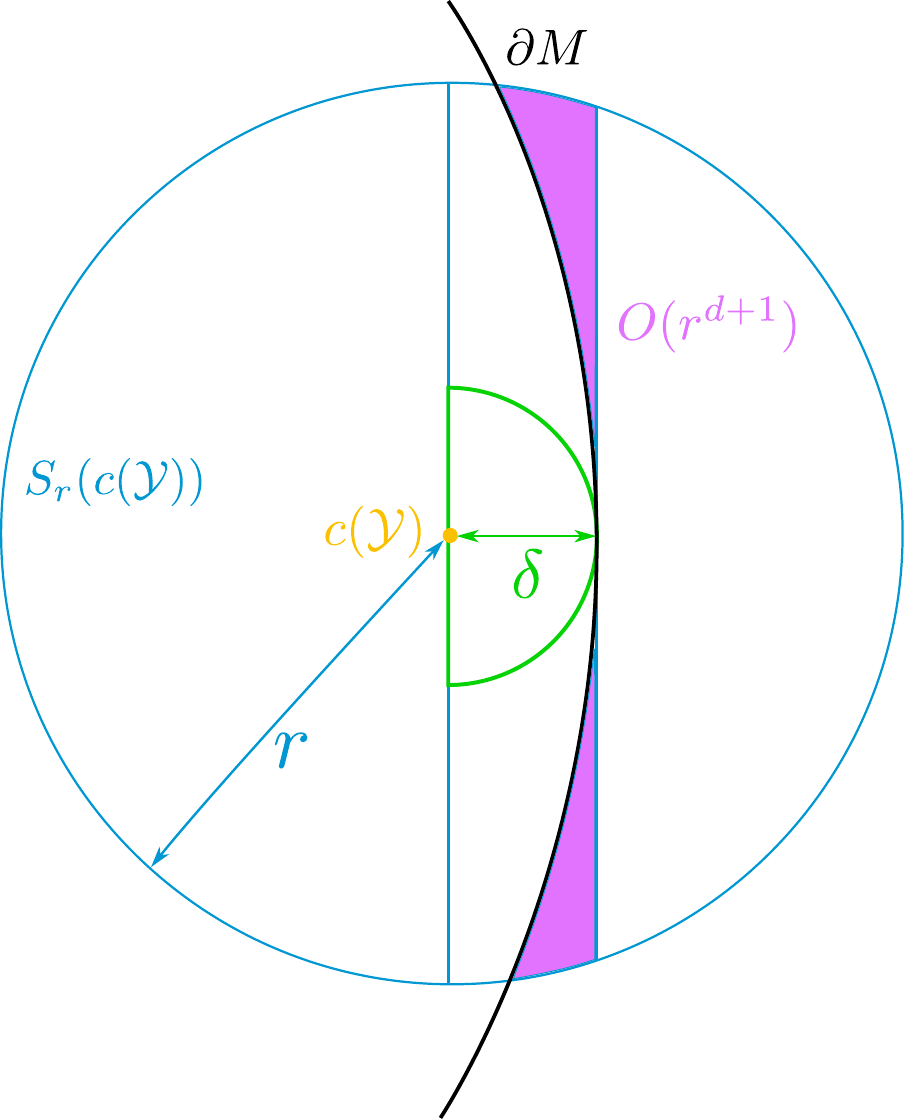}
\caption{The contributions to our lower bound for $\text{Vol}(B(\mathbf{y})\cap M)$, are illustrated in this figure. Proposition \ref{prop:Half-BallVolumes} establishes that the shaded region has volume $O(r^{d+1})$. We note that for $\delta< \tau_M$ ($\tau_M$ the reach of the manifold) then the half-ball of radius $\delta$ contributes $\frac{1}{2}\omega_d \delta^d $ to the volume.}
\label{fig:Upper-Bound-Volumes}
\end{figure}
Note that $e^{nO(r)\text{Vol}(B(\mathbf{y}))} = O(1)$ since $ \Lambda_{r_0}r \to 0$, and so we have an upper bound: 



$$ \mathbb{E}[|C_{\partial M_{r_0}}|] \leq C
\frac{n^{k+1}}{(k+1)!}
 \int_{DM^{k+1}} \mathbf{1}\{c(\mathbf{y})\in {\partial M_{r_0}}\}h_{r,r_0}(\mathbf{y})e^{-n\frac{1}{2}\left(\text{Vol}(B(\mathbf{y})) + \omega_d \delta^d \right)} |\text{dvol}_g(\mathbf{y})|$$

Using the Blaschke-Petkantschin formula from Section \ref{Blaschke} with $f(\mathbf{y}) = h_{r,r_0}(\mathbf{y})e^{-n\frac{1}{2}\left(\text{Vol}(B(\mathbf{y})) + \omega_d \delta^d \right)}$, $c=c(\mathbf{y})$, $u=\rho(\mathbf{y})$ and $\mathbf{y} = \text{exp}_c(u\mathbf{w})$ we attain the bound:

\begin{align*}
 \mathbb{E}[|C_{\partial M_{r_0}}|] \leq 
\frac{n^{k+1}}{(k+1)!}
&\int_{\partial M_{r_0}} |\text{dvol}_g(c)|
\int_r^{r_0} du \  u^{dk-1} \int_{Gr(k,T_cM)} d\mu_{k,d}(V) \\
& \times  \prod_{i=1}^{k} \left( \int_{\mathbb{S}_1(V)} \sqrt{|\text{det}(g_{\text{exp}_c(uw_i)})|} |\text{dvol}_{\mathbb{S}_1(V)}(w_i) |\right)
\mathbf{\Upsilon}^{d-k}_1(\mathbf{w})f(\text{exp}_c(u\mathbf{w}))
\end{align*}

Using compactness of the Grassmanian we can bound the final integral over $(\mathbb{S}_1(V))^k$ by some maximising subspace. The parallelogram volume $\mathbf{\Upsilon}^{d-k}_1(\mathbf{w})$ is taken over unit vectors and so bounded. Since our manifold is compact and $g$ is smooth we may bound $\sqrt{|\text{det}(g_{\text{exp}_c(uw_i)})|}$ by some constant. (For details of the bounding constants see \cite{Bobrowski2017}). Thus we yield for some constant $C$:


\begin{align*}
 \mathbb{E}[|C_{\partial M_{r_0}}|] \leq 
Cn^{k+1}\int_{\partial M_{r_0}} |\text{dvol}_g(c)|
\int_r^{r_0} du \  u^{dk-1} e^{-n\frac{1}{2}\left(\text{Vol}(B_u(c)) + \omega_d \delta^d \right)}
\end{align*}

Since $\Lambda_{r_0}r_0^2 \to 0$ we may bound the term $e^{-n\frac{1}{2}\left(\text{Vol}(B_u(c)) + \omega_d \delta^d \right)}$
using the second order Taylor expansion of $\text{exp}$ and Corollary \ref{cor: Ball Volume} :

$$ e^{-n\frac{1}{2}\left(\text{Vol}(B_u(c)) + \omega_d \delta^d \right)} \leq e^{-n\frac{1}{2}\omega_d\left( u^d+ \delta^d\right)} \left( 1 + s_{\text{max}}n\omega_dr_0^{d+2} \right) $$

Let us apply a change of variables $s = \frac{u}{r}$ and additionally separate our integral to integrate over the distance of the critical point to the boundary. We make a change of variables using the diffeomorphism $\partial M_{r_0} \cong \partial M \times [0,r_0)$ given to us from the Collar Neighbourhood Theorem. Since our manifold is compact the Jacobian term introduced will be bounded by a constant which we absorb into the constant term $C$.

\begin{align*}
\mathbb{E}[|C_{\partial M_{r_0}}|]  &\leq 
Cn^{k+1} ( 1 + s_{\text{max}}n\omega_dr_0^{d+2})
\int_{\partial M_{r_0}} e^{-n\frac{1}{2}\omega_d \delta^d} |\text{dvol}_g(c)|
\int_1^{\frac{r_0}{r}} ds \  r^{dk} s^{dk-1}  e^{-n\frac{1}{2}\omega_dr^ds^d} \\
& \leq Cn\Lambda^k( 1 + s_{\text{max}}\Lambda_{r_0}r_0^{2}) 
\int_{0}^{r_0} e^{-n\frac{1}{2}\omega_d \delta^d} d \delta
\int_1^{\frac{r_0}{r}} ds \ s^{dk-1}  e^{-\frac{1}{2} \Lambda s^d}
\end{align*}

We may bound the first integral by a term $O(n^{-\frac{1}{d}})$ by changing variables $(t = \frac{1}{2}\Lambda_\delta)$ and observing that the integral takes the form of an upper incomplete gamma function:

$$ \int_{0}^{r_0} e^{-n\frac{1}{2}\omega_d \delta^d} d \delta = C n^{-\frac{1}{d}} \int_{0}^{\frac{1}{2}\Lambda_{r_0}} t^{\frac{1}{d}-1}e^{-t} dt \leq C \Gamma(\frac{1}{d}) n^{-\frac{1}{d}}$$

The last integral also has the form of an upper incomplete gamma function.

$$ \Gamma(k,x) = \int_x^\infty t^{k-1}e^{-t}dt = (k-1)!e^{-x} \sum_{i=0}^{k-1} \frac{ x^i}{i!}$$

Let $t = \frac{1}{2} \Lambda s^d$
$$ \int_1^{\frac{r_0}{r}}   \left(\frac{1}{2} \Lambda\right)^k s^{dk-1}  e^{-\frac{1}{2} \Lambda s^d} ds =  \frac{1}{d} \int_{\frac{1}{2}\Lambda}^{\frac{1}{2} \Lambda_{r_0}}  t^{k-1} e^{-t} dt = \Gamma(k,\frac{1}{2} \Lambda) - \Gamma(k,\frac{1}{2}\Lambda_{r_0})
$$

Absorbing surplus constants into the term $C$ we attain:

\begin{align*}
\mathbb{E}[|C_{\partial M_{r_0}}|] &\leq 
  Cn^{1-\frac{1}{d}}( 1 + s_{\text{max}}\Lambda_{r_0}r_0^{2}) 
\left(  e^{-\frac{1}{2}\Lambda}\sum_{j=0}^{k-1} \frac{(\frac{1}{2}\Lambda)^j}{j!} - e^{-\frac{1}{2}\Lambda_{r_0}}\sum_{j=0}^{k-1} \frac{(\frac{1}{2}\Lambda_{r_0})^j}{j!} \right)
\end{align*}

Using the assumptions that $\Lambda_{r_0}r_0^2 \to 0$, $\Lambda \to \infty$ and $r = o(r_0)$ yields that $\mathbb{E}[|C_{\partial M_{r_0}}|] = O(n^{1-\frac{1}{d}}\Lambda^{k-1} e^{-\frac{1}{2}\Lambda})$. Thus we prove the Lemma \ref{Critical point upper bound}.

Let us now use Lemma \ref{Critical point upper bound} to establish Proposition \ref{BettiUpper}.

\begin{proof}(Proposition \ref{BettiUpper})

Having established Lemma \ref{Critical point upper bound} this proof proceeds like the proof of Proposition 6.1 in \cite{Bobrowski2017} mutatis mutandis.

Let us define $\hat{\beta}_k(r) = \beta_k(r) - \beta_k(M)$ then use Lemma \ref{Critical point upper bound} to  bound $\mathbb{E}[\hat{\beta}_k(r)]$. Let's condition on the event $E = \{M \subset B_{r_0}(\mathcal{P}_n)\}$ that we cover our manifold with the radius $r_0$ neighbourhood of our Poisson Process. Conditioning on $E$ we get:
$$ \mathbb{E}[\hat{\beta}_k(r)] = \mathbb{E}[\hat{\beta}_k(r)| E]\mathbb{P}(E) + \mathbb{E}[\hat{\beta}_k(r)| E^c]\mathbb{P}(E^c)$$

We can bound the term $\mathbb{E}[\hat{\beta}_k(r)| E]\mathbb{P}(E)$ by observing that any non-trivial $k$-cycle in $\mathcal{C}_r$ that is trivial in $M$ must be annihilated by a $(k+1)$-critical point in $C_{k+1}^{\rho_M}(r,r_0)$. Thus assuming $\Lambda \to \infty$, $r = o(r_0)$ and $ \Lambda_{r_0}r_0^2 \to 0$, using Lemma \ref{Critical point upper bound} we may bound this term by $\mathbb{E}[|C_{k+1}^{\rho_M}(r,r_0)|] = O(n^{1-\frac{1}{d}}\Lambda^{k}e^{-\frac{1}{2}\Lambda} , n\Lambda^{k}e^{-\Lambda})$

The second term, $\mathbb{E}[\hat{\beta}_k(r)| E^c]\mathbb{P}(E^c)$ can be bounded using an $\frac{r_0}{2}$-net to bound the non-coverage probability. 
Let us begin by bounding  $\mathbb{E}[\beta_k(r)| E^c]\mathbb{P}(E^c)$. A coarse upper bound for $\beta_k(r)$ is the number of $k$-dimensional faces of $\mathcal{C}_r$:

\begin{align*}
\mathbb{E}[\beta_k(r)| E^c]\mathbb{P}(E^c) & \leq \mathbb{E}\left[ \binom{|\mathcal{P}_n|}{k+1} | E^c \right]\mathbb{P}(E^c) \\
&= \sum_{k+1}^{\infty} \binom{m}{k+1} \mathbb{P}(|\mathcal{P}_n|=m \ |\ E^c)\mathbb{P}(E^c) \\
&= \sum_{k+1}^{\infty} \binom{m}{k+1} \mathbb{P}( E^c\ |\ |\mathcal{P}_n|=m )\mathbb{P}(|\mathcal{P}_n|=m)
\end{align*}

Since $\mathcal{P}_n$ is a Poisson process of intensity $n$ then $ \mathbb{P}(|\mathcal{P}_n|=m) = \frac{e^{-n}n^m}{m!}$ and conditioned on $\{|\mathcal{P}_n|=m\}$ we may write $\mathcal{P}_n$ as a set of $m$ independent uniformly distributed random variables $\chi_m = \{X_1,...,X_m\}$
$$ \mathbb{P}( E^c\ |\ |\mathcal{P}_n|=m ) = \mathbb{P}(B_{r_0}(\chi_m) \neq M)$$

Let $\mathcal{N}$ be a $\frac{r_0}{2}$-net for $M$ chosen such that $| \mathcal{N}| \leq c_dr_0^{-d}$ for some constant $c_d$ dependent only on the dimension of our manifold and the metric, and such that at least $\frac{1}{2}$ of the volume of radius $\frac{r_0}{2}$ ball lies in $M$. We then use this net to bound the non-coverage probability:

\begin{align*}
 \mathbb{P}(B_{r_0}(\chi_m) \neq M) &\leq 
\sum_{x\in\mathcal{N}}\mathbb{P}\left({\rho_{\chi_m}}(x) > \frac{r_0}{2}\right) \\
&\leq c_dr_0^{-d}\max_{x\in\mathcal{N}}\{(1-\text{Vol}(B_{r_0}(x)))^m\}
\\
&\leq c_dr_0^{-d}(1-\omega_d2^{-(d+1)}(1-s_\text{max}r_0^2)r_0^d)^m
\\
&\leq c_dr_0^{-d}(1-\kappa r_0^d)^m
\end{align*}

Where $s_\text{max}$ is a constant intrinsic to the manifold as defined in Corollary \ref{cor: Ball Volume}. Thus we yield:

\begin{align*}
\mathbb{E}[\beta_k(r)| E^c]\mathbb{P}(E^c) &\leq 
\sum_{k+1}^{\infty} \binom{m}{k+1} \mathbb{P}( E^c\ |\ |\mathcal{P}_n|=m )\mathbb{P}(|\mathcal{P}_n|=m)
\\
& \leq \sum_{k+1}^{\infty} \binom{m}{k+1} c_dr_0^{-d}(1-\kappa r_0^d)^m\frac{e^{-n}n^m}{m!}
\\
& \leq r_0^{-d}n^{k+1}(1-\kappa r_0^d)^{k+1}e^{-\kappa n r_0^d}\sum_{j=0}^{\infty} \frac{c_d}{(m-j)!}\frac{e^{-n(1-\kappa r_0^d)}n^j(1-\kappa r_0^d)^j}{j!}
\\
& \leq Cr_0^{-d}n^{k+1}(1-\kappa r_0^d)^{k+1}e^{-\kappa n r_0^d} 
\leq Cr_0^{-d}n^{k+1}e^{-\kappa n r_0^d}
\end{align*}

Calculating the asymptotic behaviour of $r_0^{-d}n^{k+1}e^{-\kappa n r_0^d}/n\Lambda^k e^{-\Lambda}$ given the conditions that $r_0 \geq r(\frac{\omega_d}{\kappa}(1+ |\log r|))^{1/d}$ and $\Lambda r \to 0$ we see:

$$ \frac{r_0^{-d}n^{k+1}e^{-\kappa n r_0^d}}{n\Lambda^k e^{-\Lambda}} = \frac{n^ke^{\Lambda - \kappa n r_0^d }}{r_0^{d}\Lambda^k} \leq \frac{n^ke^{-|\log r|\Lambda}}{r^d\Lambda^{k}(\frac{\omega_d}{\kappa}(1+ |\log r|))} 
= \frac{n^k r^{\Lambda -d}}{\Lambda^{k}(\frac{\omega_d}{\kappa}(1+ |\log r|))} \to 0 $$

Hence $\mathbb{E}[\beta_k(r)| E^c]\mathbb{P}(E^c) = o(n\Lambda^k e^{-\Lambda})$. Then we note that $\mathbb{P}(E^c) = o(n\Lambda^k e^{-\Lambda})$ by an almost identical argument and so $\mathbb{E}[\hat{\beta}_k(r)| E^c]\mathbb{P}(E^c) = \mathbb{E}[\beta_k(r)| E^c]\mathbb{P}(E^c) - \beta_k(M)\mathbb{P}(E^c) = o(n\Lambda^k e^{-\Lambda})$.

\end{proof}

%% file: lower-threshold.tex
\label{Lower Threshold}

In this section we shall find a lower bound for the expected Betti numbers of the \v Cech complex $\mathcal{C}(n,r)$. We utilise the concept of a special type of critical point called a $\Theta$-cycle defined in \cite{bobrowski_vanishing_2015}. Such a critical point is guaranteed to induce a non-trivial cycle in the homology of the resulting \v Cech complex. The lower bound shows that if the convergence of $\Lambda \to \infty $ is sufficiently slow then, w.h.p the \v Cech complex will have Betti numbers larger than that of the manifold. Thus we will attain a lower threshold for $\Lambda$.

The paper \cite{Bobrowski2017} counts $\Theta$-cycles on a general compact closed Riemannian manifold with critical values in the range $(r_1,r]$. The auxiliary radius $r_1$ is chosen to be sufficiently small that there are a large number of critical points with critical value in the range $(r_1,r]$ but simultaneously sufficiently large that a $\Theta$-cycle with critical value $r_1$ persists and remains a $\Theta$-cycle in the \v Cech complex at radius $r$.

The conditions that determine a critical point to be a $\Theta$-cycle are all local conditions. As such we can replicate the analysis in \cite{Bobrowski2017} to count the $\Theta$-cycles in a manifold with boundary which are sufficiently distant from the boundary.
It transpires that this lower bound can be improved by counting a collection of cycles which occur close to the boundary, which we shall call $\Theta$-like-cycles.


Let us proceed to define the conditions given in \cite{Bobrowski2017} that determine a critical point to be a $\Theta$-cycle. Let $\mathcal{Y} \subset \mathcal{P} \subset M$ be a generic subset of a point process on $M$ inducing centre $c(\mathcal{Y})$ and $\varepsilon \in (0,1)$.
Let us denote the closed annulus in $M$ by:
$$ A_\varepsilon(\mathcal{Y}) = \overline{B_{\rho(\mathcal{Y})}(c(\mathcal{Y}))} \setminus B_{\varepsilon\rho(\mathcal{Y})}(c(\mathcal{Y}))$$

Intuitively a $\Theta$-cycle is formed at critical point $c$ when an annulus surrounding $c$ is covered and the critical point $c$ introduces a $k$ simplex which cuts across this annulus and introduces an erroneous homological cycle. If our centre is close to the boundary, our annulus will be cut by the boundary and so we must modify the argument in \cite{Bobrowski2017} to apply to this case. 

Lemma \ref{Theta-Cycle} is a slight modification of that presented in \cite{bobrowski_vanishing_2015} and having taken into consideration that the critical point lies far from the boundary the proof follows identically.

\begin{lem}\cite{bobrowski_vanishing_2015}
Let $\mathcal{Y} \subset \mathcal{P} \subset M$ with $\mathcal{Y}$ inducing a critical point $c(\mathcal{Y}) $ of index $k$. Let us define:
$$\phi(\mathcal{Y}) = \frac{1}{2\rho(\mathcal{Y})} \min_{v\in \partial \Delta(\mathcal{Y})}|v| $$
where $\Delta(\mathcal{Y}) = Grad(\mathcal{Y})$. Suppose $\rho(\mathcal{Y}) < r_{\text{max}}$, $\text{dist}(c(\mathcal{Y}),\partial M) > \rho(\mathcal{Y})$ and $A_\phi(\mathcal{Y}) \subset B_{\rho(\mathcal{Y})}(\mathcal{P})$, then the critical point $c(\mathcal{Y})$ generates a new non-trivial cycle in $H_k(B_{\rho(\mathcal{Y})}(\mathcal{P}))$ which we call a $\Theta$-cycle.
\label{Theta-Cycle}
\end{lem}

The following technical Lemma is established by showing that a $\Theta$-cycle with critical value in the range $(r_1,r]$ will remain a $\Theta$-cycle at radius $r$.

\begin{lem}\cite{Bobrowski2017}
Let $\Theta_k^{\varepsilon}(r_1,r)$ denote the number of $\Theta$-cycles induced by subsets of $(k+1)$ points $\mathcal{Y}$ with the properties:
$$ \rho(\mathcal{Y}) \in (r_1,r] ,\ \ B_{r_2}(c(\mathcal{Y})) \cap \mathcal{P} = \mathcal{Y} ,\ \ \phi(\mathcal{Y})\geq \varepsilon, \ \ \text{dist}(c(\mathcal{Y}),\partial M) > r_2$$
Suppose $r_2>r >0$ and that $r_1 > r\sqrt{1- \frac{1}{c_g^2}(\frac{r_2}{r}-1)^2}$. For any $\varepsilon \in (0,1)$ we have that $\beta_k(r) \geq \Theta_k^{\varepsilon}(r_1,r)$.
\end{lem}

Let us use the following notation for indicator functions which track when a subset of $(k+1)$ elements $\mathcal{Y}\subset \mathcal{P}$ induces an element of $\Theta_k^{\varepsilon}(r_1,r)$.

\begin{enumerate}
\item $h(\mathcal{Y}) = \mathbf{1}\{0 \in Grad(\mathcal{Y})\}$
\item $h_{r_1,r}(\mathcal{Y}) = h(\mathcal{Y})\mathbf{1}\{r_1 < \rho(\mathcal{Y}) \leq r\}$
\item $h_{r}^\varepsilon(\mathcal{Y}) = h_{r_1,r}(\mathcal{Y})\mathbf{1}\{\phi(\mathcal{Y})\geq \varepsilon\}\mathbf{1}\{\text{dist}(c(\mathcal{Y}),\partial M) > r_2\}$
\item $g_{r}^\varepsilon(\mathcal{Y},\mathcal{P}) = h_{r}^\varepsilon(\mathcal{Y})\mathbf{1}\{ B_{r_2}(c(\mathcal{Y})) \cap (\mathcal{P}\setminus \mathcal{Y}) = \emptyset \} \mathbf{1}\{ A_\varepsilon \subset B_{\rho(\mathcal{Y})}(\mathcal{P}) \}$
\end{enumerate}

Hence we may write $\Theta_k^{\varepsilon}(r_1,r)$ as the sum:

$$ \Theta_k^{\varepsilon}(r_1,r) = \sum_{|\mathcal{Y}| = k+1} g_{r}^\varepsilon(\mathcal{Y},\mathcal{P}) $$

\begin{lem}\cite{Bobrowski2017}
Suppose $\varepsilon>0$ is sufficiently small (independent of $r$ and $n$), and $r>0$ such that $\Lambda \to \infty, \Lambda r^2 \to 0$. Then for suitably chosen $r_1,r_2$ with $r_2>r>r_1>0$ we have that:
$$ \mathbb{E}[ \Theta_k^{\varepsilon}(r_1,r)] = \Omega(n\Lambda^{k-2}e^{-\Lambda})$$
\end{lem}

The proof of the above lemma follows almost identically to the proof supplied in \cite{Bobrowski2017}.

\begin{prop}\cite{Bobrowski2017}
Suppose $\varepsilon>0$ is sufficiently small, $\omega(n) \to \infty$, $\gamma \in (0,1)$ is fixed and $\Lambda = \log n +(k-2) \log\log n -\omega(n)$, then for all $1\leq k\leq d-1$:
$$ \lim_{n\to \infty}\mathbb{P}( \Theta_k^{\varepsilon}(r_1,r) > \gamma\mathbb{E}[ \Theta_k^{\varepsilon}(r_1,r)]) = 1$$
\end{prop}

This Proposition is proven using a second moment argument based on Chebyshev's inequality and the details may be found in \cite{Bobrowski2017}. 

The results from this section yield a lower threshold for $\Lambda$:
\begin{prop}
Let $M$ be a unit volume, compact, Riemannian manifold with boundary. Let $\Lambda =  \log n +(k-2) \log\log n -\omega(n)$ then for all $1\leq k\leq d-1$:
$$ \lim_{n\to \infty} \mathbb{P}(H_k(\mathcal{C}(n,r)) \cong H_k(M)) = 
0$$
\end{prop}

This threshold is identical to that found in \cite{Bobrowski2017}. We adapted the proof from \cite{Bobrowski2017} in order that we only counted $\Theta$-cycles sufficiently far from the boundary that the local considerations did not detect the manifold had a boundary. 

%% file: lower-threshold-refined.tex
We shall now adapt the argument counting $\Theta$-cycles to establish that a large number of $\Theta$-like-cycles are present close to the boundary. From this we can deduce a greater lower threshold. Let $\Theta_k^{\varepsilon, \partial M}(r_1,r)$ denote the number $\Theta$-like-cycles of index $k$ with critical value in the range $(r_1,r]$. This type of critical point induces an erroneous homological cycle near the boundary and will be defined thoroughly below, (Definition \ref{defn: Theta-like-cycle}).

Similarly to the previous section we count $\Theta$-like-cycles with critical values in the range $(r_1,r]$ with auxilliary radius $r_2$ used as in the previous section to ensure a $\Theta$-like-cycle with critical value $r_1$ remains a $\Theta$-like-cycle at scale $r$.

\begin{lem}
Let $M$ be a unit volume compact Riemannian manifold with boundary. Suppose $n\to \infty$ and $r\to 0$ such that $\Lambda \to \infty$, $\Lambda r^2 \to 0$ and suppose $\alpha = \frac{1}{2} + O((\log n)^{-1})$. Then for all $1 \leq k \leq d-1$ we have that: $$\mathbb{E}[|\Theta_k^{\varepsilon, \partial M}(r_1,r)|] = \Omega( n\Lambda^{k-2}e^{-\alpha\Lambda} r (\log n)^{-(k+1)})$$ 
\label{Critical point lower bound}
\end{lem}
Heuristically a $\Theta$-cycle is a critical point lying in the centre of a surrounding annulus, and thus induces a new cycle. For a critical point close to the boundary, the whole annulus is not contained in the manifold and so we define a $\Theta$-like-cycle for this corresponding situation. The annuli which are cut by the boundary form a cup which for certain critical points give rise to non-trivial homological cycles. We shall identify which critical points introduce new homological cycles.

The argument to establish Lemma \ref{Critical point lower bound} in this section is fairly delicate and has the following structure. We first define a special class of critical points, $\Theta$-like-cycles, and establish that these critical points introduce new spurious homological cycles near the boundary (Lemma \ref{lem: Theta-like-cycles induce homology}).
We then wish to bound from below the expected number of $\Theta$-like-cycles. In order for a critical point to be a $\Theta$-like-cycle we require that both a partial annulus surrounding the critical point is covered and that the simplex induced by the critical point is approximately tangential to the boundary. We define the partial annulus $A^{({\varphi})}_\varepsilon(c)$ which is constructed in order that if the points inducing the centre all lie in this annulus then the simplex introduced by the critical point is  approximately tangential to the boundary. Then for suitable $\varphi$ the partial annulus is covered and so both conditions are met and the critical point $c$ thus induces a $\Theta$-like-cycle.





\begin{figure}[ht]

\centering
\includegraphics[scale=0.4]{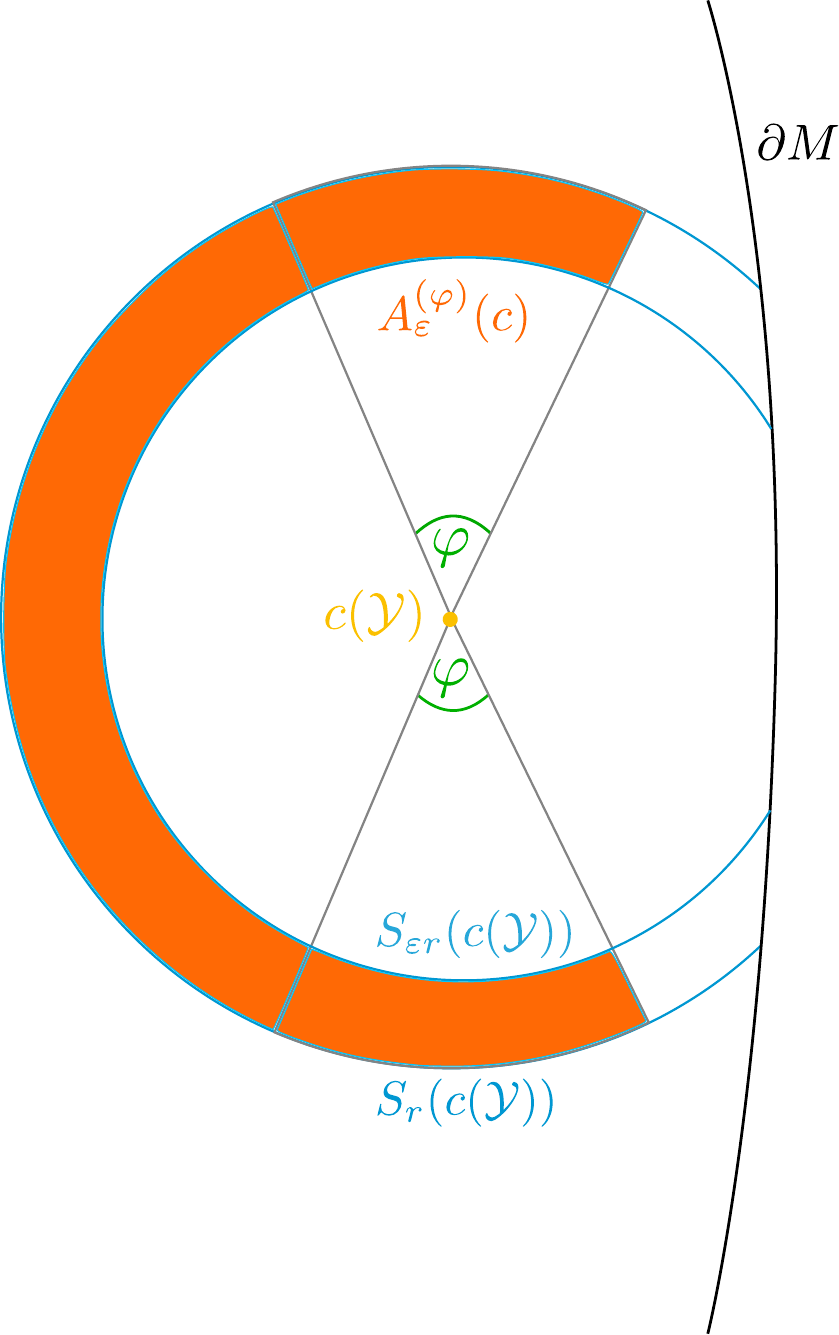}
\caption{The partial annulus $A^{({\varphi})}_\varepsilon(c)$}
\label{fig: Partial Annulus}

\caption{A sketch of the partial annulus $A^{({\varphi})}_\varepsilon(c)$. Our calculations show that if the centre $c(\mathcal{Y})$ is at distance  $\sim \frac{r}{\log n}$ from the boundary and $\varphi \sim (\log n)^{-1}$ then $A^{({\varphi})}_\varepsilon(c)$ is covered w.h.p.  }
\end{figure}

\begin{defn}(The Partial Annulus)
Let $c$ be a critical point with radius $\rho$ such that $\text{dist}(\partial M, c) < \rho$. Let $p$ be the closest point to $c$ on $\partial M$. Let $\vec{n} := \exp_c^{-1}{p}$, which can be thought of as the vector pointing in the direction of the normal to the boundary at $p$. We define the partial annulus: $$A^{(\beta)}_\varepsilon(c) = \{x \in A_\varepsilon(c) \ :\ \textrm{exp}_c^{-1} x  \textrm{ makes angle greater than } \pi/2 - \beta/2 \textrm{ with } \vec{n}\}$$

See Figure \ref{fig: Partial Annulus}.
\end{defn}

\subsection{$\Theta$-like-cycles Induce Homology}

Let us define a function $\psi(\mathcal{Y},\varphi)$ for $\Theta$-like-cycles that plays the same role as $\phi$ did for $\Theta$-cycles:

$$ \psi(\mathcal{Y},\varphi) = \frac{1}{2} \sup \{\varepsilon \geq 0 \ | \ \partial \Delta(\mathcal{Y}) \subset {A}^{(\varphi)}_\varepsilon(c(\mathcal{Y})) \}$$

Note that like $\phi$, the term $\psi(\mathcal{Y},\varphi)$, is not sensitive to the scale $\rho(\mathcal{Y})$ in the sense that $\psi(\mathcal{Y},\varphi)$ is dependent on the distribution of the points $\mathcal{Y}$ on the sphere centred at $c(\mathcal{Y}))$ and so has a lower bound which does not scale with the critical value $\rho(\mathcal{Y})$. We will sometimes abbreviate $\psi(\mathcal{Y},\varphi)$ to just $\psi$.

\begin{lem}($\Theta$-like-cycles induce non-trivial homological cycles)

Let $\mathcal{Y}\subset\mathcal{P}\subset M$ be a set of points inducing an index $k$ critical point with critical value $\rho$ and $1 \leq k\leq d-1$. Let $\psi = \psi(\mathcal{Y})>0$ and suppose ${A}^{(\varphi)}_{\psi}(c(\mathcal{Y})) \subset B_\rho(\mathcal{P})$, then $c(\mathcal{Y})$ induces a non-trivial cycle in $H_k(\mathcal{C}(\mathcal{P},\rho))$.
\label{lem: Theta-like-cycles induce homology}
\end{lem}

\begin{proof}

We may assume $\mathcal{P}$ to be generic in the sense that the critical values of each simplex are distinct, so there is some $\rho_- < \rho$ such that $\mathcal{C}(\mathcal{P},\rho) = \mathcal{C}(\mathcal{P},\rho_-) \cup \Delta$ where $\Delta$ is the $k$-simplex $\mathcal{Y}$. Moreover the boundary $\partial \Delta \in \mathcal{C}(\mathcal{P},\rho_-)$, and $\Delta$ is not the face of any higher simplex in $\mathcal{C}(\mathcal{P},\rho)$ by the construction of the \v{C}ech complex.

Suppose that $\partial \Delta$ is a boundary in $\mathcal{C}(\mathcal{P},\rho_-)$ so that $\partial \Delta = \partial \gamma$ for some $k$ chain $\gamma \in \mathcal{C}(\mathcal{P},\rho_-)$. Then clearly $\Delta - \gamma$ is a $k$-cycle in $\mathcal{C}(\mathcal{P},\rho)$. However $\Delta - \gamma$ is not a boundary or homologous to another cycle since $\Delta$ is not the face of any higher simplex in $\mathcal{C}(\mathcal{P},\rho)$, and so we introduce a new non-trivial $k$-cycle.

Thus it suffices to show that $\partial \Delta$ is indeed a boundary. Consider the natural map from the simplicial chains to the singular chains $\iota_k :C_k(\mathcal{C}(\mathcal{P},\rho_-)) \to C_k(B_{\rho_-}(\mathcal{P}))$. For sufficiently small $\rho_-$, by the Nerve Lemma the induced map on homology is an isomorphism $h_k : H_k(\mathcal{C}(\mathcal{P},\rho_-)) \cong H_k(B_{\rho_-}(\mathcal{P}))$.

If $\rho_-$ is sufficiently close to $\rho$ then given that ${A}^{\varphi}_\psi(\mathcal{Y}) \subset B_\rho(\mathcal{P})$ we also have ${A}^{\varphi}_{2\psi}(\mathcal{Y}) \subset B_{\rho_-}(\mathcal{P})$. Our $\psi$ is constructed in order that $\iota_{k-1}(\partial \Delta) \in C_{k-1}({A}^{\varphi}_{2\psi}) \subset C_{k-1}(B_{\rho_-}(\mathcal{P}))$. ${A}^{\varphi}_{2\psi}$ is homotopic to a $d$-dimensional annulus sliced by a $d-1$-hyperplane and so $H_{k-1}({A}^{\varphi}_{2\psi}) = 0$, and thus $h_{k-1}(\partial \Delta)=0$
\end{proof}

\begin{defn}($\Theta$-like-cycle)

A critical point $c(\mathcal{Y})$ which satisfies the conditions of the above lemma is a $\Theta$-like-cycle at scale $\psi(\mathcal{Y})$.
 \label{defn: Theta-like-cycle}
\end{defn}

\begin{figure}[h]
\centering
	\begin{subfigure}[b]{0.4\linewidth}
    	\includegraphics[width=\linewidth]{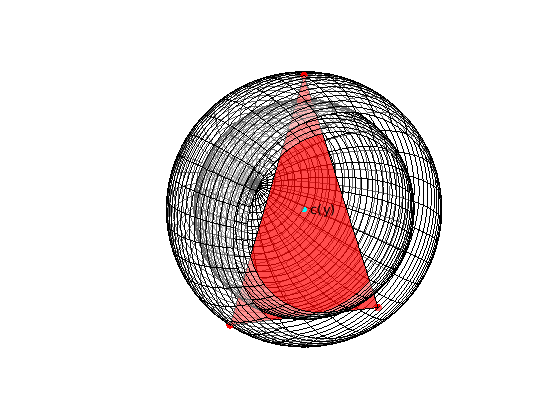}
    	\caption{Large Epsilon}
	\end{subfigure}
	\begin{subfigure}[b]{0.4\linewidth}
    	\includegraphics[width=\linewidth]{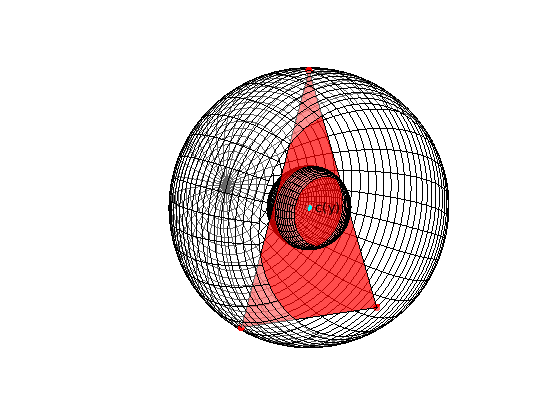}
    	\caption{Small Epsilon}
	\end{subfigure}
    \caption{A sketch of an index $2$ critical point near the boundary inducing a $\Theta$-like-cycle. The annulus ${A}_\varepsilon$ is bounded by the spheres and we note that for sufficiently small $\varepsilon$ the boundary of the $2$-simplex associated to the critical point is contained in the annulus.}
    \label{fig:Theta-Cycle}
\end{figure}

\begin{figure}[h]
\centering
	\begin{subfigure}[b]{0.4\linewidth}
    	\includegraphics[width=\linewidth]{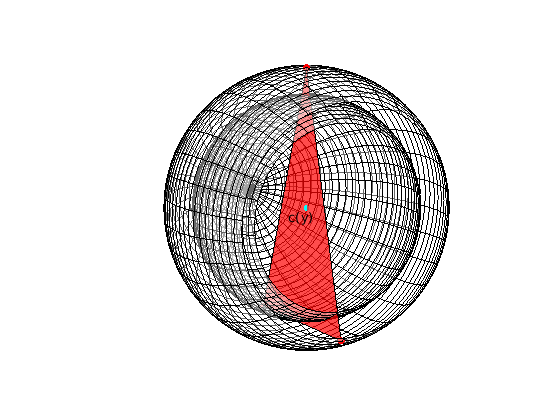}
    	\caption{Large Epsilon}
	\end{subfigure}
	\begin{subfigure}[b]{0.4\linewidth}
    	\includegraphics[width=\linewidth]{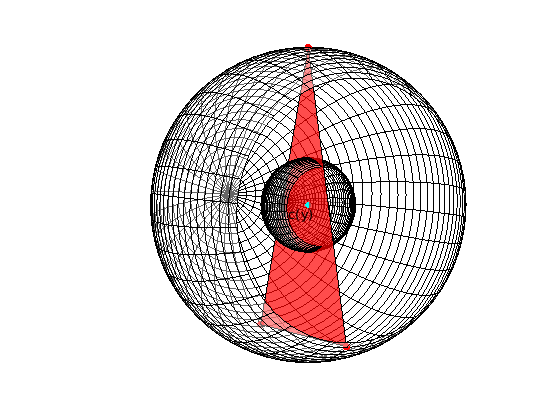}
    	\caption{Small Epsilon}
	\end{subfigure}
    \caption{A sketch of an index $2$ critical point near the boundary that does not induce a $\Theta$-like-cycle. The annulus ${A}_\varepsilon$ does not contain the boundary of the $2$-simplex associated to the critical point for $\varepsilon > 0$.}
    \label{fig:Non-Theta-Cycle}
\end{figure}

Note that if $\partial \Delta (\mathcal{Y})$ traverses the section of the annulus ${A}_\varepsilon$ that has been sliced away by the boundary for all $\varepsilon>0$ then $\psi(\mathcal{Y}) = 0$. We do not want to count the critical points induced by such $\mathcal{Y}$ since they do not introduce a non-trivial homological cycle. We sketch an example of a critical point inducing a $\Theta$-like-cycle in Figure \ref{fig:Theta-Cycle}, and in contrast also sketch an example of a critical point not inducing a $\Theta$-like-cycle in Figure \ref{fig:Non-Theta-Cycle}.

\subsection{$\Theta$-like-cycle Lower Bound Computation}

Let us modify our indicator functions from the previous section to count $\Theta$-like-cycles with critical values in the range $[r_1,r)$. We count cycles whose critical point lies at a distance to the boundary in the range $[\delta,2\delta]$, for suitably chosen $\delta$.
\begin{enumerate}
\item $h(\mathcal{Y}) = \mathbf{1}\{0 \in Grad(\mathcal{Y})\}$

\item $h_{r_1,r}(\mathcal{Y}) = h(\mathcal{Y})\mathbf{1}\{r_1 < \rho(\mathcal{Y}) \leq r\}$

\item $h_{r}^{\varepsilon, \delta}(\mathcal{Y}) = h_{r_1,r}(\mathcal{Y})\mathbf{1}\{\psi(\mathcal{Y},\varphi)\geq \varepsilon\}\mathbf{1}\{\delta \leq \text{dist}(c(\mathcal{Y}),\partial M) \leq 2\delta\}$

\item $g_{r}^{\varepsilon, \delta}(\mathcal{Y},\mathcal{P}) = h_{r}^{\varepsilon, \delta}(\mathcal{Y})\mathbf{1}\{ B_{r_2}(c(\mathcal{Y})) \cap (\mathcal{P}\setminus \mathcal{Y}) = \emptyset \} \mathbf{1}\{ {A}^{(\varphi)}_\varepsilon \subset B_{\rho(\mathcal{Y})}(\mathcal{P}) \}$

\end{enumerate}

Lemma \ref{lem: Theta-like-cycles induce homology} shows that these critical points introduce erroneous $k$ cycles.  Our $\delta$ will be chosen to be $\delta \sim \frac{r}{\log n}$. If $\mathcal{P}_n$ is our Poisson process on $M$ then we count the $\Theta$-like-cycles with the following sum of indicator functions:

$$ \Theta_k^{\varepsilon, \partial M}(r_1,r) = \sum_{|\mathcal{Y}| = k+1} g_{r}^{\varepsilon, \delta}(\mathcal{Y},\mathcal{P}_n)$$

\subsubsection{Partial $\varepsilon$-Annulus Coverage}

Let us denote the probability of covering the partial $\varepsilon$-annulus of a critical point $c(\mathbf{y})$ conditioned on the critical point having no points but $\mathbf{y}$ in an $r_2$ neighbourhood as:

$$ p_{\varepsilon,\varphi}(\mathbf{y}) = \mathbb{P}({A}^{(\varphi)}_\varepsilon(\mathcal{Y}')\subset B_{\rho(\mathcal{Y}')}(\mathcal{P}_n))\ | \ \mathcal{Y}' = \mathbf{y}, \ \mathcal{P}_n \cap B_{r_2}(c(\mathcal{Y}')) = \mathcal{Y}')$$

As in previous sections we use Palm Theory (Theorem \ref{thm:PalmTheory1}) applied to the Poisson process on $M$ to attain an integral expression:

\begin{align*} \mathbb{E}[|\Theta_k^{\varepsilon, \partial M}|] &=
\frac{n^{k+1}}{(k+1)!} \mathbb{E}[g_{r}^{\varepsilon, \delta}(\mathcal{Y}',\mathcal{Y}' \cup \mathcal{P}_n)] \\
&= \frac{n^{k+1}}{(k+1)!} \int_{DM^{k+1}} h_{r}^{\varepsilon, \delta}(\mathbf{y})p_{\varepsilon,\varphi}(\mathbf{y})e^{-n\text{Vol}(B_{r_2}(\mathbf{y})\cap M)} |\text{dvol}_g(\mathbf{y})|
\end{align*}

There is a tradeoff in the choice of $\varphi$ and $\delta$ for maximising our lower bound for the expected count of $\Theta$-like-cycles. The expectation depends on $\varphi$ in the support of the functions $\mathbf{1}\{\psi(\mathcal{Y},\varphi)\geq \varepsilon\}$, $\mathbf{1}\{ {A}^{(\varphi)}_\varepsilon \subset B_{\rho(\mathcal{Y})}(\mathcal{P}) \}$. Increasing $\varphi$, augments the support of $\mathbf{1}\{\psi(\mathcal{Y},\varphi)\geq \varepsilon\}$ and diminishes the support of $\mathbf{1}\{ {A}^{(\varphi)}_\varepsilon \subset B_{\rho(\mathcal{Y})}(\mathcal{P}) \}$. The term $\mathbf{1}\{\psi(\mathcal{Y},\varphi)\geq \varepsilon\}$ determines that the points $\mathcal{Y}$ must lie in a hyperplane approximately tangent to the boundary, and the influence of this term on the lower bound is computed in the following section. Meanwhile the term $\mathbf{1}\{ {A}^{(\varphi)}_\varepsilon \subset B_{\rho(\mathcal{Y})}(\mathcal{P}) \}$ insists that the partial annulus is covered by the point process and is tracked by $p_{\varepsilon,\varphi}(\mathbf{y})$.

With a judicious choice of $\varphi \sim (\log n)^{-1}$ and $\delta \sim \frac{r}{\log n}$ we can show that the term $p_{\varepsilon,\varphi}(\mathbf{y}) \to 1$ uniformly as $\Lambda \to \infty$ for fixed $\varepsilon$.

\begin{figure}[ht]
\centering
	\begin{subfigure}[b]{0.3\linewidth}
    	\includegraphics[width=\linewidth]{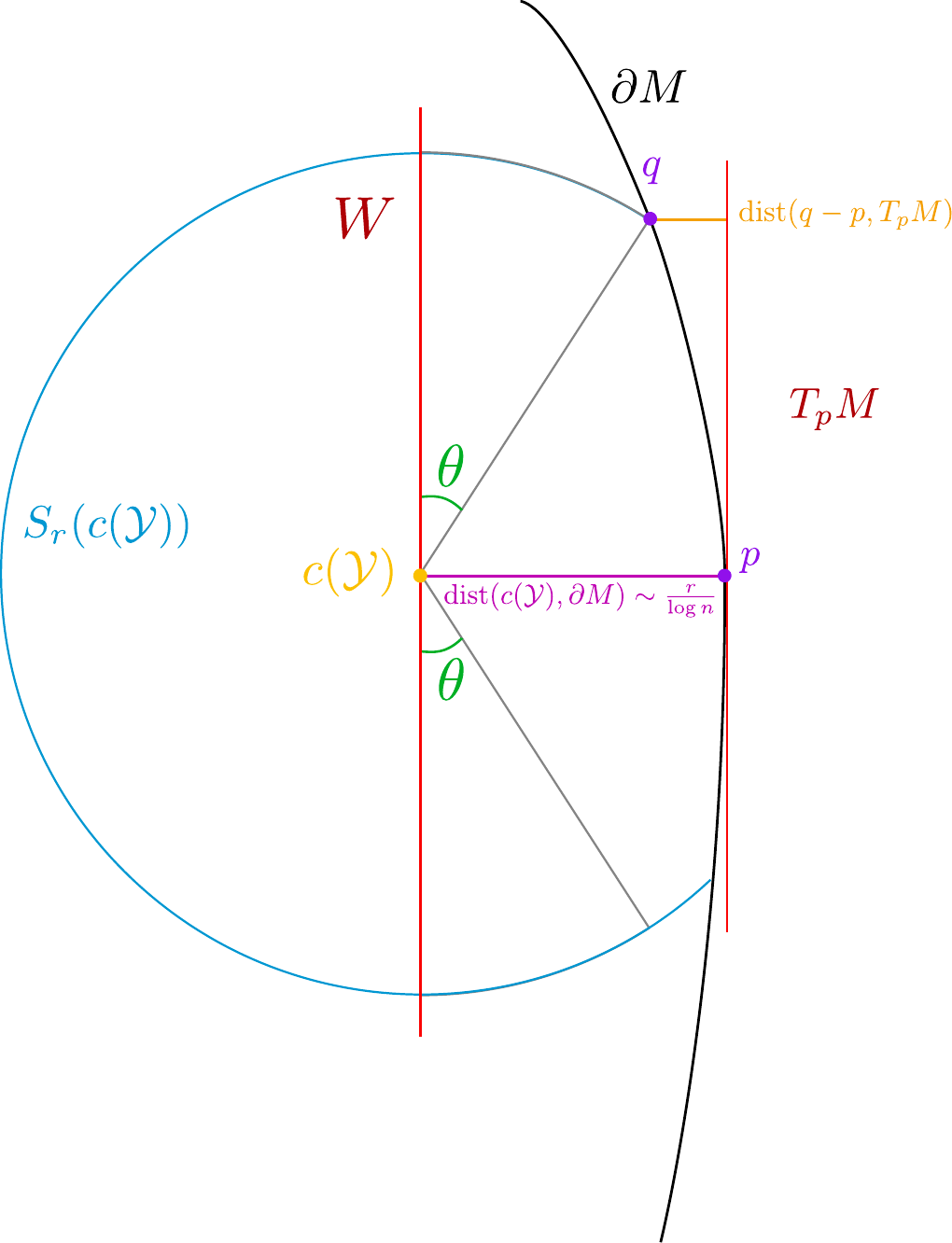}
    	\caption{}
    	\label{fig:Coverage_1}
	\end{subfigure}
	\begin{subfigure}[b]{0.3\linewidth}
    	\includegraphics[width=\linewidth]{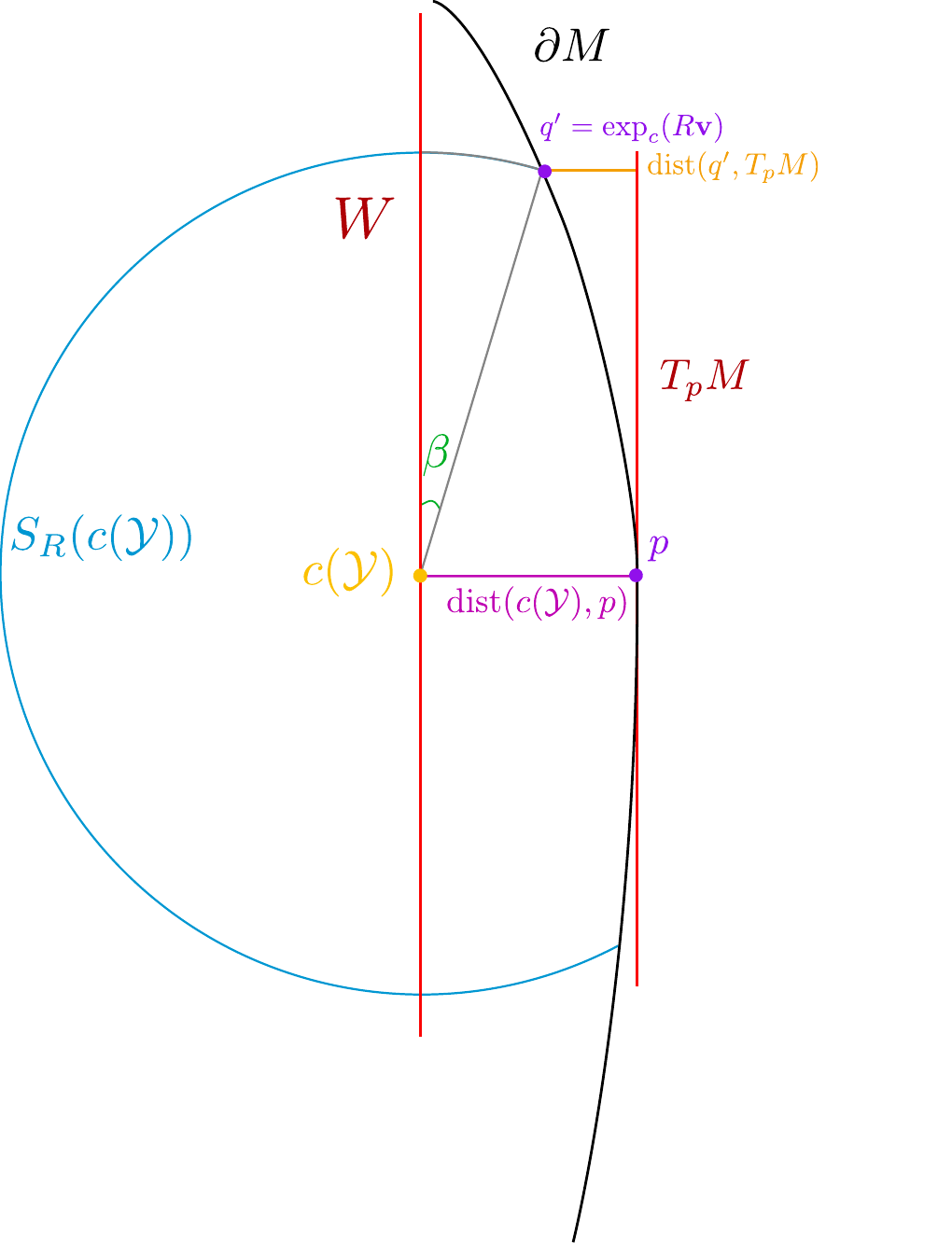}
    	\caption{}
    	\label{fig:Coverage_2}
	\end{subfigure}
	\begin{subfigure}[b]{0.25\linewidth}
    	\includegraphics[width=\linewidth]{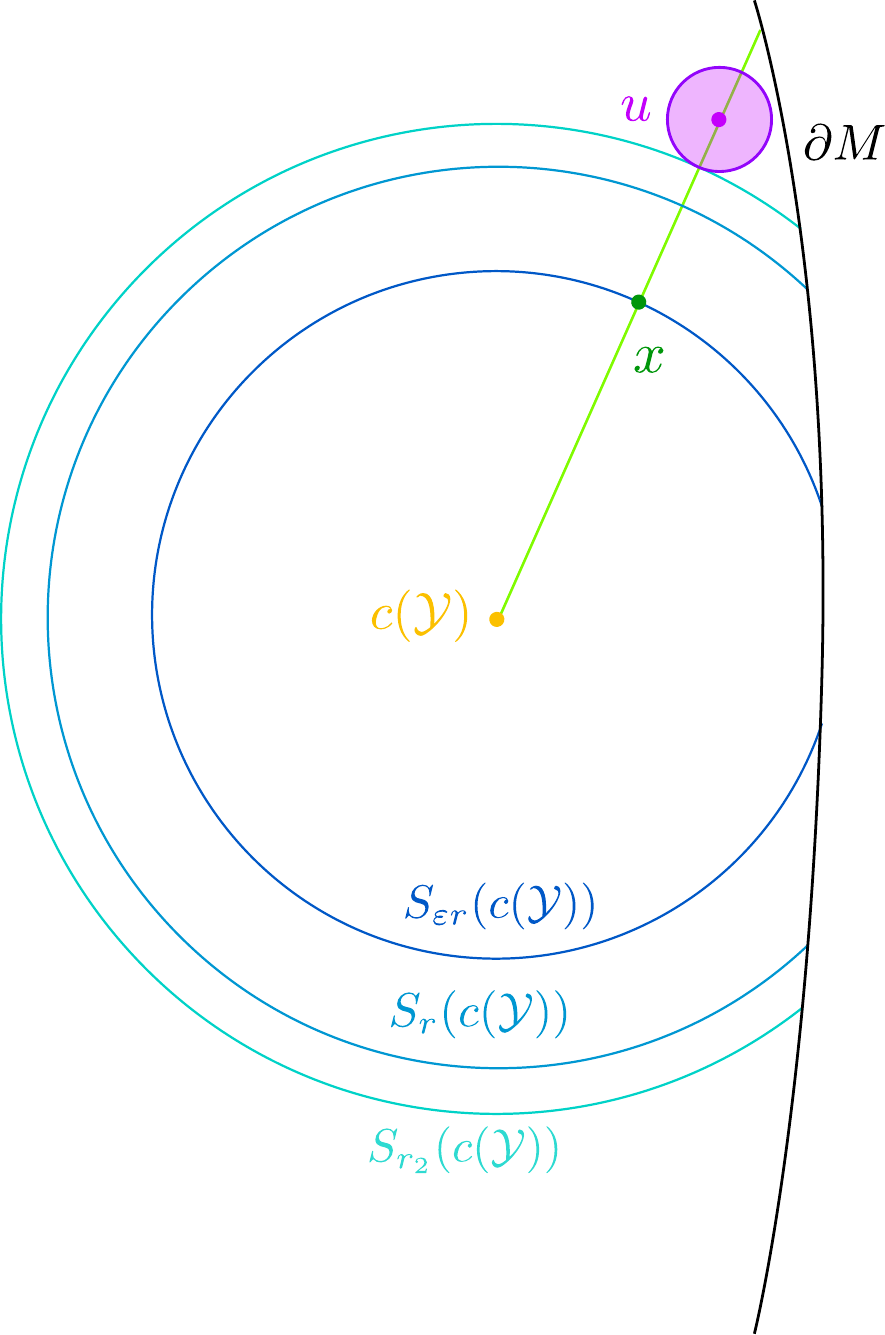}
    	\caption{}
    	\label{fig:Coverage_3}
	\end{subfigure}
    \caption{(a) We see that for $\theta = C (\log n)^{-1}$, the boundary does not meet the partial annulus $A_\varepsilon^{(\theta)}$. (b) We see that for $\beta = C (\log n)^{-1}$, the boundary does not meet $S_{R}(c(\mathcal{Y}))$. (c) A point $x$ on the inner most part of the partial annulus can be covered by any point sampled in the ball centred at $u$.}
    \label{fig:Partial-Annulus-Coverage}
\end{figure}

\begin{lem}
Let $c$ be a critical point at distance $\delta$ from the boundary with $\delta \in [\frac{r}{ \log n},\frac{2r}{ \log n}]$, then for $\varphi \sim \frac{1}{\log n}$ the partial annulus at $c$ is contained in the deformed annulus at $c$, then $p_{\varepsilon,\varphi}(\mathbf{y}) \to 1$ uniformly as $\Lambda \to \infty$ for fixed $\varepsilon$. \label{lem: Partial Annulus Coverage}
\end{lem}
\begin{proof}

Let $p\in \partial M$ be the unique point on the boundary closest to $c$, and let $\vec{n} \in T_cM$ be the \textit{normal} vector at such that $\textrm{exp}_c(\vec{n}) = p$. This induces the \textit{tangent} hyperplane $W$ at $c$; ($W \subset T_cM$ such that $\vec{n}\cdot W  =0$).

We shall first find $\theta$ such that $A_\varepsilon^{(\theta)}(c)$ meets the boundary. See Figure \ref{fig:Coverage_1}.
Let $ U  = \{ \vec{u} \in T_cM \ :\ \textrm{exp}_c(\vec{u}) \in \partial M \textrm{ and } \|u\| = r \}  = \textrm{exp}_c^{-1}(S_r(c)\cap \partial M)$, then we define $\theta_\vec{u}$ to be twice the acute angle formed between $\vec{u}$ and $W$, and let $\theta = \min_U \theta_\vec{u}$. Thus for $p$ the closest point on the boundary to $c$ and $q = \textrm{exp}_c(\vec{u})$ we attain $d(q-p,T_pM) \leq \frac{\|q-p\|^2}{2\tau_M} \sim r^2 = o(\frac{r}{\log n})$. Hence it follows that $\tan \frac{\theta}{2} \sim \frac{1}{\log n}$ and thus $\theta \sim \frac{1}{ \log n}$. Let $C_1 , C_2$ be constants such that $ \frac{C_1}{\log n } \leq \theta \leq \frac{C_2}{\log n}$.



Let $\vec{v}$ be a unit vector making acute angle $\beta$ with the \textit{tangent} hyperplane $W$ normal to $\vec{n}$. Let the line leaving $c$ in direction $\vec{v}$ meet the boundary at point $q'$, that is $q' = \text{exp}_c(R\vec{v}) \in \partial M$ for some scalar length $R$. We wish to bound $\beta$ such that $R \geq 2r$, see Figure \ref{fig:Coverage_2}.

Assume $\rho_M(p,q') < 3r$ else we are done by the triangle inequality. Using Theorem \ref{thm:Federer} we see that since  $\textrm{dist}(W,R\vec{v}) \geq \delta / 2$. Using basic trigonometry we observe that $\textrm{dist}(W, R\vec{v}) = R \sin \beta \leq R \beta$. Thus for $\beta  = \frac{\varphi}{C}$ and $C \geq 2 C_2$ we have $R \geq \frac{\delta}{2\beta} \geq \frac{rC}{C_2} \geq 2r$ as desired.

Next let us show that for all $x \in {A}^{(\frac{\theta}{C})}_\varepsilon(c)$ the volume $\textrm{Vol}( B_{r(1-\varepsilon/10)}(x)\setminus B_{r_2}(c))\geq \frac{1}{1000} \epsilon^d r^d$. It suffices to show this property for $x = \text{exp}_c(\varepsilon r \vec{v})$ since the volume $\textrm{Vol}( B_{r(1-\varepsilon/10)}(x)\setminus B_{r_2}(c))$ increases in the radial direction from $c$, see Figure \ref{fig:Coverage_3}. Let $t = \frac{1}{2}((1+\frac{9\varepsilon}{10})r -r_2)$ then the volume of the maximal radius ball centred at $u = \text{exp}_c(t\vec{v})$ contained in $ B_{r(1-\varepsilon/10)}(x)\setminus B_{r_2}(c) $ witnesses that $\textrm{Vol}( B_{r(1-\varepsilon/10)}(x)\setminus B_{r_2}(c)) > \frac{1}{1000}\varepsilon^d r^d$. 

This follows since we have shown $R \geq 2r$ and so $u=\text{exp}_c(t\vec{v}) \in M$, and moreover $r_2 \to r$ means for sufficiently large $n$ the maximal radius is at least $\frac{\varepsilon}{20}$. 
Observe that $\text{Vol}({A}^{(\varphi)}_\varepsilon) \sim \rho^d (1-\varepsilon^d)$ and so there is a constant $C$ dependent only on the metric $g$ such that there is an $\frac{\varepsilon\rho}{10}$-net of ${A}^{(\varphi)}_\varepsilon$, $\mathcal{S}$ with $|\mathcal{S}|\leq C \frac{1-\varepsilon^d}{\varepsilon^d}$. It is clear that ${A}^{(\varphi)}_\varepsilon \subset B_\rho(\mathcal{P}_n)$ if for all $s \in \mathcal{S}$ we have $\mathcal{P}_n \cap B_{\rho(1-\varepsilon/10)}(s) \neq \emptyset$. Conditioning on the event that $\{ \mathcal{P}_n \cap B_{r_2} = \emptyset\}$ we may thus bound $p_{\varepsilon,\varphi}$ from below:

$$ p_{\varepsilon,\varphi} \geq 1 - C \max_{s\in \mathcal{S}} e^{-n \text{Vol}(B_{\rho(1-\varepsilon/10)}\setminus B_{r_2})}$$

Hence given we have shown that $\textrm{Vol}( B_{r(1-\varepsilon/10)}(x)\setminus B_{r_2}(c))\geq \frac{1}{1000} \epsilon^d r^d$ for all $x \in {A}^{(\frac{\theta}{C})}_\varepsilon(c)$ we have that $ p_{\varepsilon,\varphi} \geq 1- Ce^{-n\varepsilon^dr^d/1000} \to 1$ uniformly as $\Lambda \to \infty$ for fixed $\varepsilon$. 

\end{proof}

Since $p_{\varepsilon,\varphi}(\mathbf{y}) \to 1$ uniformly as $\Lambda \to \infty$ for fixed $\varepsilon$, we may remove the term $p_{\varepsilon,\varphi}$ from the integral and replace it by a constant.

\begin{align*} \mathbb{E}[|\Theta_k^{\varepsilon, \partial M}|] &\geq C \frac{n^{k+1}}{(k+1)!} \int_{DM^{k+1}} h_{r}^{\varepsilon, \delta}(\mathbf{y})e^{-n\text{Vol}(B_{r_2}(c(\mathbf{y}))\cap M)} |\text{dvol}_g(\mathbf{y})|
\end{align*}

The next term we shall attack in the lower bound is $e^{-n\text{Vol}(B_{r_2}(c(\mathbf{y}))\cap M)}$. Using the result of Lemma \ref{prop:Half-BallVolumes} and the Blaschke Petkantschin formula we attain the lower bound:



\begin{align*} \mathbb{E}[|\Theta_k^{\varepsilon, \partial M}|] &\geq C \frac{n^{k+1}}{(k+1)!} \int_{DM^{k+1}} h_{r}^{\varepsilon, {\partial M_{\delta}}}(\mathbf{y})e^{-n\alpha\text{Vol}(B_{r_2}(c(\mathbf{y})))} |\text{dvol}_g(\mathbf{y})| \\
&= C \frac{n^{k+1}}{(k+1)!}
\int_{\partial M_{[\delta,2\delta]}} |\text{dvol}_g(c)|
\int_{r_1}^{r} du \  u^{dk-1} \int_{Gr(k,T_cM)} d\mu_{k,d}(V) \\
& \times  \prod_{i=1}^{k} \left( \int_{\mathbb{S}_1(V)} \sqrt{|\text{det}(g_{\text{exp}_c(uw_i)})|} |\text{dvol}_{\mathbb{S}_1(V)}(w_i) |\right)
\mathbf{\Upsilon}^{d-k}_1(\mathbf{w})f(\text{exp}_c(u\mathbf{w}))
\end{align*}

Where $f(\mathbf{y}) = h_{r}^{\varepsilon, \partial M_{\delta}}(\mathbf{y})e^{-n\alpha\text{Vol}(B_{r_2}(c(\mathbf{y})))}$ and $\alpha = \frac{1}{2} + O(r, \frac{\delta}{r}) = \frac{1}{2} + O((\log n )^{-1})$. We can bound from below the components of this integral using the fact that the manifold is compact. In particular:

$$ e^{-n\alpha\text{Vol}(B_{r_2}(c(\mathbf{y})))} \geq e^{-\alpha\Lambda_{r_2}}(1+ s_\text{min}\Lambda_{r_2}r^2)$$

\begin{align*} \mathbb{E}[|\Theta_k^{\varepsilon, \partial M}|] &\geq 
C \frac{n^{k+1}}{(k+1)!} e^{-\alpha\Lambda_{r_2}}(1+ s_\text{min}\Lambda_{r_2}r^2)
\int_{\partial M_{[\delta,2\delta]}} |\text{dvol}_g(c)|
\int_{r_1}^{r} du \  u^{dk-1} \int_{Gr(k,T_cM)} d\mu_{k,d}(V) \\
& \times  \prod_{i=1}^{k} \left( \int_{\mathbb{S}_1(V)} \sqrt{|\text{det}(g_{\text{exp}_c(uw_i)})|} |\text{dvol}_{\mathbb{S}_1(V)}(w_i) |\right)
\mathbf{\Upsilon}^{d-k}_1(\mathbf{w})h_{r}^{\varepsilon, \delta}(\text{exp}_c(u\mathbf{w}))
\end{align*}

\subsubsection{Grassmannian Volume}

We require subtle analysis to bound from below the contribution of the integrals over the Grassmannian. Let us define:

$$ D^\varepsilon_{k,\varphi} = \int_{Gr(k,T_cM)} d\mu_{k,d}(V)  \times  \prod_{i=1}^{k} \left( \int_{\mathbb{S}_1(V)} \sqrt{|\text{det}(g_{\text{exp}_c(uw_i)})|} |\text{dvol}_{\mathbb{S}_1(V)}(w_i) |\right)
\mathbf{\Upsilon}^{d-k}_1(\mathbf{w})h_{r}^{\varepsilon, \delta}(\text{exp}_c(u\mathbf{w}))$$

Taking sufficiently small $r$ we may assume the determinant term is arbitrarily close to $1$. Moreover given $\psi(\mathbf{y},\varphi)\geq \varepsilon$ it is clear the volume $\Upsilon(\mathbf{w})$ is bounded below by a term of order $\text{dist}({A}^{(\varphi)}_\varepsilon,c(\mathcal{Y}))^k = \Omega(\varepsilon^k)$.

\input{technical2}




We can calculate a lower bound using this estimate and with $r_1 = r(1-\frac{\xi^2}{2c_g^2}) $, $r_2 = r(1+\xi)$: 

\begin{align*} \mathbb{E}[|\Theta_k^{\varepsilon, \partial M}|] &\geq 
C (1-c_Rr^2)(1+ s_\text{min}\Lambda_{r_2}r^2) D^\varepsilon_{k,\varphi}  n^{k+1} e^{-\alpha\Lambda_{r_2}}
\int_{\partial M_{[\delta,2\delta]}} |\text{dvol}_g(c)|
\int_{r_1}^{r} du \  u^{dk-1} \\
& \geq 
C D^\varepsilon_{k,\varphi}   n^{k+1}e^{-\alpha\Lambda_{r_2}} \delta r^{dk}
\int_{r_1/r}^{1} s^{dk} \ ds \\
& \geq 
C D^\varepsilon_{k,\varphi}   n\Lambda^ke^{-\alpha\Lambda_{r_2}} \delta \xi^2
\end{align*}

For $\xi = \Lambda^{-1}$ we observe that asymptotically $e^{-\alpha\Lambda_{r_2}}= \Omega(e^{-\alpha \Lambda})$ and so: 

\begin{align*} \mathbb{E}[|\Theta_k^{\varepsilon, \partial M}|] &\geq 
C D^\varepsilon_{k,\varphi}   n\Lambda^{k-2}e^{-\alpha\Lambda_{r_2}} \delta   \\
&\geq C n\Lambda^{k-2}e^{-\alpha\Lambda} r (\log n)^{-(k+1)}
\end{align*}

Thus we see that $\mathbb{E}[|\Theta_k^{\varepsilon, \partial M}|] = \Omega( n\Lambda^{k-2}e^{-\alpha\Lambda}r (\log n)^{-(k+1)})$

In order to complete the argument for the lower threshold we will need to verify that the second-moment arguments in \cite{Bobrowski2017} carry over to our modified notion of a $\Theta$-like-cycle, and so w.h.p. the number of $\Theta$-like-cycles is bounded below by the same regime as the expected number of $\Theta$-like-cycles.

%% file: technical2.tex
\begin{figure}[h]
\centering
	\begin{subfigure}[b]{0.3\linewidth}
	    \centering
    	\includegraphics[width=\linewidth]{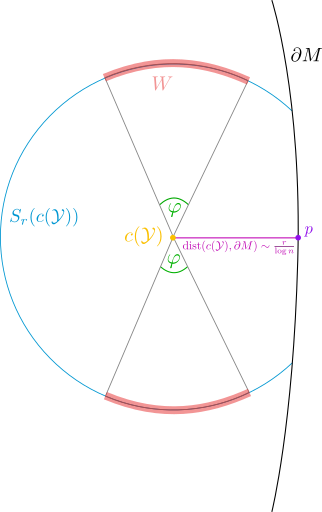}
    	\caption{A cartoon of a critical point near the boundary}
    	\label{fig: Boundary Critical Point Cartoon}
	\end{subfigure}
	\hfill
	\begin{subfigure}[b]{0.6\linewidth}
	    \centering
    	\includegraphics[width=\linewidth]{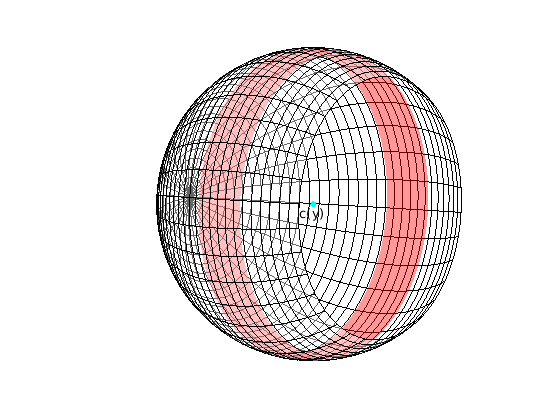}
    	\caption{A critical point near the boundary in a $3$ manifold with boundary}
    	\label{fig: Boundary Critical Point 3-manifold}
	\end{subfigure}
\caption{In Figure \ref{fig: Boundary Critical Point Cartoon} we sketch a critical point $c(\mathcal{Y})$ close to the boundary to highlight the important quantities in computing a lower bound for the number of $\Theta$-like-cycles. We highlight a subset $W \subset S_r(c(\mathcal{Y}))$ such that if every point of $\mathcal{Y}$ lies in $W$ the simplex introduced by the critical point $c$ is a $k$-dimensional lid on the sliced-annulus centred around $c$ and so induces a $\Theta$-like cycle. The corresponding subset $\mathcal{W} \subset \Gr(k,T_cM)$ consists of those $k$-planes whose intersection with the sphere lies in  $W \subset S_r(c(\mathcal{Y}))$. In Figure \ref{fig: Boundary Critical Point 3-manifold} we give a concrete example for $M$ a $3$-manifold with boundary. If the points $\mathcal{Y}$ inducing an index $2$ critical point at $c$ lie in the highlighted subset they induce a $\Theta$-like-cycle for some $\varepsilon$ as in the case of Figure \ref{fig:Theta-Cycle}.}
\label{fig:Grassmannian-sketch}
\end{figure}

Let us calculate a lower bound for the volume of the subspace of the Grassmannian $\Gr(k,d)$ for which the points inducing our centre must lie in order that  $c(\mathcal{Y})$ induces a $\Theta$-like cycle of order $k$, see Figure \ref{fig:Grassmannian-sketch}. Since we are using the invariant measure $\mu_{k,d}$ on $\Gr(k,d)$ as in the Blaschke-Petkantschin Formula we can appeal to Theorem 13.1.5 from \textit{Stochastic and Integral Geometry} \cite{schneider_stochastic_2008}, and compute a bound for this volume using the action of $SO(d)$ on $\Gr(k,d)$.

For sufficiently small $r$ there is a unique closest point to $c$ on the boundary $p\in \partial M$. Let $\vec{n} \in T_cM$ be such that $\textrm{exp}_c(\vec{n}) = p$ and let $e_d = \vec{n}/\| \vec{n}\|$ be the unit vector \textit{normal} to the boundary at $c$. 


Define $W$ to be a subset of the sphere in the tangent space at $c$ as sketched in Figure \ref{fig:Grassmannian-sketch}:
$$W := \{\vec{v} \in T_cM : \|v\| = 1, |\langle \vec{v},\vec{n} \rangle| \leq \sin (\varphi / 2) \} \subset T_{c(\mathcal{Y})}M$$

This induces a subset of the Grassmannian:

$$ \mathcal{W} := \{V \in \Gr(k,d) : V = \langle w_1,...,w_k\rangle, w_i \in W  \} \subset \Gr(k,d)$$

For our suitably chosen $\varphi$, a critical point induced by points lying in $W$ will introduce a simplex which lies approximately tangential to the boundary and thus will induce a $\Theta$-like-cycle. In order to calculate the volume of the associated subset of the Grassmannian $\mathcal{W}$ we shall use Theorem \ref{thm: Compact group Volume}.

\begin{thm}(Theorem 13.1.5)\cite{schneider_stochastic_2008}
Let $G$ be a compact group operating continuously and transitively on a Hausdorff space $E$, and suppose $G$ and $E$ have countable bases. Let $\nu$ be a Haar measure on $G$ with $\nu(G) =1 $. Then there exists a unique $G$-invariant Borel measure $\rho$ on $E$ with $\rho(E)=1$ defined by: 
$$ \rho(B) = \nu(\{g\in G \ :\  g\cdot x_0 \in B\}), \ \ B\in \mathcal{B}(E)$$ for arbitrary $x_0 \in E$.
\label{thm: Compact group Volume}
\end{thm}

We observe that the compact group $SO(d)$ acts continuously and transitively on $\Gr(k,d)$. We can parametrise $SO(d)$ with ${d \choose 2}$ angles $\{\phi_{i,j}\}_{1\leq i \leq j\leq d-1}$ \cite{hurwitz_uber_1897,diaconis_._2015}. 

Let us briefly recall this parametrization. Denote the augmentation of vectors and matrices via a superscript $a$:
$$R^a = \begin{bmatrix}
R & \vec{0} \\
\vec{0}^T & 1
\end{bmatrix} , \ \  (v_1,...,v_l)^a =  (v_1,...,v_l,0)$$

Recall the hyperspherical coordinate system of the unit $d$-sphere $S^d$ given by $d$ angles $\theta_1,...,\theta_d$, $\Sigma_d : \mathbb{R}^d \to S^d$ . Let $\pi_i$ denote the projection onto the $i^\textrm{th}$ coordinate.

$$ \pi_i(\Sigma_d(\theta_1,...,\theta_d)) = \begin{cases} 
\sin \theta_1 ... \sin \theta_d & \textrm{if } i = 1 \\
\cos \theta_{i-1} \sin \theta_i ... \sin \theta_d & \textrm{if } 2 \leq i \leq d+1
\end{cases}$$

We have an orthonormal basis $\{\vec{u}_{d,1},...,\vec{u}_{d,d}, \Sigma_d(\theta_1,...,\theta_d) \}$ where $\vec{u}_{d,i}$ is the unit vector in the direction $\frac{\partial \Sigma_d(\theta_1,...,\theta_d)}{\partial \theta_i}$. Let the rotation $R_{d+1}(\theta_1,...,\theta_d) \in SO(d+1)$ be the matrix:

$$ R_{d+1}(\theta_1,...,\theta_d) = \begin{bmatrix}
\vec{u}_{d,1} & ... & \vec{u}_{d,d} & \Sigma_d(\theta_1,...,\theta_d)
\end{bmatrix} $$ 

Observe that $R_{d+1}$ maps the standard basis vector $\vec{e}_{d+1}$ to the point $\Sigma_d(\theta_1,...,\theta_d)\in S^d$.

We inductively construct $\Phi_d \in SO(d)$ from ${d \choose 2}$ angles $\{\phi_{i,j}\}_{1\leq i \leq j\leq d-1}$ as follows:

$$ \Phi_2(\phi_{1,1}) = \begin{bmatrix} \cos \phi_{1,1} & \sin \phi_{1,1} \\
-\sin \phi_{1,1} & \cos \phi_{1,1}
\end{bmatrix}$$
$$\Phi_d(\{\phi_{i,j}\}_{1\leq i \leq j\leq d-1}) = R_{d}(\phi_{1,d-1},...,\phi_{d-1,d-1})\left(\Phi_{d-1}(\{\phi_{i,j}\}_{1\leq i \leq j\leq d-2})\right)^a $$
$$ \phi_{1,j}\in [0,2\pi] , \ \phi_{i,j} \in [0,\pi] \textrm{ for } 2\leq i \leq d-1, 1\leq j \leq d-1 $$

The Haar measure on $SO(d)$ is then given up to normalising constant by: $$d\nu = \left(\prod_{1\leq i\leq j \leq d-1} \sin^{i-1}\phi_{i,j}\right) d\phi_{1,1} ... d\phi_{d-1,d-1}$$

Let $P_k$ denote the projection onto the $k$-plane $x_0$. Observe that if $\| P_k \Phi_d e_d \| \leq \sin \frac{\varphi}{2}$ then $\Phi_d x_0 \in \mathcal{W}$. 
Recall $x_0 = \langle e_1,...,e_k \rangle$ where $e_i$ are standard basis vectors. Using the inductive construction we observe that $\| P_k \Phi_d e_d \| = \| P_k \Sigma_{d-1}(\phi_{1,d-1},...,\phi_{d-1,d-1}) \|\leq  \sin \phi_{k,d-1}$ since the first $k$ coordinates of $\Sigma_{d-1}(\phi_{1,d-1},...,\phi_{d-1,d-1})$ contain a factor of $\sin \phi_{k,d-1}$.

Thus we attain the following estimate for $\mu_{k,d}(\mathcal{W})$:

\begin{align*}
\mu_{k,d}(\mathcal{W}) &= \int_{\Phi\in SO(d)} \vec{1}\{\Phi\cdot x_0 \in \mathcal{W} \} d\nu \sim \int_{\phi_{i,j}} \vec{1}\{\Phi\cdot x_0 \in \mathcal{W} \} \left(\prod_{1\leq i\leq j \leq d-1} \sin^{i-1}\phi_{i,j}\right) d\phi_{1,1} ... d\phi_{d-1,d-1} \\
&\geq C \int_{0}^{\varphi/2} \sin ^{k-1} \phi_{k,d-1} d\phi_{k,d-1} = \Omega(\varphi^{k})
\end{align*}

Thus we have a lower bound estimate for the Grassmannian volume in terms of the angle $\varphi$.

Given a critical point $c$ induced by $\mathcal{Y}\subset A^{(\varphi)}_\varepsilon(c)$ and associated critical simplex $\Delta$, for sufficiently small $\varepsilon$ we have $\partial \Delta \subset A^{(\varphi)}_\varepsilon(c) $ and thus $c$ induces a $\Theta$-like cycle.

We can compute a lower bound for the support of the term $\mathbf{1}\{\psi(\vec{y},\varphi)\geq \varepsilon\}$ by observing that:
$$ \mathbf{1}\{\psi(\vec{y},\varphi)\geq \varepsilon \} \geq \mathbf{1}\{\phi(\vec{y}) \geq \varepsilon \} \mathbf{1}\{\vec{w} \in W \subset T_{c(\vec{y})}M\}$$

The support of the term  $\mathbf{1}\{\vec{w} \in W \subset T_{c(\vec{y})}M\}$ may be bounded below using the Grassmannian volume calculation above and as observed in \cite{Bobrowski2017} the support of the term $\mathbf{1}\{\phi(\vec{y}) \geq \varepsilon \}$ is bounded below by a constant for fixed $\varepsilon$.
Hence we conclude that for sufficiently small $\varepsilon$ we have $D^\varepsilon_{k,\varphi} = \Omega((\log n)^{-k})$.

%% file: second-moment-calculations.tex
\label{section: Second Moments}
We have attained a lower bound for the expected number of $\Theta$-like-cycles near the boundary of the manifold. Let us now seek to show that w.h.p. the number of $\Theta$-like-cycles is bounded below by this regime. Let us make the abbreviation $\Theta_k^{\varepsilon, \partial M} = \Theta_k^{\varepsilon, \partial M}(r_1,r)$. Using Chebyshev's inequality we observe that:
$$ \mathbb{P}(\Theta_k^{\varepsilon, \partial M} \leq \gamma \mathbb{E}[\Theta_k^{\varepsilon, \partial M}]) \leq \frac{\text{Var}(\Theta_k^{\varepsilon, \partial M})}{(1-\gamma)^2\mathbb{E}[\Theta_k^{\varepsilon, \partial M}]^2} $$

It suffices to show that the right-hand-side goes to zero. In bounding the variance we will be interested in the following term:

$$ \left(\Theta_k^{\varepsilon, \partial M}(r_1,r)\right)^2 = \sum_{\substack{\mathcal{Y}_1,\mathcal{Y}_2 \subset \mathcal{P}_n \\ |\mathcal{Y}_i| = k+1 }} g_{r}^{\varepsilon,\delta}(\mathcal{Y}_1,\mathcal{P}_n)g_{r}^{\varepsilon,\delta}(\mathcal{Y}_2,\mathcal{P}_n) $$

Let us separate the distinct situations in which the centres induced by $\mathcal{Y}_i$ are either close together $T_2$, or far apart $T_1$. Let us define:
\begin{align*} 
\Phi_r(\mathcal{Y}_1,\mathcal{Y}_2) &= \mathbf{1}\{B_r(c(\mathcal{Y}_1) \cap B_r(c(\mathcal{Y}_2) = \emptyset \} \\
T_1 =  \Theta_k^{\varepsilon, \partial M}(r_1,r)^2\Phi_{2r}(\mathcal{Y}_1,\mathcal{Y}_2) &= \sum_{\substack{\mathcal{Y}_1,\mathcal{Y}_2 \subset \mathcal{P}_n \\ |\mathcal{Y}_i| = k+1 }} g_{r}^{\varepsilon,\delta}(\mathcal{Y}_1,\mathcal{P}_n)g_{r}^{\varepsilon,\delta}(\mathcal{Y}_2,\mathcal{P}_n)\Phi_{2r}(\mathcal{Y}_1,\mathcal{Y}_2) \\
T_2 =  \Theta_k^{\varepsilon, \partial M}(r_1,r)^2(1-\Phi_{2r}(\mathcal{Y}_1,\mathcal{Y}_2) )&=\sum_{\substack{\mathcal{Y}_1,\mathcal{Y}_2 \subset \mathcal{P}_n \\ |\mathcal{Y}_i| = k+1 }} g_{r}^{\varepsilon,\delta}(\mathcal{Y}_1,\mathcal{P}_n)g_{r}^{\varepsilon,\delta}(\mathcal{Y}_2,\mathcal{P}_n)(1-\Phi_{2r}(\mathcal{Y}_1,\mathcal{Y}_2) )
\end{align*}

Thus we may express the variance as:

$$ \text{Var}(\Theta_k^{\varepsilon, \partial M}) = \mathbb{E}[(\Theta_k^{\varepsilon, \partial M})^2] - \mathbb{E}[\Theta_k^{\varepsilon, \partial M}]^2 = (\mathbb{E}[T_1] - \mathbb{E}[\Theta_k^{\varepsilon, \partial M}]^2 ) + \mathbb{E}[T_2] $$

Using Palm Theory \ref{thm:PalmTheory1} to simplify the sums taken over $(k+1)$-sized subsets we obtain:

$$ \mathbb{E}[\Theta_k^{\varepsilon, \partial M}]^2 = \frac{n^{2k+2}}{((k+1)!)^2} \mathbb{E}[g_{r}^{\varepsilon,\delta}(\mathcal{Y}_1',\mathcal{Y}_1' \cup \mathcal{P}_n)g_{r}^{\varepsilon,\delta}(\mathcal{Y}_2',\mathcal{Y}_2' \cup \mathcal{P}_n')]$$
$$ \mathbb{E}[T_1] = \frac{n^{2k+2}}{((k+1)!)^2} \mathbb{E}[g_{r}^{\varepsilon,\delta}(\mathcal{Y}_1',\mathcal{Y}' \cup \mathcal{P}_n)g_{r}^{\varepsilon,\delta}(\mathcal{Y}_2',\mathcal{Y}' \cup \mathcal{P}_n)]$$

Where $\mathcal{P}_n,\mathcal{P}_n'$ are i.i.d Poisson processes, $\mathcal{Y}_1', \mathcal{Y}_2'\subset M$ are i.i.d, and $\mathcal{Y}' = \mathcal{Y}_1' \cup \mathcal{Y}_2'$. Notice that when $\Phi_{2r}(\mathcal{Y}_1,\mathcal{Y}_2) \neq 0$ then $g_{r}^{\varepsilon,\delta}(\mathcal{Y}_i',\mathcal{Y}' \cup \mathcal{P}_n) = g_{r}^{\varepsilon,\delta}(\mathcal{Y}_i',\mathcal{Y}_i' \cup \mathcal{P}_n)$ since $g_{r}^{\varepsilon,\delta}$ is is only dependent on the points in the second argument which are sufficiently close to the points in the first argument. Thus we can make the following simple bound by omitting the negative contribution of the third expectation term in the first expression:

\begin{align*}
\mathbb{E}[T_1] - \mathbb{E}[\Theta_k^{\varepsilon, \partial M}]^2 &= \frac{n^{2k+2}}{((k+1)!)^2} \bigg( \mathbb{E}[g_{r}^{\varepsilon,\delta}(\mathcal{Y}_1',\mathcal{Y}' \cup \mathcal{P}_n)g_{r}^{\varepsilon,\delta}(\mathcal{Y}_2',\mathcal{Y}' \cup \mathcal{P}_n)\Phi_{2r}(\mathcal{Y}_1,\mathcal{Y}_2)] \\
& - \mathbb{E}[g_{r}^{\varepsilon,\delta}(\mathcal{Y}_1',\mathcal{Y}_1' \cup \mathcal{P}_n)g_{r}^{\varepsilon,\delta}(\mathcal{Y}_2',\mathcal{Y}_2' \cup \mathcal{P}_n')\Phi_{2r}(\mathcal{Y}_1,\mathcal{Y}_2)]\\
& - \mathbb{E}[g_{r}^{\varepsilon,\delta}(\mathcal{Y}_1',\mathcal{Y}_1' \cup \mathcal{P}_n)g_{r}^{\varepsilon,\delta}(\mathcal{Y}_2',\mathcal{Y}_2' \cup \mathcal{P}_n')(1-\Phi_{2r}(\mathcal{Y}_1,\mathcal{Y}_2))] \bigg) \\
& \leq \frac{n^{2k+2}}{((k+1)!)^2} \bigg( \mathbb{E}[g_{r}^{\varepsilon,\delta}(\mathcal{Y}_1',\mathcal{Y}_1' \cup \mathcal{P}_n)g_{r}^{\varepsilon,\delta}(\mathcal{Y}_2',\mathcal{Y}_2' \cup \mathcal{P}_n)\Phi_{2r}(\mathcal{Y}_1,\mathcal{Y}_2)] \\
& - \mathbb{E}[g_{r}^{\varepsilon,\delta}(\mathcal{Y}_1',\mathcal{Y}_1' \cup \mathcal{P}_n)g_{r}^{\varepsilon,\delta}(\mathcal{Y}_2',\mathcal{Y}_2' \cup \mathcal{P}_n')\Phi_{2r}(\mathcal{Y}_1,\mathcal{Y}_2)]\bigg) \\
& =: \frac{n^{2k+2}}{((k+1)!)^2} \mathbb{E}[\Delta g_{r}^{\varepsilon,\delta} ]
\end{align*}

Let us condition on the subsets $\mathcal{Y}_1',\mathcal{Y}_2'$ and note that whenever $\Delta g_{r}^{\varepsilon,\delta} \neq 0$ these two subsets are sufficiently separated so as not to interact. Thus by the independence property of Poisson processes:

\begin{align*}
\mathbb{E}[g_{r}^{\varepsilon,\delta}(\mathcal{Y}_1',\mathcal{Y}_1' \cup \mathcal{P}_n)g_{r}^{\varepsilon,\delta}(\mathcal{Y}_2',\mathcal{Y}_2' \cup \mathcal{P}_n) | \mathcal{Y}_1',\mathcal{Y}_2'] &=  \mathbb{E}[g_{r}^{\varepsilon,\delta}(\mathcal{Y}_1',\mathcal{Y}_1' \cup \mathcal{P}_n)|  \mathcal{Y}_1',\mathcal{Y}_2'] \mathbb{E}[g_{r}^{\varepsilon,\delta}(\mathcal{Y}_2',\mathcal{Y}_2' \cup \mathcal{P}_n) |  \mathcal{Y}_1',\mathcal{Y}_2'] \\ 
&= \mathbb{E}[g_{r}^{\varepsilon,\delta}(\mathcal{Y}_1',\mathcal{Y}_1' \cup \mathcal{P}_n)|  \mathcal{Y}_1',\mathcal{Y}_2'] \mathbb{E}[g_{r}^{\varepsilon,\delta}(\mathcal{Y}_2',\mathcal{Y}_2' \cup \mathcal{P}_n') |  \mathcal{Y}_1',\mathcal{Y}_2']
\end{align*}
\label{Second Moment Calculations}
Thus $\mathbb{E}[\Delta g_{r}^{\varepsilon,\delta} ] = \mathbb{E}[\mathbb{E}[\Delta g_{r}^{\varepsilon,\delta} |\mathcal{Y}_1',\mathcal{Y}_2']] = 0$. It remains to bound the second term of the variance $\mathbb{E}[T_2]$. Let us split the term $T_2$ into the separate cases in which $\mathcal{Y}_1, \mathcal{Y}_2$ share $j$ points:

$$ T_2 = \sum_{j=0}^{k+1} \sum_{|\mathcal{Y}_1\cap \mathcal{Y}_2| = j} g_{r}^{\varepsilon,\delta}(\mathcal{Y}_1,\mathcal{P}_n)g_{r}^{\varepsilon,\delta}(\mathcal{Y}_2,\mathcal{P}_n)(1-\Phi_{2r}(\mathcal{Y}_1,\mathcal{Y}_2) ) = \sum_{j=0}^{k+1} I_j$$


Using the Palm Theory result given in Corollary \ref{cor: PalmTheory2} yields:

\begin{align*}
\mathbb{E}[I_j] &= \sum_{|\mathcal{Y}_1\cap \mathcal{Y}_2| = j} g_{r}^{\varepsilon,\delta}(\mathcal{Y}_1,\mathcal{P}_n)g_{r}^{\varepsilon,\delta}(\mathcal{Y}_2,\mathcal{P}_n)(1-\Phi_{2r}(\mathcal{Y}_1,\mathcal{Y}_2) ) \\
& \leq  \frac{n^{2k+2-j}}{j!((k+1-j)!)^2} \int_{{DM}^{2k+2-j}} h_{r}^{\varepsilon,\delta}(\mathbf{y}_1)h_{r}^{\varepsilon,\delta}(\mathbf{y}_2)e^{-n\text{Vol}(\mathbf{y}_1, \mathbf{y}_2)} (1-\Phi_{2r}(\mathcal{Y}_1,\mathcal{Y}_2) )|\text{dvol}_g(\mathbf{y})|
\end{align*}

We bound the volume term using the results of Section \ref{subsec:RiemannianVolumes}: 
$$ \text{Vol}(\mathbf{y}_1, \mathbf{y}_2) = \text{Vol}(B_r(c(\mathbf{y}_1)\cup B_r(c(\mathbf{y}_2)) \geq (\frac{1}{2} + \alpha)(1 - (d \nu' r +\nu r^2))(1 + \frac{\omega_{d-1} d(c_1,c_2)}{\omega_{d}r} + O(\frac{d(c_1,c_2)^2}{r}))\omega_dr^d$$

Where $\alpha = O(\frac{\delta}{r}) \pm O(r)$.

We now make two separate change of variables separating the cases where $j=0$ and $j \neq 0$. We take the centre $c(\mathbf{y}_2)$ to be in polar coordinates around $c(\mathbf{y}_1)$. For details see Lemma \ref{lem: Multivariable BP Formula} and Section \ref{Blaschke} on Blaschke Petkantshin Formulae. Let us further partition our integral into two regions $\Omega_a, \Omega_b$:

$$\Omega_a = \{(\vec{y}_1,\vec{y}_2) : 0 \leq \frac{\rho(c(\vec{y}_1),c(\vec{y}_2))}{r} \leq \epsilon \}  $$
$$\Omega_b = \{(\vec{y}_1,\vec{y}_2) : \epsilon \leq \frac{\rho(c(\vec{y}_1),c(\vec{y}_2))}{r} \leq 4 \} $$
Under the assumption that $\Lambda r \to 0$ and $r_1 = (1 - \frac{1}{2c_g^2\Lambda^2})r$ we attain:

\begin{align*}
\mathbb{E}[I_0^{(a)}] & 
\leq \frac{n^{2k+2}}{((k+1)!)^2} \int_{\Omega_a} h_{r}^{\varepsilon,\delta}(\mathbf{y}_1)h_{r}^{\varepsilon,\delta}(\mathbf{y}_2)e^{-n\text{Vol}(\mathbf{y}_1, \mathbf{y}_2))} (1-\Phi_{2r}(\mathcal{Y}_1,\mathcal{Y}_2) )|\text{dvol}_g(\mathbf{y})| \\
& 
\leq C n^{2k+2} \int_{\Omega_a} h_{r}^{\varepsilon,\delta}(\mathbf{y}_1)h_{r}^{\varepsilon,\delta}(\mathbf{y}_2)e^{-\frac{\Lambda}{2}(1-(d\eta'r + \eta r^2))} (1-\Phi_{2r}(\mathcal{Y}_1,\mathcal{Y}_2) )|\text{dvol}_g(\mathbf{y})| \\
& 
\leq C n^{2k+2} e^{-\frac{\Lambda}{2}} \int_{\partial M_{[\delta,2\delta]}} |\text{dvol}_g(c_1)|
\int_{0}^{\epsilon r} ds \  \int_{\mathbb{S}_1(T_{c_1}M)} s^{d-1} \text{dvol}_{\mathbb{S}_1(T_{c_1}M)}(w)  \\
& 
\times  \prod_{i=1}^{2} \int_{r_1}^r du_i \int_{Gr(k,d)} u_i^{k(d-k)} d\mu_{k,d}(V) \int_{(\mathbb{S}_1(V))^{k+1}} u_i^{(k-1)(k+1)} |\text{dvol}_{(\mathbb{S}_1(V))^{k+1}}(\mathbf{w}_i) |
h_{r}^{\varepsilon,\delta}(\text{exp}_{c_i}(u_1\mathbf{w}_i)) \\
& 
\leq C n^{2k+2} e^{-\frac{\Lambda}{2}} \delta \epsilon^d r^{d(2k+1)} \Lambda^{-4} \leq C  n e^{-\frac{\Lambda}{2}} \Lambda^{2k-3} \epsilon^d r
\end{align*}


\begin{align*}
\mathbb{E}[I_0^{(b)}] & 
\leq \frac{n^{2k+2}}{((k+1)!)^2} \int_{\Omega_b} h_{r}^{\varepsilon,\delta}(\mathbf{y}_1)h_{r}^{\varepsilon,\delta}(\mathbf{y}_2)e^{-n\text{Vol}(\mathbf{y}_1, \mathbf{y}_2))} (1-\Phi_{2r}(\mathcal{Y}_1,\mathcal{Y}_2) )|\text{dvol}_g(\mathbf{y})| \\
& 
\leq C n^{2k+2} \int_{\Omega_b} h_{r}^{\varepsilon,\delta}(\mathbf{y}_1)h_{r}^{\varepsilon,\delta}(\mathbf{y}_2)e^{-\frac{\Lambda}{2}(1+ \frac{\epsilon \omega_d}{\omega_{d-1}})(1-(d\eta'r + \eta r^2))} (1-\Phi_{2r}(\mathcal{Y}_1,\mathcal{Y}_2) )|\text{dvol}_g(\mathbf{y})| \\
& 
\leq C n^{2k+2} e^{-\frac{\Lambda}{2}(1+ \frac{\epsilon \omega_d}{\omega_{d-1}})} \int_{\partial M_{[\delta,2\delta]}} |\text{dvol}_g(c_1)|
\int_{\epsilon r}^{4r} ds \  \int_{\mathbb{S}_1(T_{c_1}M)} s^{d-1} \text{dvol}_{\mathbb{S}_1(T_{c_1}M)}(w)  \\
& 
\times  \prod_{i=1}^{2} \int_{r_1}^r du_i \int_{Gr(k,d)} u_i^{k(d-k)} d\mu_{k,d}(V) \int_{(\mathbb{S}_1(V))^{k+1}} u_i^{(k-1)(k+1)} |\text{dvol}_{(\mathbb{S}_1(V))^{k+1}}(\mathbf{w}_i) |
h_{r}^{\varepsilon,\delta}(\text{exp}_{c_i}(u_1\mathbf{w}_i)) \\
& 
\leq C n^{2k+2} e^{-\frac{\Lambda}{2}(1+ \frac{\epsilon \omega_d}{\omega_{d-1}})} \delta  r^{d(2k+1)} \Lambda^{-4} \leq C  n e^{-\frac{\Lambda}{2}} e^{-\frac{\Lambda \epsilon \omega_d}{2 \omega_{d-1}}} \Lambda^{2k-3} r
\end{align*}

Where $c_2 = \text{exp}_{c_1}(sw)$, $\mathbf{y}_i = \text{exp}_{c_i}(\mathbf{w}_i)$ and $\mathbb{S}_1(V)$ is the unit sphere.

Similarly for $j\neq 0$ we attain bounds:

\begin{align*}
\mathbb{E}[I_j^{(a)}] & 
\leq \frac{n^{2k+2-j}}{j!((k+1-j)!)^2} \int_{\Omega_a} h_{r}^{\varepsilon,\delta}(\mathbf{y}_1)h_{r}^{\varepsilon,\delta}(\mathbf{y}_2)e^{-n\text{Vol}(\mathbf{y}_1, \mathbf{y}_2))} (1-\Phi_{2r}(\mathcal{Y}_1,\mathcal{Y}_2) )|\text{dvol}_g(\mathbf{y})| \\
& 
\leq n^{2k+2-j} \int_{\Omega_a} h_{r}^{\varepsilon,\delta}(\mathbf{y}_1)h_{r}^{\varepsilon,\delta}(\mathbf{y}_2)e^{-\frac{\Lambda}{2}(1-(d\eta'r + \eta r^2))} (1-\Phi_{2r}(\mathcal{Y}_1,\mathcal{Y}_2) )|\text{dvol}_g(\mathbf{y})| \\
& 
\leq C n^{2k+2-j} e^{-\frac{\Lambda}{2}} \int_{\partial M_{[\delta,2\delta]}} |\text{dvol}_g(c_1)| \\
& 
\times \int_{r_1}^r du_1 \ \int_{Gr(k,d)} u_1^{k(d-k)} d\mu_{k,d}(V_1) \int_{(\mathbb{S}_1(V_1))^{k+1}} u_i^{(k-1)(k+1)} \text{dvol}_{(\mathbb{S}_1(V_1))^{k+1}} \\
& 
\times \int_{0}^{\epsilon r} ds \  \int_{\mathbb{S}_1(T_{c_1}E)} s^{d-j}  \text{dvol}_{\mathbb{S}^{d-1}}(w) \int_{r_1}^r du_2  \int_{Gr(k-j,d)} u_i^{(k-j)(d-(k-j))} d\mu_{k-j,d}(W)\\
& 
\times  \int_{(\mathbb{S}_1(V_2))^{k+1-j}} u_2^{(k-1)(k+1-j)} \text{dvol}_{(\mathbb{S}_1(V_2))^{k+1-j}}
h_{r}^{\varepsilon,\delta}(\text{exp}_{c_1}(u_1\mathbf{w}_1))h_{r}^{\varepsilon,\delta}(\text{exp}_{c_2}(u_2\mathbf{w}_2)) \\
& 
\leq C n^{2k+2-j} e^{-\frac{\Lambda}{2}} \delta \epsilon^{d-j+1} r^{d(2k+1-j)+j(k-j)+1} \Lambda^{-4} \leq C n e^{-\frac{\Lambda}{2}} \Lambda^{2k-3} \epsilon^{d-j+1} r^{j(k-j)+1}
\end{align*}

\begin{align*}
\mathbb{E}[I_j^{(b)}] & 
\leq \frac{n^{2k+2-j}}{j!((k+1-j)!)^2} \int_{\Omega_b} h_{r}^{\varepsilon,\delta}(\mathbf{y}_1)h_{r}^{\varepsilon,\delta}(\mathbf{y}_2)e^{-n\text{Vol}(\mathbf{y}_1, \mathbf{y}_2))} (1-\Phi_{2r}(\mathcal{Y}_1,\mathcal{Y}_2) )|\text{dvol}_g(\mathbf{y})| \\
& 
\leq n^{2k+2-j} \int_{\Omega_b} h_{r}^{\varepsilon,\delta}(\mathbf{y}_1)h_{r}^{\varepsilon,\delta}(\mathbf{y}_2)e^{-\frac{\Lambda}{2}(1+ \frac{\epsilon \omega_d}{\omega_{d-1}})(1-(d\eta'r + \eta r^2))} (1-\Phi_{2r}(\mathcal{Y}_1,\mathcal{Y}_2) )|\text{dvol}_g(\mathbf{y})| \\
& 
\leq C n^{2k+2-j} e^{-\frac{\Lambda}{2}(1+ \frac{\epsilon \omega_d}{\omega_{d-1}})} \int_{\partial M_{[\delta,2\delta]}} |\text{dvol}_g(c_1)| \\
& 
\times \int_{r_1}^r du_1 \ \int_{Gr(k,d)} u_1^{k(d-k)} d\mu_{k,d}(V_1) \int_{(\mathbb{S}_1(V_1))^{k+1}} u_i^{(k-1)(k+1)} \text{dvol}_{(\mathbb{S}_1(V_1))^{k+1}} \\
& 
\times \int_{\epsilon r}^{4 r} ds \  \int_{\mathbb{S}_1(T_{c_1}E)} s^{d-j}  \text{dvol}_{\mathbb{S}^{d-1}}(w) \int_{r_1}^r du_2  \int_{Gr(k-j,d)} u_i^{(k-j)(d-(k-j))} d\mu_{k-j,d}(W)\\
& 
\times  \int_{(\mathbb{S}_1(V_2))^{k+1-j}} u_2^{(k-1)(k+1-j)} \text{dvol}_{(\mathbb{S}_1(V_2))^{k+1-j}}
h_{r}^{\varepsilon,\delta}(\text{exp}_{c_1}(u_1\mathbf{w}_1))h_{r}^{\varepsilon,\delta}(\text{exp}_{c_2}(u_2\mathbf{w}_2)) \\
& 
\leq C n^{2k+2-j} e^{-\frac{\Lambda}{2}(1+ \frac{\epsilon \omega_d}{\omega_{d-1}})} \delta r^{d(2k+1-j)+j(k-j)+1} \Lambda^{-4} \leq C n e^{-\frac{\Lambda}{2}} e^{-\frac{\Lambda \epsilon \omega_d}{2 \omega_{d-1}}} \Lambda^{2k-3} r^{j(k-j)+1}
\end{align*}

We note that our largest bounding function is for the expectation $\mathbb{E}[I_0]$ and so we have the following bound for all $j$ and for arbitrary $\epsilon \in (0,1)$:

$$ \mathbb{E}[I_j] \leq C n  e^{-\frac{\Lambda}{2}} \Lambda^{2k-3}r(e^{-\frac{\Lambda \epsilon \omega_d}{2 \omega_{d-1}}} + \epsilon^d) $$

We may therefore use this same regime to bound $\mathbb{E}[T_2]$. Let us choose $\epsilon = \frac{2(2k-1)\omega_{d-1}}{\omega_d} \frac{\log\log n}{\log n}$.
Then for $\Lambda = (2-\frac{2}{d})\log n + 2(k-2 -(k+1- \frac{1}{d})) \log \log n - \omega (n)$, (recalling that $\mathbb{E}[\Theta_k^{\varepsilon, \partial M}] = \Omega(e^{-\alpha\Lambda} n\Lambda^{k-2}r(\log n)^{-(k+1)})$ and that $\alpha = \frac{1}{2} + O((\log n )^{-1})$), we calculate:
\begin{align*}
 \frac{\text{Var}(\Theta_k^{\varepsilon, \partial M})}{\mathbb{E}[\Theta_k^{\varepsilon, \partial M}]^2} &=  \frac{\mathbb{E}[T_2]}{\mathbb{E}[\Theta_k^{\varepsilon, \partial M}]^2} 
 \leq C \frac{n  e^{-\frac{\Lambda}{2}} \Lambda^{2k-3}r (e^{-\frac{\Lambda \epsilon \omega_d}{2 \omega_{d-1}}} + \epsilon^d)}{n^2 \Lambda^{2k-4}e^{-2\alpha\Lambda}r^{2}(\log n)^{-2(k+1)}} 
 \sim  \frac{\Lambda e^{\frac{\Lambda}{2}}(\log n)^{2(k+1)}(e^{-\frac{\Lambda \epsilon \omega_d}{2 \omega_{d-1}}} + \epsilon^d)}{nr} \\ 
 &\sim e^{-\frac{\omega(n)}{2}} (\log n)^{2k}(e^{-\frac{\Lambda \epsilon \omega_d}{2 \omega_{d-1}}} + \epsilon^d) 
 \sim e^{-\frac{\omega(n)}{2}} \left( \frac{1}{\log n} + \frac{(\log \log n)^d }{(\log n)^d}\right)\to 0 
\end{align*}

Thus we conclude that w.h.p. $\Theta_k^{\varepsilon, \partial M} \geq \frac{1}{2} \mathbb{E}[\Theta_k^{\varepsilon, \partial M}]$.

%% file: conclusion.tex
\label{Conclusion}
Collecting the results contained in the previous sections we attain thresholds in terms of $\Lambda$ for the $k^\text{th}$ homological connectivity of a compact manifold with boundary $M$. For $\Lambda$ greater than the upper threshold we recover the $k^\text{th}$ homology of $M$ with high probability, and for $\Lambda$ less than the lower threshold we do not recover the  $k^\text{th}$ homology of $M$ with high probability.


A recent coverage result from \cite{Wei18} 
yields a sharp coverage threshold for Riemannian manifolds with boundary. For $\Lambda = (2-\frac{2}{d})\log n + 2(d-2) \log \log n + w(n)$ the manifold is covered with high probability. Our result shows that the homological connectivity threshold for lower homology groups occurs before coverage.


\begin{thm}(Homological Connectivity of Riemannian Manifold with Boundary) 

Let $M$ be a unit volume compact Riemannian manifold with smooth non-empty boundary. Let $d\geq 2$ be the dimension of $M$, $\Lambda = \omega_d n r^d$ and $\mathcal{P}_n$ a Poisson process of intensity $n$ on $M$. Then for $1 \leq k \leq d-1$ 

$$ \lim_{n\to \infty} \mathbb{P}(H_k(\mathcal{C}(n,r)) \cong H_k(M)) = \begin{cases}
1 & \Lambda = (2-\frac{2}{d})\log n + 2k \log \log n + w (n),\\
0 & \Lambda = (2-\frac{2}{d})\log n + 2(k - 2 -(k+1- \frac{1}{d})) \log \log n - w (n), 
\end{cases}$$
\label{thm:Homological Connectivity Boundary}
\end{thm}
\begin{proof}(Upper Threshold)

Let $r_0 = r(\frac{\omega_d}{\kappa}(1+|\log r |))^{1/d}$ so that the conditions of Lemma \ref{Critical point upper bound} are met.

Asymptotically $\Lambda_{\frac{r_0}{2}} \geq 2\log n $. Using our Asymptotic Coverage result Theorem \ref{Coverage} we observe that $M \subset B_{r_0}(\mathcal{P}_n)$ w.h.p. and so for sufficiently small $r_0$, by the Nerve Lemma $(H_k(\mathcal{C}(n,r_0)) \cong H_k(M))$ w.h.p.

Moreover for $\Lambda$ as described in the upper threshold we observe that $n \Lambda^k e^{-\Lambda}, n^{1-\frac{1}{d}}\Lambda^{k}e^{-\frac{1}{2}\Lambda} \to 0$ so by Lemma \ref{Critical point upper bound} the expected number of $k$-critical and $(k+1)$-critical points with critical value in the range $(r,r_0]$ goes to zero: 

\begin{align*}
 \Lambda = (2-\frac{2}{d})\log n + 2k \log \log n + w (n), \implies n \Lambda^{k} e^{-\Lambda}= O\left(\frac{n (\log n)^k}{n^{(2-2/d)}(\log n)^{2k}e^{w (n)} }\right) \to 0
\end{align*}
\begin{align*}
 \Lambda = (2-\frac{2}{d})\log n + 2k \log \log n + w (n), \implies n^{1-\frac{1}{d}}\Lambda^{k}e^{-\frac{1}{2}\Lambda}&= O\left(\frac{n^{1-\frac{1}{d}} (\log n)^k}{n^{1-1/d}(\log n)^{k} e^{w(n)} }\right) \to 0
\end{align*}

Using Morse Theory for manifolds with boundary we yield that $H_k(\mathcal{C}(n,r))\cong H_k(\mathcal{C}(n,r_0))$ w.h.p. since the probability of there being a critical point in the range $(r,r_0]$ tends to zero, and thus  $H_k(\mathcal{C}(n,r))\cong H_k(M)$ w.h.p. for $\Lambda = (2-\frac{2}{d})\log n + 2k \log \log n + w (n)$

\end{proof}

\begin{proof}(Lower Threshold)

Recall that by Lemma \ref{Critical point lower bound} we have $\mathbb{E}[\Theta_k^{\varepsilon, \partial M}] = \Omega(e^{-\alpha\Lambda} n\Lambda^{k-2}r(\log n)^{-(k+1)})$ for $\alpha = \frac{1}{2} + O((\log n )^{-1})$. The Second Moment Calculations verify that w.h.p. $|\beta_k(r)| = \Omega(e^{-\alpha\Lambda} n\Lambda^{k-2}r(\log n)^{-(k+1)})$.

Then we note that for $\Lambda = (2-\frac{2}{d})\log n + 2(k - 2 - (k+1 -\frac{1}{d}) ) \log \log n - w (n)$:
\begin{align*}
e^{-\alpha\Lambda} n\Lambda^{k-2}r(\log n)^{-(k+1)} \sim \frac{n^{1-\frac{1}{d}} (\log n)^{k-2+\frac{1}{d}-(k+1)}}{e^{\frac{\Lambda}{2}}} \sim \frac{n^{1-\frac{1}{d}} (\log n)^{k-2+\frac{1}{d}-(k+1)}}{n^{1-\frac{1}{d}} (\log n)^{k-2 +\frac{1}{d}-(k+1)} e^{-\frac{w(n)}{2}}} \to \infty
\end{align*}

Thus for $\Lambda$ in this regime the Betti numbers tend to infinity w.h.p. and so we do not recover the homology of $M$.

\end{proof}

Having adapted the techniques from \cite{Bobrowski2017} in order to count critical points near to the boundary we have attained thresholds similar to those of the Homological Connectivity Theorem from \cite{Bobrowski2017}, applicable to compact, closed Riemannian Manifolds. 

\begin{thm}\cite{Bobrowski2017}(Homological Connectivity Thresholds for Compact Manifolds without Boundary)

Let $M$ be a unit volume compact, Riemannian manifold without boundary.
Suppose that as $n \to \infty$, $w(n) \to \infty$. Then for $1\leq k\leq d-1$
$$ \lim_{n\to \infty} \mathbb{P}(H_k(\mathcal{C}(n,r)) \cong H_k(M)) = \begin{cases}
1 & \Lambda = \log n + k \log \log n + w (n),\\
0 & \Lambda = \log n + (k-2) \log \log n - w (n),
\end{cases}$$

\label{thm:Homological Connectivity Closed}
\end{thm}

The first point of similarity to note between Theorem \ref{thm:Homological Connectivity Boundary} and Theorem \ref{thm:Homological Connectivity Closed} is that whilst neither Theorem identifies a sharp threshold, both identify the first order term for the transition to homological connectivity: $(2-\frac{2}{d})\log n$ and $\log n $ respectively. It is worth noting how the geometric differences between building a \v{C}ech complex on a manifold with boundary, rather than a closed manifold, inform the differences in the coefficients of the terms in the homological connectivity thresholds.

Our analysis shows that the distance function of a Poisson point process yields a large number of critical points near to the boundary, and this results in the disparity between the constant factor of the first order terms. For a collection of points $\mathcal{Y}$ of the Poisson process to induce a critical point at the centre of these points $c(\mathcal{Y})$, one requires that no other point of the Poisson point process lies in the ball of radius $\rho(\mathcal{Y})$ centred at $c(\mathcal{Y})$. The existence of critical points near to the boundary is made more likely by the fact that the ball of radius $\rho(\mathcal{Y})$ centred at $c(\mathcal{Y})$ is cut by the boundary, and so it is more likely that no other point of the Poisson point process lies in this cut ball. In the most extreme case the ball's volume is cut in half by the boundary and this introduces the factor of $2$ in the leading term of the threshold. Since this phenomenon occurs only for critical points near to the boundary we have to scale our count by the volume of an $r$-collar neighbourhood of the boundary. This volume behaves like $\text{Vol}(\partial M_r) \sim r \sim (\frac{\log n}{n})^{\frac{1}{d}}$, and in particular the factor $n^{-\frac{1}{d}}$ introduces the term $-\frac{1}{d}$ to the leading coefficient in Theorem \ref{thm:Homological Connectivity Boundary}.

The second order disparity is identified to lie in the range $[(k-2)\log \log n, k \log \log n]$ by Theorem \ref{thm:Homological Connectivity Closed} and in the range $[2(\frac{1}{d}-3)\log \log n,2k\log \log n]$ by Theorem \ref{thm:Homological Connectivity Boundary}. The factor of $2$ in the second order terms is again induced by the effect of the boundary cutting volumes in half. The second order term of the lower threshold is affected by counting special erroneous cycles which occur near the boundary which we call $\Theta$-like-cycles. Since we only count points near the boundary, again a $\frac{1}{d}$ term is introduced to the coefficient of the second order term. The coefficient of the second order term is further impacted by a lower bound for the volume of a subset of a Grassmannian, and may be able to be improved if this bound is sharpened.


\vskip .2in
\noindent
{\bf Relative homology.} Here we studied the homological connectivity of a compact Riemannian manifold with boundary in terms of the absolute homology. Equally we could have chosen to study the homology of the manifold relative to its boundary. We could perform the same analysis to attain connectivity thresholds for the relative homology by counting critical points of a Morse function induced by the Point process. Recall that the Morse complex of the Morse function recovers absolute homology if our Morse function attains a maximum on the boundary, and relative homology if our Morse function attains a minimum on the boundary (Theorem \ref{thm: Morse Theory}). Thus the negative of the distance function from the point process would induce a Morse complex which calculates the relative homology. 

Taking the negative of the distance function converts index $k$ critical points into index $d - k$ critical points. With this setup we lose the geometric interpretation that the Morse complex at scale $r$ computes the homology of the union of radius $r$ balls built around the point process. This correspondence is a result analogous to Lefschetz Duality combined with the Universal Coefficient Theorem for Cohomology for calculating the relative homology, which makes it clear that given $H_k(M)$ and $H_{k-1}(M)$ we can compute $H_{d-k}(M,\partial M)$.

There is a geometric interpretation for the dual to our $\Theta$-like-cycles when we consider the Morse complex of the negative of the distance function. For a $\Theta$-like-cycle of index $k$ we require a $k$-simplex to be introduced approximately tangential to the boundary. The corresponding orthogonal dual $(d-k)$-simplex crosses the part of the partial annulus cut by the boundary and introduces a new homological cycle when we take homology relative to the boundary. See Figure \ref{fig:Non-Theta-Cycle}, for an illustration of such a dual simplex which introduces a relative cycle.


%% file: list-of-symbols.tex
\begin{tabular}{cp{0.9\textwidth}}
  $M$ & Compact Riemannian Manifold with boundary \\
  $g$ & Smooth Riemannian metric \\
  $\rho(\cdot,\cdot)$ & Distance induced by Riemannian Metric \\
  $\tau_M$ & The reach of the manifold $M$\\
  $\partial M_r $ & An $r$ neighbourhood of the boundary \\
  $P$ & A finite sample of points of a manifold \\
  $\mathcal{P},\mathcal{Q}$ & A Poisson point process \\ 
  $B_r(P)$ & The union of radius $r$ balls centred at each $p \in P$\\
  $ \mathcal{C}(n,r)$ & \v{C}ech complex at radius $r$ on a point process of intensity $n$\\
  $c(\mathcal{Y}), c(\vec{y})$ & The centre of a finite collection of points \\
  $\rho(\mathcal{Y}), \rho(\vec{y})$ & The critical value of a finite collection of points \\
  $B(\mathcal{Y}), B(\vec{y})$ & The ball centred at $c(\mathcal{Y})$ with radius $\rho(\mathcal{Y})$ \\
  $\beta_k(r)$ & The $k$-th Betti number associated to a \v{C}ech complex at radius $r$ \\
  $\omega_d$ & Volume of a unit radius $d$-dimensional ball \\
  $ \Lambda$ & Expected number of points lying in an $r$-ball $\Lambda = \omega_d n r^d$ \\
  $C$ & Constant factor, a product of constant terms used to simplify expressions in inequalities \\
  $ C^{\rho_M}_k(r,r_0)$ & Set of index $k$ critical points with critical values in the range $[r,r_0)$ \\
  $ \Theta^{\varepsilon}_k(r_1,r)$ & Set of index $k$ critical points with critical values in the range $[r_1,r)$ inducing $\Theta$-cycles \\
  $ \Theta^{\varepsilon, \partial M}_k(r_1,r)$ & Set of index $k$ critical points with critical values in the range $[r_1,r)$ inducing $\Theta$-like-cycles \\
  $ A_\varepsilon(c)$ & $\varepsilon$-annulus about centre $c(\vec{y})$ of radius $\rho(\vec{y})$ \\
  $ A^{(\varphi)}_\varepsilon(c)$ & Partial $\varepsilon$-annulus about centre $c(\vec{y})$ of radius $\rho(\vec{y})$\\
\end{tabular}\\